\RequirePackage{silence}
\documentclass[a4paper,svgnames]{amsart}
\usepackage[hmarginratio={1:1},vmarginratio={1:1},lmargin=65.0pt,tmargin=60.0pt]{geometry}

\allowdisplaybreaks

\usepackage[T1]{fontenc}
\usepackage{soul,latexsym,exscale,mathtools,textcomp,bbm,aliascnt,atableaux}
\usepackage{amsmath,amsthm,amsfonts,amssymb,amscd,verbatim}
\usepackage{mathrsfs,stackrel}
\usepackage{tikz,tikz-cd}
\usetikzlibrary{matrix,arrows.meta,decorations.markings,shapes.multipart,decorations.pathreplacing,decorations.pathmorphing}
\usepackage{enumitem}
\setlist[enumerate]{itemsep=0.15cm,label=\emph{\upshape(\alph*)}}
\setlist[enumerate,2]{itemsep=0.15cm,label=\emph{\upshape(\roman*)}}
\usepackage{ytableau}
\usepackage{xparse}
\usepackage{cite}

\usepackage{array}
\newcolumntype{C}{>{$}c<{$}}
\usepackage{xltabular,booktabs}

\usepackage{mfabacus}

\usepackage{tabulary}
\usepackage{booktabs}

\definecolor{mygray}{gray}{0.6}
\definecolor{mygraydark}{gray}{0.4}
\definecolor{mygraylight}{gray}{0.85}
\definecolor{spinach}{RGB}{46,139,87}
\definecolor{tomato}{RGB}{255,99,71}
\definecolor{orchid}{RGB}{143,40,194}
\definecolor{neon}{RGB}{77,77,255}
\definecolor{pumpkin}{RGB}{224,180,80}
\definecolor{citron}{RGB}{190,180,90}

\definecolor{lava}{RGB}{207,16,32}
\definecolor{cream}{RGB}{255,253,208}
\definecolor{verdigris}{RGB}{67,179,174}
\definecolor{Black}{RGB}{0,0,0}
\definecolor{mydarkblue}{RGB}{10,10,170}
\definecolor{darkspinach}{RGB}{20,70,20}
\definecolor{darktomato}{RGB}{155,40,30}
\definecolor{darkorchid}{RGB}{50,10,100}
\definecolor{darklava}{RGB}{150,8,16}

\let\emph\relax
\DeclareTextFontCommand{\emph}{\bfseries\em}

\newcommand{\eg}{\text{e.g.}}
\newcommand{\cf}{\text{cf.}}

\newcommand{\vive}{\text{vice versa}}

\usepackage{todonotes}

\newcommand\floor[2]{\lfloor\tfrac{#1}{#2}\rfloor}

\newcommand\bsig{{\boldsymbol\sigma}}

\newcommand\blam{{\boldsymbol\lambda}}

\newcommand\brho{{\boldsymbol\rho}}
\newcommand\bmu{{\boldsymbol\mu}}
\newcommand\bnu{{\boldsymbol\nu}}

\newcommand\charge{{\boldsymbol\kappa}}
\newcommand\chargetwo{{\boldsymbol\nu}}
\newcommand\chargethree{{\boldsymbol\upsilon}}

\newcommand\ba{{\boldsymbol{a}}}
\newcommand\bb{{\boldsymbol{b}}}

\newcommand\bx{\boldsymbol{x}}
\newcommand\by{\boldsymbol{y}}
\newcommand\bz{\boldsymbol{z}}
\newcommand\bS{{\boldsymbol{S}}}
\newcommand\bT{{\boldsymbol{T}}}
\newcommand\bU{{\boldsymbol{U}}}
\newcommand\bV{{\boldsymbol{V}}}

\newcommand\Pcal{\mathcal{P}}

\newcommand\str[1][s]{\mathsf{#1}}
\newcommand\Sub[1]{\overline{#1}}
\newcommand\Submap{\overline{S}}

\newcommand\Web{\mathbb{W}}

\NewDocumentCommand\Webaa{ O{\bx,\bi} O{\by,\bj} }{\Web_{#1}^{#2}}
\NewDocumentCommand\Webab{ O{(\charge,\bx),\bi} O{(\chargetwo,\by),\bj} }{\Web_{#1}^{#2}}
\NewDocumentCommand\Webabs{ O{\bx} O{\by} O{\bi} O{\bj} }{\Web_{\Sub{#1},\Sub{#3}}^{\Sub{#2},\Sub{#4}}}

\NewDocumentCommand\WA{ O{\beta} O{\brho} }{\mathscr{W}_{#1}^{#2}}
\NewDocumentCommand\WAc{ O{\beta} O{\brho} }{\mathscr{R}_{#1}^{#2}}
\NewDocumentCommand\WAs{ O{\beta} O{\brho} }{\mathscr{W}_{\Sub{#1}}^{{#2}}}
\NewDocumentCommand\WAsc{ O{\beta} O{\brho} }{\mathscr{R}^{\Sub{#1}}^{\Sub{#2}}}
\newcommand\WABasis{\mathcal{B}_{\beta}}
\NewDocumentCommand\SA{ O{\beta} O{\brho} }{\overline{\mathscr{W}}_{#1}^{#2}}

\NewDocumentCommand\TA{ O{\beta} O{\brho} }{\mathcal{W}_{#1}^{#2}}
\NewDocumentCommand\TAc{ O{\beta} O{\brho} }{\mathcal{R}_{#1}^{#2}}

\NewDocumentCommand\WAlam{ s O{\blam} } {\mathscr{W}_{n}^{\IfBooleanTF{#1}{\gdom}{\gedom}#2}}

\newcommand\sand[1][\blam]{\mathbb{S}_{#1}}

\newcommand\Oneg{\1_{\Sub{\Gamma},\Gamma}}

\newcommand\affine[1]{\underline{#1}}
\newcommand\Sym[1][n]{{\mathfrak{S}_{#1}}}
\newcommand\hell{\affine{\ell}}
\NewDocumentCommand\Parts{ D(){\ell} O{n} }{\mathrm{P}_{#1,#2}^{\brho}}
\newcommand\hParts[1][n]{\affine{\mathrm{P}}_{\hell,#1}^{\affine{\brho}}}

\newcommand\Nodes[1][n]{\mathscr{N}_{\ell,#1}}
\newcommand\hcoord{\mathtt{x}_{\affine{\charge}}}

\NewDocumentCommand\res{sO{}}{\mathop{\rm res}\nolimits_{#2\rho}}
\newcommand\Affch[1][\blam]{A^{y}(#1)}

\newcommand\BX[1][{\WA[n](X)}]{\mathrm{D}_{#1}}
\newcommand\BXc[1][{\WAc[n](X)}]{\mathrm{D}_{#1}}
\newcommand\EX[1][{\TA[n](X)}]{\mathrm{E}_{#1}}
\newcommand\EXc[1][{\TAc[n](X)}]{\mathrm{E}_{#1}}

\def\Item(#1){\item\textbf{\upshape(#1\upshape)}}
\newcounter{relation}

\def\relationautorefname~#1\null{\upshape(W$_{#1}$\upshape)}
\let\ref\autoref
\let\eqref\autoref

\newcommand\SStd[1][\charge]{\mathop{\rm SStd}\nolimits_{#1}}
\newcommand\hSStd[1][\underline{\charge}]{\mathop{\rm SStd}\nolimits_{#1}}

\newcommand\N{\mathbb{Z}_{\geq 0}}

\newcommand\Z{\mathbb{Z}}

\def\y{\mathbf{y}}

\let\<=\langle
\let\>=\rangle

\NewDocumentCommand\mI{ sd() }{\mathtt{m}_{\IfBooleanT{#1}{\Sub}\Gamma}\IfNoValueF{#2}{(#2)}}
\def\pmod#1{\space(\text{mod }#1)}

\newcommand\bi{\mathbf{i}}
\newcommand\bj{\mathbf{j}}
\newcommand\bk{\mathbf{k}}

\newcommand\R{\mathbb{R}}

\def\map#1#2{\,{:}\,#1\!\longrightarrow\!#2}

\let\gedom=\trianglerighteq
\let\gdom=\vartriangleright
\let\ledom=\trianglelefteq
\let\ldom=\vartriangleleft

\usepackage{etoolbox}
\newcommand{\DeclareMyOperator}[1]{%
\expandafter\DeclareMathOperator\csname #1\endcsname{#1}
}
\forcsvlist{\DeclareMyOperator}{Shape,defect,gr,End}

\def\1{\mathbf{1}}
\def\height{\mathop{\rm ht}\nolimits}

\def\Item(#1){\item\textbf{\upshape(#1\upshape)}}

\DeclareMathOperator\Ab{Ab}
\tikzset{
anchorbase/.style={baseline={([yshift=#1]current bounding box.center)}},
anchorbase/.default={-0.5ex},
dot colour/.initial=black,
dot colour/.default=black,
tinynodes/.style={font=\tiny,text height=0.25ex,text depth=0.05ex},
smallnodes/.style={font=\scriptsize,text height=0.75ex,text depth=0.15ex},
mor/.style={line width=0.75,color=black,fill=cream},
dots/.style={line width=1pt,line cap=round, gray, dash pattern=on 0pt off 2\pgflinewidth},
redstring/.style = {draw=red!50,fill=none,line width=0.35mm,preaction={draw=red,line width=2.5pt,-},nodes={color=red}},
affine/.style= {draw=citron!50,fill=none,
line width=0.35mm,preaction={draw=citron,line width=2.5pt,-},nodes={color=citron}},
solid/.style = {draw=blue,fill=none,dot colour=blue,line width=0.4mm,nodes={color=blue}},
ghost/.style = {draw=darkgray,fill=none,dot colour=darkgray,
densely dashed,line width=0.4mm,nodes={color=darkgray}},
dghost/.style = {draw=darkgray,double,fill=none,dot colour=darkgray,
densely dashed,line width=0.4mm,nodes={color=darkgray}},
crossline/.style={preaction={draw=white,line width=4.75pt,-},preaction={draw=black,line width=0.9pt,-}},
dot/.style = {
decoration={markings,
post length=0.25mm,
pre length=0.25mm,
mark=at position #1 with {\node[circle,radius=0.3cm,inner sep=-2.0pt,color=\pgfkeysvalueof{/tikz/dot colour},fill=\pgfkeysvalueof{/tikz/dot colour}]{};}
},
postaction={decorate}
},
dot/.default=0.5,
}

\tikzstyle directed=[postaction={decorate,decoration={markings,
mark=at position #1 with {\arrow[line width=0.25mm, black]{>}}}}]

\NewDocumentCommand\crossing{ O{1ex} mmO{} mmO{} } {
\tikz[centered=#1]{
\draw[#2, dot/.list={#4}](0,0)node[below]{$#3$}--++(1,1);
\draw[#5, dot/.list={#7}](1,0)node[below]{$#6$}--++(-1,1);
}%
}

\NewDocumentCommand\DottedIdempotentC{ D(){1.5} omm} {
\begin{tikzpicture}[anchorbase]
\def\residues{{#4}}
\foreach \res [count=\c,
evaluate=\res as \pos using {\res<=#3 ? #1*\res-\c*#1/15 : #1*(\res-2)-\c*#1/15}
] in {#4} {
\coordinate (\c) at (\pos,0);
\draw[solid](\c)node[below]{$\res$}--++(0,1);
\ifnum\res=#3\relax\else
\draw[ghost](\pos+#1,0)--++(0,1)node[above]{$\res$};
\fi
\ifnum\c=1\relax
\draw[redstring](\pos+#1/15,0)node[below]{$\res$}--++(0,1);
\fi
}
\foreach \pt [evaluate=\pt as \good using {\residues[\pt-1]==#3 ? 0 : 1}
] in {#2} {
\draw[solid,dot](\pt)--++(0,1);
\ifnum\good=1
\draw[ghost,dot]([shift={(#1,0)}]\pt)--++(0,1);
\fi
}
\end{tikzpicture}
}

\NewDocumentCommand\DottedIdempotentAA{ s D(){1.5} O{} mm} {
\begin{tikzpicture}[anchorbase]
\def\residues{{#5}}
\foreach \res [count=\c,
evaluate=\res as \pos using {\res<=#4 ? #2*\res-\c*#2/15 : #2*(\res-2)-\c*#2/15}
] in {#5} {
\coordinate (\c) at (\pos,0);
\draw[solid](\c)node[below]{$\res$}--++(0,1);
\ifnum\res=#4\relax\else
\draw[ghost](\pos+#2,0)--++(0,1)node[above]{$\res$};
\fi
\ifnum\c=1\relax
\IfBooleanTF{#1}{\draw[affine](\pos+#2/15,0)node[below]{$\res$}--++(0,1);}
{\draw[redstring](\pos+#2/15,0)node[below]{$\res$}--++(0,1);}
\fi
}
\foreach \pt [evaluate=\pt as \good using {\residues[\pt-1]==#4 ? 0 : 1}
] in {#3} {
\draw[solid,dot](\pt)--++(0,1);
\ifnum\good=1
\draw[ghost,dot]([shift={(#2,0)}]\pt)--++(0,1);
\fi
}
\end{tikzpicture}
}

\NewDocumentCommand\DottedIdempotentD{ s D(){3} O{} mm} {
\begin{tikzpicture}[anchorbase]
\def\residues{{#5}}
\foreach \res [count=\c, evaluate=\res as \pos using {\res==0 ? #2*\res-\c*#2/15 : #2*(\res-2)-\c*#2/15} ] in {#5} {
\coordinate (\c) at (\pos,0);
\draw[solid](\c)node[below]{$\res$}--++(0,1);
\ifnum\res=#4\relax\else
\ifnum\res=1
\draw[double,ghost](\pos+#2,0)--++(0,1)node[above]{$\res$};
\else\ifnum\res>0
\draw[ghost](\pos+#2,0)--++(0,1)node[above]{$\res$};
\fi
\fi
\fi
\ifnum\c=1\relax
\IfBooleanTF{#1}{\draw[affine](\pos+#2/15,0)node[below]{$\res$}--++(0,1);}
{\draw[redstring](\pos+#2/15,0)node[below]{$\res$}--++(0,1);}
\fi
}
\foreach \pt [evaluate=\pt as \res using {\residues[\pt-1]}] in {#3} {
\draw[solid,dot](\pt)--++(0,1);
\ifnum\res>1
\ifnum\res<#4 \draw[ghost,dot]([shift={(#2,0)}]\pt)--++(0,1);\fi
\fi
}
\end{tikzpicture}
}

\NewDocumentCommand\DoubleCrossing{ O{1ex} mmO{} mmO{} } {
\begin{tikzpicture}[centered=#1]
\draw[#2, dot/.list={#4}] (0,0)node[below]{$#3$} .. controls (1,0.5) .. (0,1);
\draw[#5, dot/.list={#7}] (1,0)node[below]{$#6$} .. controls (0,0.5) .. (1,1);
\end{tikzpicture}\bigskip
}

\newcommand\Ddots[1]{
\tikz[scale=0.9,centered,anchorbase]{\foreach \y in {1,...,#1} {
\draw[directed=\y/#1] (0,0) to (0.925,0);}
\foreach \x in {0,...,#1} {
\node[circle,inner sep=1.8pt,fill=DarkBlue] at (\x/#1,0){};}
}}

\NewDocumentCommand\dotstring{ O{1ex} mm}{%
\tikz[centered=#1]{\draw[#2,dot](0,0)node[below]{$#3$}--++(0,1);}}

\NewDocumentCommand\stringdot{ O{1ex} mmmm}{%
\tikz[centered=#1]{
\draw[#2](0,0)node[below]{$#3$}--++(0,1);
\draw[#4,dot](1,0)node[below]{$#5$}--++(0,1);
}%
}%

\DeclarePairedDelimiterX{\set}[1]{\{}{\}}{\setargs{#1}}
\NewDocumentCommand{\setargs}{>{\SplitArgument{1}{|}}m}{\setargsaux#1}
\NewDocumentCommand{\setargsaux}{mm}
{\IfNoValueTF{#2}{#1} {#1\,\delimsize|\,\mathopen{}#2}}

\def\NewTheorem#1{%
\newaliascnt{#1}{equation}%
\newtheorem{#1}[#1]{#1}%
\aliascntresetthe{#1}%
\expandafter\def\csname #1autorefname\endcsname{#1}%
}
\def\equationautorefname~#1\null{(#1)\null}

\numberwithin{equation}{subsection}

\NewTheorem{Proposition}
\NewTheorem{Theorem}
\NewTheorem{Corollary}
\NewTheorem{Lemma}
\theoremstyle{definition}
\NewTheorem{Definition}
\NewTheorem{Notation}
\NewTheorem{Conjecture}
\NewTheorem{Example}
\AtEndEnvironment{Example}{\null\hfill$\Diamond$}%
\NewTheorem{Examples}
\AtEndEnvironment{Examples}{\vskip-10mm\null\hfill$\Diamond$}%

\theoremstyle{remark}
\NewTheorem{Remark}
\NewTheorem{Assumption}

\usepackage[hypertexnames=false]{hyperref}
\usepackage{bookmark}
\hypersetup{
pdftoolbar=true,
pdfmenubar=true,
pdffitwindow=false,
pdfstartview={FitH},
pdftitle={Subdivision for KLRW algebras of affine type A},
pdfauthor={Tao Qin},
pdfsubject={},
pdfcreator={Tao Qin},
pdfproducer={Tao Qin},
pdfkeywords={KLRW algebras},
pdfnewwindow=true,
colorlinks=true,
linkcolor=mydarkblue,
citecolor=teal,
filecolor=magenta,
urlcolor=orchid,
linkbordercolor=lava,
citebordercolor=teal,
urlbordercolor=orchid,
linktocpage=true
}
\newcommand{\nnfootnote}[1]{%
\begin{NoHyper}
\renewcommand\thefootnote{}\footnote{#1}%
\addtocounter{footnote}{-1}%
\end{NoHyper}
}

\def\makeautorefname#1#2{\csdef{#1autorefname}{#2}}
\makeautorefname{section}{Section}%
\makeautorefname{subsection}{Section}%
\makeautorefname{subsubsection}{Section}%
\begin{document}

\title[Subdivision of KLRW Algebras in Affine Type A]{Subdivision of KLRW Algebras in Affine Type A}
\author[Tao Qin]{Tao Qin}

\address{The University of Sydney, School of Mathematics and Statistics F07, Office Carslaw 807, NSW 2006, Australia}
\email{tqin5133@uni.sydney.edu.au}

\begin{abstract}
In this paper, we consider the subdivision map between two KLRW algebras of type $A^{(1)}_e$ and $A^{(1)}_{e+1}$. We show that the image of an idempotent indexed by a partition under this map is still an idempotent indexed by a partition, and give the form of this new partition. Moreover, we give an equality of some graded decomposition numbers.
\end{abstract}

\nnfootnote{\textit{Mathematics Subject Classification 2020.} Primary: 16G99, 20C08; Secondary: 20C30, 20G43.}
\nnfootnote{\textit{Keywords.} KLR algebras, diagram algebras, cellular bases, Hecke and Schur algebras.}

\addtocontents{toc}{\protect\setcounter{tocdepth}{1}}

\maketitle

\tableofcontents

\ytableausetup{centertableaux,mathmode,boxsize=0.64cm}
\section{Introduction}
The KLR algebras, or quiver Hecke algebras, were introduced by Khovanov and Lauda \cite{KhLa-cat-quantum-sln-first, KhLa-cat-quantum-sln-second}, Rouquier \cite{Ro-2-kac-moody}. These algebras categorify the
negative part of the quantum groups. In type $A^{(1)}_{e}$, Brundan and Kleshchev proved that the cyclotomic KLR algebras are isomorphic to  Ariki-Koike algebras, showing that these algebras have a $\Z$-grading. Since then, there has been a major paradigm shift in representation theory: a focus on KLR algebras.

Webster introduced the (weighted) KLRW algebra in \cite{We-knot-invariants}, as a generalization of KLR algebras and categorifies the tensor products of simple highest weight modules. Bowman \cite{Bo-many-cellular-structures} proved that the cyclotomic KLRW algebras of type $A^{(1)}_e$ are cellular and used idempotent truncation to construct a family of (diagrammatic) cellular bases for corresponding cyclotomic KLR algebras. Later, Mathas and Tubbenhauer proved the cellularity of the cyclotomic KLRW algebras of type $C^{(1)}_e$in~\cite{MaTu-klrw-algebras}, and sandwiched cellularity for types $B_{\Z\geq 0},A^{(2)}_{2e},D^{(2)}_{e+1}$ in \cite{MaTu-klrw-algebras-bad}. Recently, they proved the cellularity of cyclotomic KLRW algebras of finite types in \cite{MaTu-klrw-algebras-crystal}.

In \cite[Section 4]{MaTu-klrw-algebras}, Mathas and Tubbenhauer introduced an interesting map subdivision map between KLRW algebras attached to different quivers. A special case of the subdivision construction gives a graded isomorphism between a KLRW algebra $\R$ of types $A^{(1)}_e$ and a quotient of idempotent truncation of a KLRW algebra of type $A^{(1)}_{e+1}$.  This makes it possible to compare the structures of both sides of this isomorphism.
moreover, it is appealing to do inductions of the length of quivers. 

This paper describes the connections between the simple modules, cell modules, and graded decomposition numbers of these two algebras. In particular, we show: 
\begin{itemize}
    \item The cell module $\triangle(\blam)$ is mapped to cell module $\triangle^+(\blam^+)$, see \autoref{idempotentequal}.
    \item The graded decomposition numbers are equal under certain conditions (expected to be weakened or even removed), that is $[\triangle(\blam),L(\bmu)]_q=[\triangle^+(\blam^+),L^+(\bmu^+)]_q$, see \autoref{decompositionequal}.
\end{itemize}
The structure of the paper is as follows:
\begin{itemize}
    \item In \autoref{klrw-algebras} and \autoref{CellularityResult}, we introduce the KLRW algebras and their cellularity results, directly from \cite{MaTu-klrw-algebras} with few changes. The experts can just skip these sections with a glance at \autoref{cellularitymorita}.
    \item In \autoref{Subdivision}, the subdivision map between KLRW algebras is introduced, focusing on the subdivision at the edge $0\mapsto 1$. The general subdivision at arbitrary edge is briefly discussed in \autoref{generalcase}.
    \item In \autoref{equalityofdecompositionnumber}, we give two definitions of the image of $1_\blam$ under the subdivision map and show they're equivalent. Then we show, under certain conditions, there is an equality of decomposition numbers.
    \item In \autoref{furthertopics}, we give some insight on further interesting questions not discussed in this paper.
\end{itemize}

\noindent\textbf{Acknowledgments.}
Most Grateful to Andrew Mathas. Also thanks to Daniel Tubbenhauer, Huang Lin, Ben Webster for useful discussions of KLRW algebras. 
 
The code for drawing diagrams comes from Andrew Mathas and the code for drawing abacuses comes from Matthew Fayers. The author is partly supported by the Australian Research Council.
\section{KLRW Algebras}\label{klrw-algebras}
In this section, we introduce the definition of KLRW algebras and some basic results of them. However, since \cite[Section 2]{MaTu-klrw-algebras} is already a perfect introduction to this, we shamelessly copy their results here, deleting some irrelevant results. Of course, the reader is encouraged to read the original text to have a full understanding. Another great introduction to KLRW algebras is \cite{Bo-many-cellular-structures} with slightly different notations. 
\subsection{Quiver Combinatorics}\label{SS:Quiver}
\begin{Definition}\label{D:Quiver}
An oriented quiver $\Gamma=(I,E)$ with countable vertex set $I$ and
countable edge set $E$ is \emph{symmetrizable} if it arises from a
symmetrizable generalized Cartan matrix for $a_{ij}a_{ji}<4$
by choosing an orientation on
the simply laced edges.
\end{Definition}

\begin{Notation}\label{N:Quiver}
If not stated otherwise, we fix an oriented symmetrizable quiver $\Gamma$ and $n,\ell\in\N$. (By convention, whenever $n$ or $\ell$ are zero then the notions involving them are vacuous.)
Let $e+1=\#I$ and $\#E$ denote the respective sizes of the vertex and edge sets, respectively. We allow $e$ and $\#E$ to be infinite. A \emph{residue} is an element $i\in I$ and an $n$-tuple $\bi\in I^{n}$ is a \emph{residue sequence}.
\end{Notation}

\begin{Notation}
We have three types of edges $i\to j$, $i\Rightarrow j$ and $i\Rrightarrow j$. Write $i\rightsquigarrow j$ when the multiplicity is unimportant. All of these count as one edge. In particular, the quivers arising in this way include all Dynkin quivers, except affine type $A_{1}$. We include the affine $A_{1}$ quiver by using the convention is that there are two arrows, $0\to 1$ and $1\to 0$, between the vertices $0$ and $1$.
\end{Notation}

\begin{Example}
In this paper we restrict our attention to the following quiver:
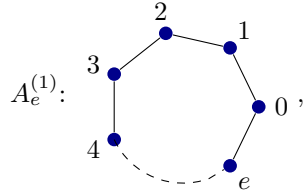
\begin{center}\label{E:Quivers}
$A^{(1)}_{e}$:
\begin{tikzpicture}[anchorbase]
\foreach \r [remember=\r as \rr] in {0,...,4} {
\node[circle,inner sep=1.8pt,fill=DarkBlue] (\r) at (360/7*\r:1){};
\node at (360/7*\r:1.3){$\r$};
\ifnum\r>0\draw(\rr)--(\r);\fi
}
\node[circle,inner sep=1.8pt,fill=DarkBlue] (6) at (360/7*6:1){};
\node at (360/7*6:1.3){$e$};
\draw[dashed](4) arc [start angle=205, end angle=290, radius=1] -- (6);
\draw(6)--(0);
\end{tikzpicture},
\\[1mm]
\end{center}

When we refer to this quiver then we fix an orientation on the simply laced edges. It is called affine type $A$.

An explicit example is:
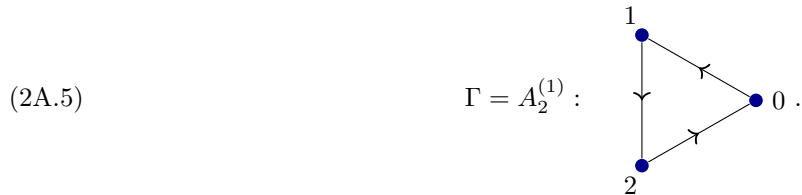
\begin{gather}\label{E:ExampleQuiver}
\Gamma=A^{(1)}_{2}:\quad
\begin{tikzpicture}[anchorbase]
\node[circle,inner sep=1.8pt,fill=DarkBlue] (0) at (360/3*0:1){};
\node at (360/3*0:1.3){$0$};
\node[circle,inner sep=1.8pt,fill=DarkBlue] (1) at (360/3*1:1){};
\node at (360/3*1:1.3){$1$};
\node[circle,inner sep=1.8pt,fill=DarkBlue] (2) at (360/3*2:1){};
\node at (360/3*2:1.3){$2$};
\draw[directed=0.5](0)--(1);
\draw[directed=0.5](1)--(2);
\draw[directed=0.5](2)--(0);
\end{tikzpicture}
.
\end{gather}
We will use this quiver in several examples below.
\end{Example}

Let $\Sym$ be the symmetric group on $\set{1,2,\dots,n}$, viewed as a Coxeter group via the presentation
\begin{gather*}
\Sym=\<s_{1},\dots,s_{n-1}|s_{j}^{2}=1,
s_{j}s_{k}=s_{k}s_{j}\text{ if }|j-k|>1,
s_{j}s_{j+1}s_{j}=s_{j+1}s_{j}s_{j+1}\>,
\end{gather*}
for all admissible $j,k$. Let $s_{j}=(j,j+1)$ and let
$(k,l)\in\Sym$ be the transposition that swaps $k$ and $l$.

Let $Q^{+}=\bigoplus_{i\in I}\N\alpha_{i}$ be the \emph{positive root lattice} of the Kac--Moody algebra determined by $\Gamma$, where $\set{\alpha_{i}|i\in I}$ are the \emph{simple roots}.
The symmetric group $\Sym$ acts on $I^{n}$ by place permutations. The
\emph{height} of
$\beta=\sum_{i\in I}b_{i}\alpha_{i}\in Q^{+}$ is $\height\beta=\sum_{i\in I}b_{i}\geq 0$. Let
$Q^{+}_{n}=\set{\beta\in Q^{+}|\height\beta=n}$. Each
$\beta\in Q^{+}_{n}$ determines the $\Sym$-orbit
$I^\beta=\set{\bi\in I^{n}|\beta=\sum_{k=1}^{n}\alpha_{i_{k}}}$ and
$I^{n}=\bigcup_{\beta\in Q^{+}_{n}}I^{\beta}$. Finally, let $\<{}_{-},{}_{-}\>$ be the Cartan pairing
associated to the quiver $\Gamma$.

\begin{Example}\label{Ex:Beta}
Consider the quiver $A_{I}$ for $I=\Z$. Let $\set{e_{i}|i\in I}$ be the standard
basis of $\R^{I}$, considered as a vector space with inner product
determined by $\<e_{i},e_{j}\>=\delta_{ij}$. Then the simple roots
$\set{\alpha_{i}|i\in I}$ can be defined as $\alpha_{i}=e_{i}-e_{i+1}$ for
$i\in I$. Taking $n=2$, we have
$Q^{+}_{2}=\{\beta_{ij}=\alpha_{i}+\alpha_{j},\beta_{i}=2\alpha_{i}|i,j\in I,i\neq j\}$,
$I^{\beta_{ij}}=\set{(i,j),(j,i)|i,j\in I,i\neq j}$ and
$I^{\beta_{i}}=\set{(i,i)|i\in I}$.
\end{Example}


\subsection{Weighted KLRW diagrams}\label{SS:Diagrams}


The main algebras considered in this paper are defined in terms of
weighted KLRW diagrams, which are the subject of this section.

\begin{Notation}\label{N:ReadingDiagrams}
In illustrations we read diagrams from bottom to top, and the concatenation $E\circ D$ of two diagrams will be viewed as stacking $E$ on top of $D$:
\begin{gather*}
E\circ D
=
\begin{tikzpicture}[anchorbase,smallnodes,rounded corners]
\node[rectangle,draw,minimum width=0.5cm,minimum height=0.5cm,ultra thick] at(0,0){\raisebox{-0.05cm}{$D$}};
\node[rectangle,draw,minimum width=0.5cm,minimum height=0.5cm,ultra thick] at(0,0.5){\raisebox{-0.05cm}{$E$}};
\end{tikzpicture}
.
\end{gather*}
In particular, left actions and left modules are given by acting from the top.
\end{Notation}

For $i\in I$, an \emph{$i$-string} is a smooth embedding
$\str\map{[0,1]}\R\times[0,1]$ such that the image of $t\in[0,1]$
belongs to $\R\times\set{t}$. We think of $\str[t]$ as a string embedded
in $\R\times[0,1]$ that is labeled by $i$. A labeled \emph{string
diagram} is an embedding of finitely many $i$-strings, for possibly
different $i\in I$, in $\R^{2}$ such that each point on these string has
a local neighborhood of one of the following two forms:
\begin{gather}\label{E:GoodCrossings}
\begin{tikzpicture}[anchorbase,smallnodes,rounded corners]
\draw[solid] (0,0)node[below]{$i$}node[below]{$\phantom{j}$} to (0,0.5)node[above,yshift=-1pt]{$\phantom{i}$};
\end{tikzpicture}
\qquad\text{or}\qquad
\begin{tikzpicture}[anchorbase,smallnodes,rounded corners]
\draw[solid] (0,0)node[below]{$i$} to (0.5,0.5);
\draw[solid] (0.5,0)node[below]{$j$} to (0,0.5)node[above,yshift=-1pt]{$\phantom{i}$};
\end{tikzpicture}
\qquad\text{for }i,i\in I.
\end{gather}
A \emph{crossing} in a diagram is a point where two strings
intersect. The right-hand diagram in \autoref{E:GoodCrossings} shows how we draw crossings.

A \emph{dot} on a string is a distinguished point on the
string that is not on any crossing or on either end point of the string.
We illustrate dots on strings as follows:
\begin{gather*}
\begin{tikzpicture}[anchorbase,smallnodes,rounded corners]
\draw[solid,dot] (0,0)node[below]{$i$} to (0,0.5)node[above,yshift=-1pt]{$\phantom{i}$};
\end{tikzpicture}.
\end{gather*}
If $\str\map{[0,1]}\R\times[0,1]$ is a string with
$\str(t)=\bigl(\str^{\prime}(t),t\bigr)$ and $\sigma\in\R$, then the
\emph{$\sigma$-shift} of $\str$ is the
string $\str+\sigma\map{[0,1]}\R\times[0,1]$ given
by $(\str+\sigma)(t)=\bigl(\str^{\prime}(t)+\sigma,t\bigr)$, for $t\in[0,1]$.

We apply this terminology below to solid, ghost and red strings, which
we now define. Before we can do this, we fix notation to ensure that the
boundary points of our strings are distinct.

\begin{Definition}\label{D:DataShifts}
Recall that we have fixed $n$, $\ell$ and a quiver $\Gamma=(I,E)$.
\begin{enumerate}

\item A \emph{solid positioning} is a $n$-tuple
$\bx=(x_{1},\dots,x_{n})\in\R^{n}$.

\item A \emph{ghost shift} for $\Gamma$ is a function
$\bsig\map{E}{\R_{\neq 0}},\epsilon\mapsto\sigma_{\epsilon}$.

\item A \emph{charge}, or red positioning, is a tuple
$\charge=(\kappa_{1},\dots,\kappa_{\ell})\in\R^{\ell}$ such that
$\kappa_{1}<\dots<\kappa_{\ell}$.

\item A \emph{loading} for $(\Gamma,\bsig)$ is a pair $(\charge,\bx)$ where
$\charge$ is a charge, $\bx$ is a solid positioning and
the numbers $x_{i}$, $x_{j}+\sigma_{\epsilon}$, $x_{j}-\sigma_{\epsilon}$ and $\kappa_{k}$ are pairwise
distinct, where $1\leq i,j\leq n$, $\epsilon\in E$, $1\leq k\leq\ell$.

\end{enumerate}
\end{Definition}

As we will see shortly, we use $\bx$, $\bsig$ and $\charge$ to determine
the boundary points of solid, ghost and red strings in the diagrams that
we consider.

\begin{Remark}\label{R:Weighting}
In \cite{We-weighted-klr} the ghost shift $\bsig$ is called a
\emph{weighting}, similar to the corresponding terminology from graph theory. We draw weighted graphs where the weights are the $\sigma_{\epsilon}$.
\end{Remark}

\begin{Notation}\label{N:Rho}
Unless otherwise stated, we fix $\bsig$ and $\brho=(\rho_{1},\dots,\rho_{\ell})\in I^{\ell}$. (Set $\rho=\rho_{1}$.) Everything below depends on these choices.
\end{Notation}

There are two extreme cases of ghost shifts that play an important role: the \emph{infinitesimal case}, where $\bsig=(\varepsilon,\dots,\varepsilon)$ for $0<\varepsilon\ll 1$, and the \emph{asymptotic case}, where $\bsig=(1/\varepsilon,\dots,1/\varepsilon)$ for $0<\varepsilon\ll 1$.

\begin{Definition}\label{D:WeightedKLRW}({See \cite[Definition 2.3]{We-weighted-klr}, \cite[Definition 4.1]{We-rouquier-dia-algebra}}.)
\label{D:WebsterDiagram}
Fix a pair $(\Gamma,\bsig)$, where $\Gamma$ is a quiver and $\bsig$ is a ghost shift. Suppose
that $\bi\in I^{n}$ and that $(\charge,\bx)$ and
$(\chargetwo,\by)$ are loadings for $(\Gamma,\bsig)$. A
\emph{weighted KLRW diagram} $D$ of
type $(\charge,\bx)\text{-}(\chargetwo,\by)$ and residue $\bi$
is a string diagram consisting of:
\begin{enumerate}

\item \emph{Solid} strings $\str_{1},\dots,\str_{n}$ such that
$\str_{k}$ is an $i_{k}$-string with $\str_k(0)=(x_{k},0)$ and
$\str_{k}(1)=(x^{\prime}_{w(k)},1)$, for some $w\in\Sym$ and $1\leq k\leq n$.

\item For each $\epsilon:i\rightsquigarrow j\in E$ with $\sigma_{\epsilon}>0$ and each $\epsilon^{\prime}:k\rightsquigarrow i\in E$ with $\sigma_{\epsilon^{\prime}}<0$, every
solid $i$-string $\str$ has \emph{ghost} $i$-strings
$\str[g]_{\epsilon}=\str+\sigma_{\epsilon}$
and $\str[g]_{\epsilon^{\prime}}=\str-\sigma_{\epsilon^{\prime}}$.

\item \emph{Red} strings $\str[r]_{1},\dots,\str[r]_{\ell}$ such that
$\str[r]_{k}$ is a $\rho_{k}$-string with
$\str[r]_k(t)=(t\kappa_{k}^{\prime}+(1-t)\kappa_k,t)$,
for $t\in[0,1]$ and $1\leq k\leq\ell$.

\item Solid strings can be decorated with finitely many dots, and
ghost strings with finitely many ghost dots, such that a dot appears
at position $\str(t)$ if and only if ghost dots appear at positions $\str[g]_{\epsilon}(t)$, for all relevant edges $\epsilon\in E$.
\end{enumerate}
\end{Definition}

We will usually simply call a weighted KLRW diagram a \emph{diagram}.
We warn the reader that the weighted KLRW diagrams have a left-right bias in that ghost are always shifted to the right.

\begin{Remark}
Recall that $i\Rightarrow j$ and $i\Rrightarrow j$ count as a single edge,
so \autoref{D:WeightedKLRW} gives only one ghost strings for such edges.
\end{Remark}

\begin{Remark}
\autoref{D:WeightedKLRW} does not allow ghost shifts of outgoing edges to be equal, such as
\begin{gather*}
\begin{tikzpicture}[anchorbase]
\node[circle,inner sep=1.8pt,fill=DarkBlue] (0) at (0,0){};
\node[circle,inner sep=1.8pt,fill=DarkBlue] (1) at (1,0){};
\node[circle,inner sep=1.8pt,fill=DarkBlue] (2) at (-1,0){};
\draw[directed=0.5](0) to node[above,yshift=-1pt]{$\sigma>0$} (1);
\draw[directed=0.5](0) to node[above,yshift=-1pt]{$\sigma>0$} (2);
\end{tikzpicture}
.
\end{gather*}
This is because loadings do not exist for such a pair $(\Gamma,\bsig)$. There are other cases
that are not allowed, depending on the sign of $\sigma$ and the orientation of the quiver.

In particular, unlike \cite{We-weighted-klr}, we do not allow ghost
shifts to be zero because this would mean that solid strings overlap with their ghost strings, which we do not want.
This said, the zero ghost shift case is captured by the infinitesimal case.
\end{Remark}

To help distinguish between the different types of strings in diagrams we draw solid strings as in \autoref{E:GoodCrossings}, ghost strings as dashed gray strings (with their labels illustrated at the top)
and red strings as thick red strings, {\cf} \autoref{E:DrawingGhost}. By
\autoref{D:WeightedKLRW}, red strings do not cross each other because,
in contrast to solid and ghost strings, we do not allow a permutation of
their endpoints. Consequently, locally, a diagram is always of one of the following forms.
\begin{gather}\label{E:DrawingGhost}
\begin{tikzpicture}[anchorbase,smallnodes,rounded corners]
\draw[solid] (0,0)node[below]{$i$}node[below]{$\phantom{j}$} to (0,0.5)node[above,yshift=-1pt]{$\phantom{i}$};
\end{tikzpicture}
,\,
\begin{tikzpicture}[anchorbase,smallnodes,rounded corners]
\draw[ghost] (0,0)node[below]{$\phantom{i}$}node[below]{$\phantom{j}$} to (0,0.5)node[above,yshift=-1pt]{$i$};
\end{tikzpicture}
,\,
\begin{tikzpicture}[anchorbase,smallnodes,rounded corners]
\draw[solid,dot] (0,0)node[below]{$i$}node[below]{$\phantom{j}$} to (0,0.5)node[above,yshift=-1pt]{$\phantom{i}$};
\end{tikzpicture}
,\,
\begin{tikzpicture}[anchorbase,smallnodes,rounded corners]
\draw[ghost,dot] (0,0)node[below]{$\phantom{i}$}node[below]{$\phantom{j}$} to (0,0.5)node[above,yshift=-1pt]{$i$};
\end{tikzpicture}
,\,
\begin{tikzpicture}[anchorbase,smallnodes,rounded corners]
\draw[redstring] (0,0)node[below]{$\rho$}node[below]{$\phantom{j}$} to (0,0.5)node[above,yshift=-1pt]{$\phantom{i}$};
\end{tikzpicture}
,\,
\begin{tikzpicture}[anchorbase,smallnodes,rounded corners]
\draw[solid] (0,0)node[below]{$i$} to (0.5,0.5);
\draw[solid] (0.5,0)node[below]{$j$} to (0,0.5)node[above,yshift=-1pt]{$\phantom{i}$};
\end{tikzpicture}
,\,
\begin{tikzpicture}[anchorbase,smallnodes,rounded corners]
\draw[ghost] (0,0)node[below]{$\phantom{i}$} to (0.5,0.5)node[above,yshift=-1pt]{$i$};
\draw[ghost] (0.5,0)node[below]{$\phantom{i}$} to (0,0.5)node[above,yshift=-1pt]{$j$};
\end{tikzpicture}
,\,
\begin{tikzpicture}[anchorbase,smallnodes,rounded corners]
\draw[ghost] (0.5,0)node[below]{\phantom{i}} to (0,0.5)node[above,yshift=-1pt]{$j$};
\draw[solid] (0,0)node[below]{$i$} to (0.5,0.5);
\end{tikzpicture}
,\,
\begin{tikzpicture}[anchorbase,smallnodes,rounded corners]
\draw[ghost] (0,0)node[below]{\phantom{i}} to (0.5,0.5)node[above,yshift=-1pt]{$i$};
\draw[solid] (0.5,0)node[below]{$j$} to (0,0.5)node[above,yshift=-1pt]{$\phantom{i}$};
\end{tikzpicture}
,\,
\begin{tikzpicture}[anchorbase,smallnodes,rounded corners]
\draw[solid] (0,0)node[below]{$i$} to (0.5,0.5);
\draw[redstring] (0.5,0)node[below]{$\rho$} to (0,0.5)node[above,yshift=-1pt]{$\phantom{i}$};
\end{tikzpicture}
,\,
\begin{tikzpicture}[anchorbase,smallnodes,rounded corners]
\draw[solid] (0.5,0)node[below]{$i$} to (0,0.5)node[above,yshift=-1pt]{$\phantom{i}$};
\draw[redstring] (0,0)node[below]{$\rho$} to (0.5,0.5);
\end{tikzpicture}
,\,
\begin{tikzpicture}[anchorbase,smallnodes,rounded corners]
\draw[ghost] (0,0)node[below]{$\phantom{i}$} to (0.5,0.5)node[above,yshift=-1pt]{$i$};
\draw[redstring] (0.5,0)node[below]{$\rho$} to (0,0.5)node[above,yshift=-1pt]{$\phantom{i}$};
\end{tikzpicture}
,\,
\begin{tikzpicture}[anchorbase,smallnodes,rounded corners]
\draw[ghost] (0.5,0)node[below]{$\phantom{i}$} to (0,0.5)node[above,yshift=-1pt]{$i$};
\draw[redstring] (0,0)node[below]{$\rho$} to (0.5,0.5);
\end{tikzpicture}
.
\end{gather}

\begin{Example}
Take $n=1=\ell$ and $\charge=\chargetwo$ and let $0<\varepsilon\ll 1$. Typical diagrams
for the following quivers and ghost shifts are:
\begin{gather*}
\begin{tikzpicture}[anchorbase]
\node[circle,inner sep=1.8pt,fill=DarkBlue] (0) at (0,0){};
\node at (0,-0.25){$i$};
\node[circle,inner sep=1.8pt,fill=DarkBlue] (1) at (1,0){};
\node at (1,-0.28){$j$};
\draw[directed=0.5](0) to node[above,yshift=-1pt]{$-1$} node[below]{$\phantom{-1}$} (1);
\end{tikzpicture}
:
\begin{tikzpicture}[anchorbase,smallnodes,rounded corners]
\draw[ghost](0.5,0)node[below]{$\phantom{i}$}--++(0,1)node[above,yshift=-1pt]{$j$};
\draw[solid,dot=0.25](0,0)node[below]{$i$}--++(1,1);
\draw[solid](-0.5,0)node[below]{$j$}--++(0,1);
\draw[redstring](0.75,0)node[below]{$\rho$}--++(0,1);
\end{tikzpicture}
,\quad
\begin{tikzpicture}[anchorbase]
\node[circle,inner sep=1.8pt,fill=DarkBlue] (0) at (0,0){};
\node at (0,-0.25){$i$};
\node[circle,inner sep=1.8pt,fill=DarkBlue] (1) at (1,0){};
\node at (1,-0.28){$j$};
\draw[directed=0.5](0) to node[above,yshift=-1pt]{$1$} node[below]{$\phantom{1}$} (1);
\end{tikzpicture}
:
\begin{tikzpicture}[anchorbase,smallnodes,rounded corners]
\draw[ghost,dot=0.25](1,0)node[below]{$\phantom{i}$}--++ (1,1)node[above,yshift=-1pt]{$i$};
\draw[solid,dot=0.25](0,0)node[below]{$i$}--++(1,1);
\draw[solid](-0.5,0)node[below]{$j$}--++(0,1);
\draw[redstring](0.75,0)node[below]{$\rho$}--++(0,1);
\end{tikzpicture}
,\quad
\begin{tikzpicture}[anchorbase]
\node[circle,inner sep=1.8pt,fill=DarkBlue] (0) at (0,0){};
\node at (0,-0.25){$i$};
\node[circle,inner sep=1.8pt,fill=DarkBlue] (1) at (1,0){};
\node at (1,-0.28){$j$};
\draw[directed=0.5](0) to node[above,yshift=-1pt]{$\varepsilon$} node[below]{$\phantom{\epsilon}$} (1);
\end{tikzpicture}
:
\begin{tikzpicture}[anchorbase,smallnodes,rounded corners]
\draw[ghost,dot=0.25](0.2,0)node[below]{$\phantom{i}$}--++ (1,1)node[above,yshift=-1pt]{$i$};
\draw[solid,dot=0.25](0,0)node[below]{$i$}--++(1,1);
\draw[solid](-0.5,0)node[below]{$j$}--++(0,1);
\draw[redstring](0.75,0)node[below]{$\rho$}--++(0,1);
\end{tikzpicture}
.
\end{gather*}
From left to right,
$\sigma_{\epsilon}=-1$, $\sigma_{\epsilon}=1$ and $\sigma_{\epsilon}=\varepsilon$, all for one edge $\epsilon:i\to j\in E$.
When there is more than one edge we obtain more ghost strings.
For example,
\begin{gather*}
\begin{tikzpicture}[anchorbase]
\node[circle,inner sep=1.8pt,fill=DarkBlue] (0) at (0,0){};
\node at (0,-0.25){$i$};
\node[circle,inner sep=1.8pt,fill=DarkBlue] (1) at (1,0.5){};
\node at (1.25,0.5){$j$};
\node[circle,inner sep=1.8pt,fill=DarkBlue] (2) at (1,0){};
\node at (1.25,0){$k$};
\node[circle,inner sep=1.8pt,fill=DarkBlue] (3) at (1,-0.5){};
\node at (1.25,-0.5){$l$};
\draw[directed=0.5](0) to[out=45,in=180] node[above,yshift=-1pt]{$1$} (1);
\draw[directed=0.5](0) to node[above,xshift=0.05cm,yshift=-0.075cm]{$-1$} (2);
\draw[directed=0.5](0) to[out=315,in=180] node[below]{$\varepsilon$} (3);
\end{tikzpicture}
:
\begin{tikzpicture}[anchorbase,smallnodes,rounded corners]
\draw[ghost,dot=0.25](0.2,0)node[below]{$\phantom{i}$}--++ (1,1)node[above,yshift=-1pt]{$i$};
\draw[ghost,dot=0.25](1,0)node[below]{$\phantom{i}$}--++ (1,1)node[above,yshift=-1pt]{$i$};
\draw[ghost,dot=0.25,dot=0.75](-1.25,0)node[below]{$\phantom{i}$}--++(-0.5,1)node[above,yshift=-1pt]{$k$};
\draw[solid,dot=0.25](0,0)node[below]{$i$}--++(1,1);
\draw[solid](-0.5,0)node[below]{$j$}--++(0,1);
\draw[solid,dot=0.25,dot=0.75](-2.25,0)node[below]{$k$}--++(-0.5,1);
\draw[solid](2,0)node[below]{$l$}--++(-0.5,1);
\draw[redstring](0.75,0)node[below]{$\rho$}--++(0,1);
\end{tikzpicture}
,
\end{gather*}
is a diagram for the illustrated pair
$(\Gamma,\bsig)$ and the corresponding loadings.
\end{Example}

\begin{Notation}\label{N:ExampleGhost}
Unless otherwise stated, in examples we usually take
$\sigma_{\epsilon}=1$.
\end{Notation}

The following simple, yet important, classes of diagrams are used throughout this paper.

\begin{Definition}
An \emph{idempotent diagram} is any diagram
with no dots and no crossings and fixed $x$-coordinate of all strings.
A \emph{straight line diagram} is any diagram
with no dots and no crossings.
\end{Definition}

\begin{Example}\label{E:StraightLine}
Prototypical examples of idempotent straight line diagrams are:
\begin{gather*}
\begin{tikzpicture}[anchorbase,smallnodes,rounded corners]
\draw[ghost](0,0)node[below]{$\phantom{i}$}--++(0,1)node[above,yshift=-1pt]{$i_{1}$};
\draw[ghost](1.2,0)node[below]{$\phantom{i}$}--++(0,1)node[above,yshift=-1pt]{$i_{2}$};
\draw[ghost](2.85,0)node[below]{$\phantom{i}$}--++(0,1)node[above,yshift=-1pt]{$i_{3}$};
\draw[ghost](4.25,0)node[below]{$\phantom{i}$}--++(0,1)node[above,yshift=-1pt]{$i_{4}$};
\draw[solid](-1,0)node[below]{$i_{1}$}--++(0,1) node[above,yshift=-1pt]{$\phantom{i_{1}}$};
\draw[solid](0.2,0)node[below]{$i_{2}$}--++(0,1);
\draw[solid](1.85,0)node[below]{$i_{3}$}--++(0,1);
\draw[solid](3.25,0)node[below]{$i_{4}$}--++(0,1);
\draw[redstring](-0.5,0)node[below]{$\rho$}--++(0,1);
\draw[redstring](1.5,0)node[below]{$\rho_{2}$}--++(0,1);
\draw[redstring](2.25,0)node[below]{$\rho_{3}$}--++(0,1);
\end{tikzpicture}
,\quad
\begin{tikzpicture}[anchorbase,smallnodes,rounded corners]
\draw[ghost](0,0)node[below]{$\phantom{i}$}--++(-0.2,1)node[above,yshift=-1pt]{$i_{1}$};
\draw[ghost](1.2,0)node[below]{$\phantom{i}$}--++(-0.3,1)node[above,yshift=-1pt]{$i_{2}$};
\draw[ghost](2.85,0)node[below]{$\phantom{i}$}--++(0.2,1)node[above,yshift=-1pt]{$i_{3}$};
\draw[ghost](4.25,0)node[below]{$\phantom{i}$}--++(0.1,1)node[above,yshift=-1pt]{$i_{4}$};
\draw[solid](-1,0)node[below]{$i_{1}$}--++(-0.2,1) node[above,yshift=-1pt]{$\phantom{i_{1}}$};
\draw[solid](0.2,0)node[below]{$i_{2}$}--++(-0.3,1);
\draw[solid](1.85,0)node[below]{$i_{3}$}--++(0.2,1);
\draw[solid](3.25,0)node[below]{$i_{4}$}--++(0.1,1);
\draw[redstring](-0.5,0)node[below]{$\rho$}--++(-0.1,1);
\draw[redstring](1.5,0)node[below]{$\rho_{2}$}--++(-0.1,1);
\draw[redstring](2.25,0)node[below]{$\rho_{3}$}--++(0.2,1);
\end{tikzpicture}
.
\end{gather*}
The left diagram is a both an idempotent and a straight line diagram. The right diagram is only a straight line diagram.
\end{Example}

\begin{Definition}\label{D:Bjstrings}
Let $\Webab$ be the set of diagrams of
type $(\charge,\bx)\text{-}(\chargetwo,\by)$ and residue $\bi$ such
that the residues of the
strings at the top of the diagrams,
when read in order from left to right, are given by $\bj$.
Whenever $D\in\Webab$ we
assume that $\bj$ is the permutation of $\bi$ determined by $D$.
\end{Definition}

\begin{Definition}
Two diagrams $D$ and $D^{\prime}$ in $\Webab$ are
equivalent if they differ by an isotopy, which is a smooth
deformation $\theta\map{[0,1]}\Webab$ such that
$\theta(0)=D$ and $\theta(1)=D^{\prime}$.
\end{Definition}

Note that isotopies cannot separate
strings, cannot change the number of dots on any string and cannot change residues. We abuse notation and write
$\Webab$ for the corresponding set of equivalence
classes under isotopy. Note that, up to isotopy, there is a unique idempotent diagram
$\1_{(\charge,\bx),\bi}\in\Webab[(\charge,\bx),\bi][(\charge,\bx),\bi]$
that has no dots and no crossings.

\begin{Example}
Let $n=1=\ell$. Then we have one solid (thus, also one ghost by \autoref{N:ExampleGhost}), and one red string, which we set to be at position $\charge=\chargetwo=(0)$. Let $\bx=\by=(-1)$. Then, for $i\in I$,
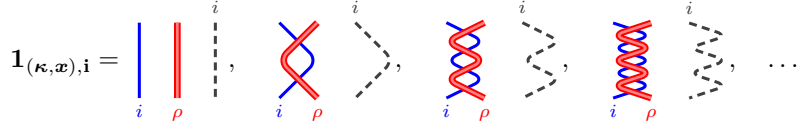
\begin{gather*}
\1_{(\charge,\bx),\bi}
=
\begin{tikzpicture}[anchorbase,smallnodes,rounded corners]
\draw[ghost](0.5,1)node[above,yshift=-1pt]{$i$}--++(0,-1)node[below]{$\phantom{i}$};
\draw[solid](-0.5,1)--++(0,-1)node[below]{$i$};
\draw[redstring](0,0)node[below]{$\rho$}--++(0,1);
\end{tikzpicture}
,\quad
\begin{tikzpicture}[anchorbase,smallnodes,rounded corners]
\draw[ghost](0.5,1)node[above,yshift=-1pt]{$i$}--++(0.5,-0.5)--++(-0.5,-0.5)node[below]{$\phantom{i}$};
\draw[solid](-0.5,1)--++(0.5,-0.5)--++(-0.5,-0.5)node[below]{$i$};
\draw[redstring](0,1)--++(-0.5,-0.5)--++(0.5,-0.5)node[below]{$\rho$};
\end{tikzpicture}
,\quad
\begin{tikzpicture}[anchorbase,smallnodes,rounded corners]
\draw[ghost](0.5,1)node[above,yshift=-1pt]{$i$}--++(0.5,-0.25)--++(-0.5,-0.25)--++(0.5,-0.25)--++(-0.5,-0.25)node[below]{$\phantom{i}$};
\draw[solid](-0.5,1)--++(0.5,-0.25)--++(-0.5,-0.25)--++(0.5,-0.25)--++(-0.5,-0.25)node[below]{$i$};
\draw[redstring](0,1)--++(-0.5,-0.25)--++(0.5,-0.25)--++(-0.5,-0.25)--++(0.5,-0.25)node[below]{$\rho$};
\end{tikzpicture}
,\quad
\begin{tikzpicture}[anchorbase,smallnodes,rounded corners]
\draw[ghost](0.5,1)node[above,yshift=-1pt]{$i$}--++(0.5,-0.166)--++(-0.5,-0.166)--++(0.5,-0.166)--++(-0.5,-0.166)--++(0.5,-0.166)--++(-0.5,-0.166)node[below]{$\phantom{i}$};
\draw[solid](-0.5,1)--++(0.5,-0.166)--++(-0.5,-0.166)--++(0.5,-0.166)--++(-0.5,-0.166)--++(0.5,-0.166)--++(-0.5,-0.166)node[below]{$i$};
\draw[redstring](0,1)--++(-0.5,-0.166)--++(0.5,-0.166)--++(-0.5,-0.166)--++(0.5,-0.166)--++(-0.5,-0.166)--++(0.5,-0.166)node[below]{$\rho$};
\end{tikzpicture}
,\quad
\dots
\end{gather*}
are examples of diagrams in $\Webab$ that are not equivalent.
\end{Example}

Let $D\in\Webab[(\charge,\bx),\bi][(\chargetwo,\by),\bj]$ and $E\in\Webab[(\chargetwo,\by),\bj][(\chargethree,\bz),\bk]$
be diagrams. Then $E\circ D\in\Webab[(\charge,\bx),\bi][(\chargethree,\bz),\bk]$
is obtained by gluing $D$ under $E$ (see \autoref{N:ReadingDiagrams}) and then rescaling. If $D\in\Webab[(\charge,\bx),\bi][(\chargetwo,\by),\bj]$, then $D=\1_{(\chargetwo,\by),\bj}\circ D\circ\1_{(\charge,\bx),\bi}$.


\subsection{Weighted KLRW algebras}\label{SS:WebsterAlgebras}


Recall that we have fixed {\eg} $\bsig$ and $\brho$ as in \autoref{N:Rho} (the ghost shift and the labels for the red strings). We additionally need:

\begin{Notation}
Fix a commutative integral domain $R$, for example $R=\Z$.
Throughout the rest of this section
we also fix $\beta\in Q^{+}_{n}$ of height $n$ as in \autoref{SS:Quiver} (the labels for the solid and ghost strings), and a finite non-empty set $X$ set of loadings, called the \emph{positioning}.
\end{Notation}

We define
\begin{gather*}
\Web_{\beta}^{\brho}(X)
=
\bigcup_{\bx,\by\in X}
\bigcup_{\bi,\bj\in I^{\beta}}
\Webaa.
\end{gather*}
In particular, $\1_{\bx,\bi}\in\Web_{\beta}^{\brho}(X)$, whenever
$\bx\in X$ and $\bi\in I^{\beta}$.

For $D\in\Webaa$ define $y_{r}D$ to be the diagram obtained from $D$
by concatenating with a dotted idempotent on top that has a dot on
the $r$th solid string and a ghost dot on the $r$th ghost string. Extend this notation so that $f(y_{1},\dots,y_{n})D$ is
the evident linear combination of diagrams for any polynomial
$f(u_{1},\dots,u_{n})\in R[u_{1},\dots,u_{n}]$.

\begin{Example}
For example,
\begin{gather}\label{E:Unsteady}
y_{3}y_{4}^{2}\1_{\bx,\bi}
=
\begin{tikzpicture}[anchorbase,smallnodes,rounded corners]
\draw[ghost](0,1)node[above,yshift=-1pt]{$i_{1}$}--++(0,-1) node[below]{$\phantom{i}$};
\draw[ghost](1.2,1)node[above,yshift=-1pt]{$i_{2}$}--++(0,-1) node[below]{$\phantom{i}$};
\draw[ghost,dot](2.85,1)node[above,yshift=-1pt]{$i_{3}$}--++(0,-1) node[below]{$\phantom{i}$};
\draw[ghost,dot=0.25,dot](4.25,1)node[above,yshift=-1pt]{$i_{4}$}--++(0,-1) node[below]{$\phantom{i}$};
\draw[solid](-1,1)node[above,yshift=-1pt]{$\phantom{i}$}--++(0,-1) node[below]{$i_{1}$};
\draw[solid](0.2,1)--++(0,-1) node[below]{$i_{2}$};
\draw[solid,dot](1.85,1)--++(0,-1) node[below]{$i_{3}$};
\draw[solid,dot=0.25,dot](3.25,1)--++(0,-1) node[below]{$i_{4}$};
\draw[redstring](-0.5,1)--++(0,-1) node[below]{$\rho_{1}$};
\draw[redstring](1.5,1)--++(0,-1) node[below]{$\rho_{2}$};
\draw[redstring](2.25,1)--++(0,-1) node[below]{$\rho_{3}$};
\end{tikzpicture}
,
\end{gather}
where $\bx$ and $\bi$ can be read-off from the illustration \autoref{E:Unsteady}.
\end{Example}

Let $\mathbf{d}\in(\N)^{e}$ be the symmetrizer of the Kac--Moody data associated to $\Gamma$.
For $i,j\in I$ fix polynomials $Q_{ij}(u,v)\in R[u,v]$, called \emph{$Q$-polynomials}, where $u$ and
$v$ are indeterminates (all of our variables appearing in polynomial rings will be indeterminates) of degrees
$2d_{i}$ respectively $2d_{j}$, such that:

\begin{enumerate}

\item We assume that $Q_{ii}(u,v)=0$ and that $Q_{ij}(u,v)$ is invertible if $i$ and $j$ are not connected by an edge in $\Gamma$.

\item For $i\neq j$ we assume that $Q_{ij}(u,v)$
is homogeneous of degree $2\<\alpha_{i},\alpha_{j}\>$ and the coefficients of all monomials are units.

\item For $i\neq j$ we assume that $Q_{ij}(u,v)=Q_{ji}(v,u)$.

\end{enumerate}

Similar conditions appear in \cite[Section 3.2.3]{Ro-2-kac-moody}, or \cite{Ro-quiver-hecke}, and \cite[Section 2.1]{We-weighted-klr}.

\begin{Example}\label{Ex:QPoly}
For the quivers in \autoref{E:Quivers} standard choices for $i\neq j$ are
\begin{gather}\label{E:QPoly}
Q_{ij}(u,v)=
\begin{cases*}
u-v & if $i\rightarrow j$,
\\
v-u & if $i\leftarrow j$,
\\
(u-v^{2}) & if $i\Rightarrow j$,
\\
(v-u^{2}) & if $i\Leftarrow j$,
\\
0 & if $i=j$,
\\
1 & otherwise.
\end{cases*}
\end{gather}
\end{Example}

Further, define polynomials $Q_{i,j,k}(u,v,w)\in R[u,v,w]$ by
\begin{gather*}
Q_{i,j,k}(u,v,w)=
\begin{cases*}
\frac{Q_{ij}(u,v)-Q_{kj}(u,w)}{w-v}& if $i=k$,\\
0& otherwise.
\end{cases*}
\end{gather*}
Below we abuse notation and write
$Q_{i,j}(\y)D=Q_{i,j}(y_{r},y_{s})D$ and
$Q_{i,j,k}(\y)D=Q_{i,j,k}(y_{r},y_{s},y_{t})D$ for the linear combination of diagrams obtained by putting dots on the corresponding strings $r,s$ and $t$ of residues $i_{r}=i$, $i_{s}=j$ and $i_{t}=k$, respectively.

We are almost ready to define the algebras that we are interested in.
Since each solid string can have several ghosts, the relations that we
use are \emph{bilocal} in the following sense:
We need to simultaneously
apply the relations in local neighborhoods around the solid strings
and in the corresponding local neighborhoods around the ghost strings.
We can only apply the relations if they can be simultaneously applied in
all neighborhoods to give new diagrams.

To ease notation we sometimes
omit solid strings or ghosts strings from relations or diagrams. In such cases the missing strings are implicit because solid and ghost strings always occur together.

\begin{Notation}\label{N:KLRWdata}
As in \autoref{D:DataShifts}, fix a quiver $\Gamma$, a ghost shift $\bsig$ and a charge $\charge$. In addition, fix $\beta\in Q^{+}_{n}$, polynomials $Q_{ij}(u,v)$ and a positioning set $X$. Even though our notation does not reflect this, the following algebras depend on all of these choices.
\end{Notation}

The following is our formulation of
\cite[Definition 2.4]{We-weighted-klr}, \cite[Definition 4.2]{We-rouquier-dia-algebra}. Recall that we write $i\rightsquigarrow j$ if there is an edge from $i$ to $j$, of any multiplicity, in $\Gamma$.

\begin{Definition}\label{D:RationalCherednik}
The \emph{weighted KLRW algebra} $\WA(X)$
is the unital associative $R$-algebra generated (as an algebra) by the diagrams
in $\Webab[\beta][\brho](X)$ with multiplication given by
\begin{gather*}
ED=
\begin{cases*}
E\circ D & if $D\in\Webaa$
and $E\in\Webaa[\by,\bj][\bz,\bk]$,
\\
0 & otherwise,
\end{cases*}
\end{gather*}
and subject to the following bilocal relations.
\begin{enumerate}

\item \label{I:Cross}
The \emph{dot sliding relation} holds, that is,
solid and ghost dots can pass through any crossing except:
\begin{gather}\label{R:DotCrossing}
\begin{tikzpicture}[anchorbase,smallnodes,rounded corners]
\draw[solid](0.5,0.5)node[above,yshift=-1pt]{$\phantom{i}$}--(0,0) node[below]{$i$};
\draw[solid,dot=0.25](0,0.5)--(0.5,0) node[below]{$i$};
\end{tikzpicture}
-
\begin{tikzpicture}[anchorbase,smallnodes,rounded corners]
\draw[solid](0.5,0.5)node[above,yshift=-1pt]{$\phantom{i}$}--(0,0) node[below]{$i$};
\draw[solid,dot=0.75](0,0.5)--(0.5,0) node[below]{$i$};
\end{tikzpicture}
=
\begin{tikzpicture}[anchorbase,smallnodes,rounded corners]
\draw[solid](0,0.5)node[above,yshift=-1pt]{$\phantom{i}$}--(0,0) node[below]{$i$};
\draw[solid](0.5,0.5)--(0.5,0) node[below]{$i$};
\end{tikzpicture}
=
\begin{tikzpicture}[anchorbase,smallnodes,rounded corners]
\draw[solid,dot=0.75](0.5,0.5)node[above,yshift=-1pt]{$\phantom{i}$}--(0,0) node[below]{$i$};
\draw[solid](0,0.5)--(0.5,0) node[below]{$i$};
\end{tikzpicture}
-
\begin{tikzpicture}[anchorbase,smallnodes,rounded corners]
\draw[solid,dot=0.25](0.5,0.5)node[above,yshift=-1pt]{$\phantom{i}$}--(0,0) node[below]{$i$};
\draw[solid](0,0.5)--(0.5,0) node[below]{$i$};
\end{tikzpicture}
.
\end{gather}

\item \label{I:DoubleCross}
The \emph{Reidemeister II relation} holds except in the following cases:
\begin{gather}\label{R:SolidSolid}
\begin{tikzpicture}[anchorbase,smallnodes,rounded corners]
\draw[solid](0,1)--++(0.5,-0.5)--++(-0.5,-0.5) node[below]{$i$};
\draw[solid](0.5,1)node[above,yshift=-1pt]{$\phantom{i}$}--++(-0.5,-0.5)--++(0.5,-0.5) node[below]{$i$};
\end{tikzpicture}
=0
.
\end{gather}
\begin{gather}\label{R:GhostSolid}
\begin{tikzpicture}[anchorbase,smallnodes,rounded corners]
\draw[ghost](0,1)node[above,yshift=-1pt]{$i$}--++(0.5,-0.5)--++(-0.5,-0.5) node[below]{$\phantom{i}$};
\draw[solid](0.5,1)--++(-0.5,-0.5)--++(0.5,-0.5) node[below]{$j$};
\end{tikzpicture}
=Q_{ij}(\y)
\begin{tikzpicture}[anchorbase,smallnodes,rounded corners]
\draw[ghost](0,1)node[above,yshift=-1pt]{$i$}--++(0,-1)node[below]{$\phantom{i}$};
\draw[solid](0.5,1)--++(0,-1)node[below]{$j$};
\end{tikzpicture}
,\quad
\begin{tikzpicture}[anchorbase,smallnodes,rounded corners]
\draw[ghost](0.5,1)node[above,yshift=-1pt]{$i$}--++(-0.5,-0.5)--++(0.5,-0.5) node[below]{$\phantom{i}$};
\draw[solid](0,1)--++(0.5,-0.5)--++(-0.5,-0.5) node[below]{$j$};
\end{tikzpicture}
=Q_{ji}(\y)
\begin{tikzpicture}[anchorbase,smallnodes,rounded corners]
\draw[ghost](0.5,1)node[above,yshift=-1pt]{$i$}--++(0,-1)node[below]{$\phantom{i}$};
\draw[solid](0,1)--++(0,-1)node[below]{$j$};
\end{tikzpicture}
\quad
\text{if $i\rightsquigarrow j$}
.
\end{gather}
\begin{gather}\label{R:RedSolid}
\begin{tikzpicture}[anchorbase,smallnodes,rounded corners]
\draw[solid](0.5,1)node[above,yshift=-1pt]{$\phantom{i}$}--++(-0.5,-0.5)--++(0.5,-0.5) node[below]{$i$};
\draw[redstring](0,1)--++(0.5,-0.5)--++(-0.5,-0.5) node[below]{$i$};
\end{tikzpicture}
=
\begin{tikzpicture}[anchorbase,smallnodes,rounded corners]
\draw[solid,dot](0.5,0)node[below]{$i$}--++(0,1)node[above,yshift=-1pt]{$\phantom{i}$};
\draw[redstring](0,0)node[below]{$i$}--++(0,1);
\end{tikzpicture}
,\quad
\begin{tikzpicture}[anchorbase,smallnodes,rounded corners]
\draw[solid](0,1)node[above,yshift=-1pt]{$\phantom{i}$}--++(0.5,-0.5)--++(-0.5,-0.5) node[below]{$i$};
\draw[redstring](0.5,1)--++(-0.5,-0.5)--++(0.5,-0.5) node[below]{$i$};
\end{tikzpicture}
=
\begin{tikzpicture}[anchorbase,smallnodes,rounded corners]
\draw[solid,dot](0,0)node[below]{$i$}--++(0,1)node[above,yshift=-1pt]{$\phantom{i}$};
\draw[redstring](0.5,0)node[below]{$i$}--++(0,1);
\end{tikzpicture}.
\end{gather}

\item \label{I:TripleCrossings}
The \emph{Reidemeister III relation} holds except in the following cases:
\begin{gather}\label{R:BraidGSG}
\begin{tikzpicture}[anchorbase,smallnodes,rounded corners]
\draw[ghost](1,1)node[above,yshift=-1pt]{$i$}--++(1,-1)node[below]{$\phantom{i}$};
\draw[ghost](2,1)node[above,yshift=-1pt]{$i$}--++(-1,-1)node[below]{$\phantom{i}$};
\draw[solid,smallnodes,rounded corners](1.5,1)--++(-0.5,-0.5)--++(0.5,-0.5)node[below]{$j$};
\end{tikzpicture}
{=}
\begin{tikzpicture}[anchorbase,smallnodes,rounded corners]
\draw[ghost](3,1)node[above,yshift=-1pt]{$\phantom{i}$}--++(1,-1)node[below]{$\phantom{i}$};
\draw[ghost](4,1)node[above,yshift=-1pt]{$i$}--++(-1,-1)node[below]{$\phantom{i}$};
\draw[solid,smallnodes,rounded corners](3.5,1)--++(0.5,-0.5)--++(-0.5,-0.5)node[below]{$j$};
\end{tikzpicture}
-Q_{i,j,i}(\y)
\begin{tikzpicture}[anchorbase,smallnodes,rounded corners]
\draw[ghost](6.2,1)node[above,yshift=-1pt]{$i$}--++(0,-1)node[below]{$\phantom{i}$};
\draw[ghost](7.2,1)node[above,yshift=-1pt]{$i$}--++(0,-1)node[below]{$\phantom{i}$};
\draw[solid](6.7,1)--++(0,-1)node[below]{$j$};
\end{tikzpicture}
,\quad
\begin{tikzpicture}[anchorbase,smallnodes,rounded corners]
\draw[solid](1,1)node[above,yshift=-1pt]{$\phantom{i}$}--++(1,-1)node[below]{$j$};
\draw[solid](2,1)--++(-1,-1)node[below]{$j$};
\draw[ghost,smallnodes,rounded corners](1.5,1)node[above,yshift=-1pt]{$i$}--++(-0.5,-0.5)--++(0.5,-0.5)node[below]{$\phantom{i}$};
\end{tikzpicture}
{=}
\begin{tikzpicture}[anchorbase,smallnodes,rounded corners]
\draw[solid](3,1)node[above,yshift=-1pt]{$\phantom{i}$}--++(1,-1)node[below]{$j$};
\draw[solid](4,1)--++(-1,-1)node[below]{$j$};
\draw[ghost,smallnodes,rounded corners](3.5,1)node[above,yshift=-1pt]{$i$}--++(0.5,-0.5)--++(-0.5,-0.5)node[below]{$\phantom{i}$};
\end{tikzpicture}
+Q_{i,j,i}(\y)
\begin{tikzpicture}[anchorbase,smallnodes,rounded corners]
\draw[solid](7.2,1)node[above,yshift=-1pt]{$\phantom{i}$}--++(0,-1)node[below]{$j$};
\draw[solid](8.2,1)--++(0,-1)node[below]{$j$};
\draw[ghost](7.7,1)node[above,yshift=-1pt]{$i$}--++(0,-1)node[below]{$\phantom{i}$};
\end{tikzpicture}
\quad
\text{if $i\rightsquigarrow j$}
.
\end{gather}
\begin{gather}\label{R:BraidSRS}
\begin{tikzpicture}[anchorbase,smallnodes,rounded corners]
\draw[solid](1,1)node[above,yshift=-1pt]{$\phantom{i}$}--++(1,-1)node[below]{$i$};
\draw[solid](2,1)--++(-1,-1)node[below]{$i$};
\draw[redstring](1.5,1)--++(-0.5,-0.5)--++(0.5,-0.5)node[below]{$i$};
\end{tikzpicture}
=
\begin{tikzpicture}[anchorbase,smallnodes,rounded corners]
\draw[solid](3,1)node[above,yshift=-1pt]{$\phantom{i}$}--++(1,-1)node[below]{$i$};
\draw[solid](4,1)--++(-1,-1)node[below]{$i$};
\draw[redstring](3.5,1)--++(0.5,-0.5)--++(-0.5,-0.5)node[below]{$i$};
\end{tikzpicture}
-
\begin{tikzpicture}[anchorbase,smallnodes,rounded corners]
\draw[solid](5,1)node[above,yshift=-1pt]{$\phantom{i}$}--++(0,-1)node[below]{$i$};
\draw[solid](6,1)--++(0,-1)node[below]{$i$};
\draw[redstring](5.5,1)--++(0,-1)node[below]{$i$};
\end{tikzpicture}
.
\end{gather}
\end{enumerate}
\end{Definition}

\begin{Example}
Bilocal relations can be tricky to apply. For example, it looks as if we can apply \autoref{R:DotCrossing} to the left-hand side of the following diagram:
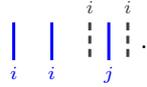
\begin{gather*}
\begin{tikzpicture}[anchorbase,smallnodes,rounded corners]
\draw[ghost](1,0.5)node[above,yshift=-1pt]{$i$}--(1,0) node[below]{$\phantom{i}$};
\draw[ghost](1.5,0.5)node[above,yshift=-1pt]{$i$}--(1.5,0) node[below]{$\phantom{i}$};
\draw[solid](0,0.5)node[above,yshift=-1pt]{$\phantom{i}$}--(0,0) node[below]{$i$};
\draw[solid](0.5,0.5)--(0.5,0) node[below]{$i$};
\draw[solid](1.25,0.5)--(1.25,0) node[below]{$j$};
\end{tikzpicture}
.
\end{gather*}
However, \autoref{R:DotCrossing} can not be used here because the two solid
$i$-strings can not be pulled arbitrarily close to one another since their ghost strings are separated by another string.
\end{Example}

\begin{Remark}\label{R:BetaVsN}
Note that we fixed $\beta\in Q^{+}_{n}$
of height $\mathrm{ht}(\beta)=n$, which
amounts to fixing the labels of the
solid strings. We sometimes write
\begin{gather*}
\WA[n](X)=
\bigoplus_{\beta\in Q^{+}_{n}}\WA(X).
\end{gather*}
Of course, anything we say about $\WA(X)$ has a corresponding statement for $\WA[n](X)$, which we will not usually state explicitly.
\end{Remark}


\subsection{The grading on \texorpdfstring{$\WA(X)$}{W}}\label{SS:WebsterDegrees}


\begin{Notation}
In this paper a \emph{graded} algebra or module
will always mean a $\Z$-graded algebra or module.
\end{Notation}

\begin{Definition}\label{D:Grading}
We endow the algebra $\WA(X)$ with a grading as follows.
The grading is defined
on the diagrams by summing over the contributions from each dot and crossing in the diagram according to the following local (not bilocal) rules:
\begin{gather*}
\deg\begin{tikzpicture}[anchorbase,smallnodes,rounded corners]
\draw[solid,dot] (0,0)node[below]{$i$} to (0,0.5)node[above,yshift=-1pt]{$\phantom{i}$};
\end{tikzpicture}
=2d_{i}
,\quad
\deg\begin{tikzpicture}[anchorbase,smallnodes,rounded corners]
\draw[ghost,dot] (0,0)node[below]{$\phantom{i}$} to (0,0.5)node[above,yshift=-1pt]{$i$};
\end{tikzpicture}=0,
\quad
\deg\begin{tikzpicture}[anchorbase,smallnodes,rounded corners]
\draw[solid] (0,0)node[below]{$i$} to (0.5,0.5);
\draw[solid] (0.5,0)node[below]{$j$} to (0,0.5)node[above,yshift=-1pt]{$\phantom{i}$};
\end{tikzpicture}
=-\delta_{i,j}2d_{i}
,\quad
\deg\begin{tikzpicture}[anchorbase,smallnodes,rounded corners]
\draw[ghost] (0,0)node[below]{$\phantom{i}$} to (0.5,0.5)node[above,yshift=-1pt]{$i$};
\draw[ghost] (0.5,0)node[below]{$\phantom{i}$} to (0,0.5)node[above,yshift=-1pt]{$j$};
\end{tikzpicture}
=0
,\\
\deg\begin{tikzpicture}[anchorbase,smallnodes,rounded corners]
\draw[ghost] (0,0)node[below]{$\phantom{i}$} to (0.5,0.5)node[above,yshift=-1pt]{$i$};
\draw[solid] (0.5,0)node[below]{$j$} to (0,0.5)node[above,yshift=-1pt]{$\phantom{i}$};
\end{tikzpicture}
=
\deg\begin{tikzpicture}[anchorbase,smallnodes,rounded corners]
\draw[ghost] (0.5,0)node[below]{$\phantom{i}$} to (0,0.5)node[above,yshift=-1pt]{$i$};
\draw[solid] (0,0)node[below]{$j$} to (0.5,0.5)node[above,yshift=-1pt]{$\phantom{i}$};
\end{tikzpicture}
=
\begin{cases}
\<\alpha_{i},\alpha_{j}\>
&\text{if $i\rightsquigarrow j$},
\\
0&\text{else},
\end{cases}
\\
\deg\begin{tikzpicture}[anchorbase,smallnodes,rounded corners]
\draw[solid] (0,0)node[below]{$i$} to (0.5,0.5)node[above,yshift=-1pt]{$\phantom{i}$};
\draw[redstring] (0.5,0)node[below]{$j$} to (0,0.5);
\end{tikzpicture}
=
\deg\begin{tikzpicture}[anchorbase,smallnodes,rounded corners]
\draw[solid] (0.5,0)node[below]{$i$} to (0,0.5)node[above,yshift=-1pt]{$\phantom{i}$};
\draw[redstring] (0,0)node[below]{$j$} to (0.5,0.5);
\end{tikzpicture}
=\tfrac{1}{2}\delta_{i,j}\<\alpha_{i},\alpha_{i}\>
,\quad
\deg\begin{tikzpicture}[anchorbase,smallnodes,rounded corners]
\draw[ghost] (0,0)node[below]{$\phantom{i}$} to (0.5,0.5)node[above,yshift=-1pt]{$i$};
\draw[redstring] (0.5,0)node[below]{$j$} to (0,0.5)node[above,yshift=-1pt]{$\phantom{i}$};
\end{tikzpicture}
=
\deg\begin{tikzpicture}[anchorbase,smallnodes,rounded corners]
\draw[ghost] (0.5,0)node[below]{$\phantom{i}$} to (0,0.5)node[above,yshift=-1pt]{$i$};
\draw[redstring] (0,0)node[below]{$j$} to (0.5,0.5)node[above,yshift=-1pt]{$\phantom{i}$};
\end{tikzpicture}
=0.
\end{gather*}
\end{Definition}

\begin{Lemma}
\autoref{D:Grading} endows $\WA(X)$ with the structure of a graded algebra.
\end{Lemma}

\begin{proof}
The algebra $\WA(X)$ is graded using
these degrees because all of the relations in
\autoref{D:RationalCherednik} are homogeneous with respect to the degree function from \autoref{D:Grading}.
\end{proof}

From now on, we consider $\WA(X)$ as a graded algebra using \autoref{D:Grading}.

\begin{Notation}
Let $\WA(X)\text{-}\mathbf{Mod}_{\Z}$ be the \emph{category of graded $\WA(X)$-modules}, where a \emph{graded category} means a category with hom-spaces enriched in the category of
graded $R$-modules. The corresponding categories of right $\WA(X)$-modules is
$\mathbf{Mod}_{\Z}\text{-}\WA(X)$.
\end{Notation}
There is a standard basis of this algebra.
\begin{Proposition}\label{standardbasis}
    $\WA(X)$ is a free module over $R$ with homogeneous basis $\WABasis$, where 
    \begin{gather}\label{E:AffineBasis}
        \WABasis=\set{\1_{\by,\bj}D(w)y_{1}^{a_{1}}\dots y_{n}^{a_{n}}\1_{\bx,\bi}|\bx,\by\in X,\bi,\bj\in I^{\beta},w\in\Sym,a_{1},\dots,a_{n}\in\N}.
    \end{gather}
\end{Proposition}
\begin{proof}
    See \cite[Lemma 3B.3, Proposition 3B.12.]{MaTu-klrw-algebras}.
\end{proof}
We will need the cyclotomic quotient of KLRW algebras, which is more interesting because it is finite dimensional and have categorification properties.

\begin{Definition}
A diagram $\1_{\bx,\bi}$ is (right) \emph{unsteady}
if it contains a solid string that can be pulled arbitrarily far to the right when the red strings are bounded by $X$.
\end{Definition}

\begin{Definition}
The \emph{cyclotomic weighted KLRW algebra} $\WAc(X)$
is the quotient of $\WA(X)$ by the two-sided ideal generated by all
diagrams that factor through an unsteady
idempotent diagram.
\end{Definition}

Similarly, we can define left unsteady diagrams. Reflecting diagrams shows that a quotient algebra that is defined by factoring out by the two-sided ideal of diagrams that factor through some left unsteady idempotent diagram is isomorphic to some cyclotomic weighted KLRW algebra and {\vive}.
We work with right unsteady diagrams because we already have a left-right bias for the ghosts strings in \autoref{D:WeightedKLRW}.

\begin{Example}
The following diagrams are unsteady and steady, respectively.
\begin{gather*}
\text{Unsteady}\colon
\begin{tikzpicture}[anchorbase,smallnodes,rounded corners]
\draw[ghost](1,0)--++(0,1)node[above,yshift=-1pt]{$i$};
\draw[solid](0.5,0)node[below]{$i$}--++(0,1)node[above,yshift=-1pt]{$\phantom{i}$};
\draw[redstring](0.25,0)node[below]{$\rho$}--++(0,1)node[above,yshift=-1pt]{$\phantom{i}$};
\end{tikzpicture}
,\quad
\text{steady}\colon
\begin{tikzpicture}[anchorbase,smallnodes,rounded corners]
\draw[ghost](0.5,0)--++(0,1)node[above,yshift=-1pt]{$i$};
\draw[solid](0,0)node[below]{$i$}--++(0,1)node[above,yshift=-1pt]{$\phantom{i}$};
\draw[redstring](0.25,0)node[below]{$\rho$}--++(0,1)node[above,yshift=-1pt]{$\phantom{i}$};
\end{tikzpicture}
.
\end{gather*}
In the right-hand diagram the solid string can not be pulled further
to the right because it cannot be pulled past the red string, so the diagram is steady. However for $\rho\neq i$
we can use \autoref{R:RedSolid}:
\begin{gather*}
\begin{tikzpicture}[anchorbase,smallnodes,rounded corners]
\draw[ghost](0.5,0)--++(0,1)node[above,yshift=-1pt]{$i$};
\draw[solid](0,0)node[below]{$i$}--++(0,1)node[above,yshift=-1pt]{$\phantom{i}$};
\draw[redstring](0.25,0)node[below]{$\rho$}--++(0,1)node[above,yshift=-1pt]{$\phantom{i}$};
\end{tikzpicture}
=
\begin{tikzpicture}[anchorbase,smallnodes,rounded corners]
\draw[ghost](0.5,0)--++(0.5,0.35)--++(0,0.3)--++(-0.5,0.35)node[above,yshift=-1pt]{$i$};
\draw[solid](0,0)node[below]{$i$}--++(0.5,0.35)--++(0,0.3)--++(-0.5,0.35)node[above,yshift=-1pt]{$\phantom{i}$};
\draw[redstring](0.25,0)node[below]{$\rho$}--++(0,1)node[above,yshift=-1pt]{$\phantom{i}$};
\end{tikzpicture}
=
\begin{tikzpicture}[anchorbase,smallnodes,rounded corners]
\draw[ghost](0.5,0)--++(0,0.1)--++(2,0.25)--++(0,0.3)--++(-2,0.25)--++(0,0.1)node[above,yshift=-1pt]{$i$};
\draw[solid](0,0)node[below]{$i$}--++(0,0.1)--++(2,0.25)--++(0,0.3)--++(-2,0.25)--++(0,0.1)node[above,yshift=-1pt]{$\phantom{i}$};
\draw[redstring](0.25,0)node[below]{$\rho$}--++(0,1)node[above,yshift=-1pt]{$\phantom{i}$};
\end{tikzpicture}
\end{gather*}
so the diagram is still zero in $\WAc(X)$.
\end{Example}

We abuse notation and identify the elements of $\WA(X)$ with their images in $\WAc(X)$.

\begin{Proposition}
The algebra $\WAc(X)$ is finite dimensional.
\end{Proposition}
\begin{proof}
    See \cite[Proposition 3D.4.]{MaTu-klrw-algebras}.
\end{proof}
\section{Cellularity Result of KLRW in Type A}\label{CellularityResult}
\subsection{(Graded) Cellular Algebra}
Cellular algebras were introduced by Graham and Lehrer \cite{GrLe-cellular}, later K{\"o}nig and Xi \cite{KoXi-affine-cellular} generalized their definition to the affine case. In \cite{HuMa-klr-basis}, Hu and Mathas defined the graded cellular algebras to suit the graded world (of KLR algebras). In this section, we give the definition and show KLRW algebras of type $A^{(1)}_e$ is cellular in next subsection.  
\begin{Definition}\label{D:CellularAlgebra}(See \cite[5B.5]{MaTu-klrw-algebras} or \cite[Definition 2.1]{HuMa-klr-basis})
Let $R$ be a unital commutative ring and $A$ be a locally unital graded $R$-algebra.
A \emph{graded affine cell datum} for $A$ is a quintuple $(\Pcal,T,\sand[],C,\deg)$, where:
\begin{itemize}
\item $\Pcal=(\Pcal,\leq)$ is a poset,
\item $T=\bigcup_{\lambda\in\Pcal}T(\lambda)$ is a collection of finite sets,
\item $\sand[]=\bigoplus_{\lambda\in\Pcal}\sand[\lambda]$ is a direct sum of
quotients of polynomial rings such that $B(\lambda)$ is a homogeneous
basis of $\sand[\lambda]$ (we write $\deg$ for the degree function on $\sand[\lambda]$),
\item $C\map{\coprod_{\lambda\in\Pcal}T(\lambda)\times B(\lambda)\times
T(\lambda)}{A};(S,\ba,T)\mapsto C^{\ba}_{ST}$ is an injective map,
\item $\deg\map{\coprod_{\lambda\in\Pcal}T(\lambda)}{\Z}$ is a function,
\end{itemize}
such that:
\begin{enumerate}[label=\upshape(GC${}_{\arabic*}$\upshape)]
\item For $\lambda\in\Pcal$, $S,T\in T(\lambda)$ and $\ba\in B(\lambda)$, $C^{\ba}_{ST}$
is homogeneous of degree $\deg(S)+\deg\ba+\deg(T)$.

\item The set $\set{C^{\ba}_{ST}|\lambda\in\Pcal,S,T\in T(\lambda),\ba\in B(\lambda)}$
is a basis of $A$.

\item For all $x\in A$ there exist scalars $r_{SU}\in R$ that do not depend
on $T$ or on $\ba$, such that
\begin{gather*}
xC^{\ba}_{ST}\equiv
\sum_{U\in T(\lambda)}r_{SU}C^{\ba}_{UT}\pmod{A^{>\lambda}},
\end{gather*}
where $A^{>\lambda}$ is the $R$-submodule of $A$ spanned by
$\set{C^{\bb}_{UV}|\mu\in\Pcal,\mu>\lambda,U,V\in T(\mu),\bb\in B(\mu)}$.

\item Let $A(\lambda)=A^{\geq\lambda}/A^{>\lambda}$, where
$A^{\geq\lambda}$ is the $R$-submodule of $A$ spanned by $\set{C^{\bb}_{UV}|\mu\in\Pcal,
\mu\geq\lambda,U,V\in T(\mu),\bb\in B(\mu)}$. Then $A(\lambda)$ is isomorphic to
$\Delta(\lambda)\otimes_{\sand[\lambda]}\nabla(\lambda)$ for free graded right and left
$\sand[\lambda]$-modules $\Delta(\lambda)$ and $\nabla(\lambda)$, respectively.

\item There is an antiinvolution $({}_{-}){}^{\star}\map{A}{A}$ of~$A$ such that
$(C^{\ba}_{ST})^{\star}\equiv C^{\ba}_{TS}\pmod{A^{>\lambda}}$, for all $S,T\in T(\lambda)$
and $\ba\in B(\lambda)$, for $\lambda\in P$. This antiinvolution identifies $\Delta(\lambda)$
and $\nabla(\lambda)$.

\end{enumerate}
The algebra $A$ is a \emph{graded affine cellular algebra} if it
has a graded affine cell datum and it is an \emph{affine cellular algebra}
if $\deg(S)=0$ for all $S\in T$. The image of $C$ is an \emph{homogeneous affine cellular basis} of $A$.

A graded \emph{cell datum} for $A$ is a graded affine cell datum such
that $\sand[\lambda]\cong R$, for all $\lambda\in\Pcal$. In this case
the image of $C$ is an \emph{homogeneous cellular basis} of $A$.
\end{Definition}
The importance of cellular algebras is that we can use them to give an easy complete classification of simple modules. We are not going to discuss details, the readers can refer to \cite[Section 2]{HuMa-klr-basis}, \cite{Ma-Iwahori-hecke},\cite{GrLe-cellular}, \cite{KoXi-affine-cellular} for more details.
\subsection{General Result}\label{generalresult}
There are two standard papers \cite{Bo-many-cellular-structures} and \cite{MaTu-klrw-algebras} establishing the cellularity results of KLRW algebras of type $A^{(1)}_{e}$. Here we adopt the conventions from the latter, with some changes. The readers are encouraged to refer to the these two papers.

Let $\Parts$ be the set of $\ell$-partition of $n$. For any $\blam\in\Parts$, we can form its box configuration (or young diagram), which is denoted by $[\blam]$. Moreover, since we have a charge $(\rho_1,\cdots,\rho_\ell)=\brho\in I^\ell$, we may fill $[\blam]$ with residues with regard to $\brho$. That is, the starting residue of each component is $\rho_i$, then increase along rows and decrease along columns.
\begin{Example}\label{eg1}
    Let $\blam=(3,1,1|4)$, $\brho=(0,3)$, $e=3$ and $n=8$. Then $[\blam]$ filled with residues (with regard to $\brho$) is as follows:
    \begin{center}
        \Tableau[scale=0.5]{{0,1,2},{3},{2}}\qquad
        \Tableau[scale=0.5,top]{{3,0,1,2}}
    \end{center}
\end{Example}
Set $\hell=\ell+n(e+1)$, and $\hell=1$ for $n=\ell=0$, so that $\Parts(\hell)$ is the
set of $\hell$-partitions of $n$.
If $n>0$, we will see that a proper subset of
$\Parts(\hell)$ indexes the cells of $\WA[n](X)$, not
$\Parts(\hell)$ itself. Identify $\Parts$ with the
left-adjusted subset of $\Parts(\hell)$, given by having empty partitions from $\ell+1$ onward.

Define
$\affine{\charge}=
(\affine{\kappa}_{1},\dots,\affine{\kappa}_{\hell})\in\Z^{\hell}$
and
$\affine{\brho}=(\affine{\rho}_{1},\dots,\affine{\rho}_{\hell})\in I^{\hell}$ by
\begin{gather*}
\affine{\kappa}_{m}=
\begin{cases*}
\kappa_{m} & if $1\leq m\leq\ell$,\\
\kappa_{\ell}+2n(m-\ell) & otherwise,
\end{cases*}
\quad\text{and}\quad
\affine{\rho}_{m}=
\begin{cases*}
\rho_{m} & if $1\leq m\leq\ell$,\\
\floor{m-\ell-1}{n}+(e+1)\Z & otherwise.
\end{cases*}
\end{gather*}
If $\ell=0$, then $\kappa_{\ell}=0$ in the above formula, by convention.
Note that $\affine{\kappa}_{1}<\dots<\affine{\kappa}_{\hell}$.
Define the \emph{($\affine{\brho}$-)residue} of the node $(m,r,c)$ to be
$\res(m,r,c)=c-r+\affine{\rho}_{m}+(e+1)\Z\in I$.

Fix $0<\varepsilon<\tfrac{1}{2n\hell}$.
Motivated by \cite[Section 5]{Bo-many-cellular-structures}, define a
\emph{positioning function}:

\begin{Definition}\label{D:Acoordinates}
Let $\hcoord\map{\Nodes}{\R}$ be the map
\begin{equation}\label{positioningfunction}
\hcoord(m,r,c)=\affine{\kappa}_{m}+(c-r)-\tfrac{m}{\hell}-(c+r)\varepsilon.
\end{equation}
For $\blam\in\Parts(\hell)$ let
$\hcoord(\blam)=\set{\hcoord(m,r,c)|(m,r,c)\in\blam}$ and write
$\hcoord(\blam)=\set{x^{\blam}_{1}<\dots<x^{\blam}_{n}}$.
The \emph{suspension point} of $\lambda^{(m)}$ is $\affine{\kappa}_{m}$,
for $1\leq m\leq\hell$.
\end{Definition}

The residue of the string corresponds
to the node $(m,r,c)$ is $\res(m,r,c)=c-r+\rho_m$. Add ghost strings for any solid string by moving the position $\epsilon=1$ units to the right. At last, add $\ell$ red strings at position $\kappa_m$ with residue $\rho_m$ where $1\leq m\leq \ell$.
\begin{Example}
    Continue with \autoref{eg1}, assume $\charge=(0,3)$, we have $1_\blam$ as follows:
    
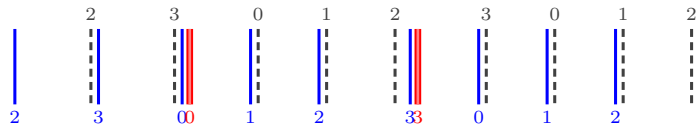
\begin{center}
        \begin{tikzpicture}[anchorbase,smallnodes,rounded corners]
            \draw[solid](-2.3,1)node[above,yshift=-1pt]{$\phantom{i}$}--++(0,-1)node[below]{$2$};
            \draw[solid](-1.2,1)node[above,yshift=-1pt]{$\phantom{i}$}--++(0,-1)node[below]{$3$};
            \draw[ghost](-1.3,1)node[above,yshift=-1pt]{$2$}--++(0,-1)node[below]{$\phantom{i}$};
            \draw[ghost](-0.2,1)node[above,yshift=-1pt]{$3$}--++(0,-1)node[below]{$\phantom{i}$};
            \draw[solid](-0.1,1)node[above,yshift=-1pt]{$\phantom{i}$}--++(0,-1)node[below]{$0$};
            \draw[redstring](0,1)--++(0,-1)node[below]{$0$};
            \draw[solid](0.8,1)node[above,yshift=-1pt]{$\phantom{i}$}--++(0,-1)node[below]{$1$};
            \draw[solid](1.7,1)node[above,yshift=-1pt]{$\phantom{i}$}--++(0,-1)node[below]{$2$};
            \draw[ghost](0.9,1)node[above,yshift=-1pt]{$0$}--++(0,-1)node[below]{$\phantom{i}$};
            \draw[ghost](1.8,1)node[above,yshift=-1pt]{$1$}--++(0,-1)node[below]{$\phantom{i}$};
            \draw[ghost](2.7,1)node[above,yshift=-1pt]{$2$}--++(0,-1)node[below]{$\phantom{i}$};
            \draw[redstring](3,1)--++(0,-1)node[below]{$3$};
            \draw[solid](2.9,1)node[above,yshift=-1pt]{$\phantom{i}$}--++(0,-1)node[below]{$3$};
            \draw[solid](3.8,1)node[above,yshift=-1pt]{$\phantom{i}$}--++(0,-1)node[below]{$0$};
            \draw[solid](4.7,1)node[above,yshift=-1pt]{$\phantom{i}$}--++(0,-1)node[below]{$1$};
            \draw[solid](5.6,1)node[above,yshift=-1pt]{$\phantom{i}$}--++(0,-1)node[below]{$2$};
            \draw[ghost](3.9,1)node[above,yshift=-1pt]{$3$}--++(0,-1)node[below]{$\phantom{i}$};
            \draw[ghost](4.8,1)node[above,yshift=-1pt]{$0$}--++(0,-1)node[below]{$\phantom{i}$};
            \draw[ghost](5.7,1)node[above,yshift=-1pt]{$1$}--++(0,-1)node[below]{$\phantom{i}$};
            \draw[ghost](6.6,1)node[above,yshift=-1pt]{$2$}--++(0,-1)node[below]{$\phantom{i}$};
        \end{tikzpicture}
    \end{center}
\end{Example}
The position function \autoref{positioningfunction} implies that the construction of $1_\blam$ coincides with the 'pulling to the right' principle. This means, for any component $\lambda^{(m)}$ of $\blam$, we first put a red string of residue $\rho_m$ at position $\charge_m$, then for the first row, put a solid string corresponding to the first node to the left of the red string, then use KLRW relations to pull it as right as possible, then the second, third. For the second row, the ghost string corresponding to the first node is blocked by the leftmost solid string of the first row, so it can't be pulled further to the right. Do the same for other nodes, we can complete the second row. All other rows are exactly the same with the second row.

The position function \autoref{positioningfunction} implies that the construction of $1_\blam$ coincides with the 'pulling to the right' principle. This means that for any component $\lambda^{(m)}$ of $\blam$, we first place a red string of residue $\rho_m$ at position $\charge_m$. Then, for the first row, we place a solid string corresponding to the first node to the left of the red string and use KLRW relations to pull it as far right as possible. We repeat this process for the second and third rows. In the second row, the ghost string corresponding to the first node is blocked by the leftmost solid string of the first row, preventing it from being pulled further to the right. By following the same procedure for the other nodes, we can complete the second row. All subsequent rows follow the same pattern as the second row.

Our definition of semistandard tableaux is
the same as in \cite[Definition 2.11]{We-rouquier-dia-algebra} or \cite[Section 1.3]{Bo-many-cellular-structures}:

\begin{Definition}
\label{D:Semistandard}
Let $\blam,\bmu\in\hParts$. A \emph{$\blam$-tableau of
type $\bmu$} is a bijection $\bT\map\blam\hcoord(\bmu)$.
Such a tableau is \emph{semistandard} if:
\begin{enumerate}

\item We have $\bT(m,1,1)\leq\kappa_{m}$ for $1\leq m\leq\ell$.

\item We have $\bT(m,r,c)+1<\bT(m,r-1,c)$ for all $(m,r,c),(m,r-1,c)\in\blam$.

\item We have $\bT(m,r,c)<\bT(m,r,c-1)+1$ for all
$(m,r,c),(m,r,c-1)\in\blam$.

\end{enumerate}
Let $\hSStd(\blam,\bmu)$ be the set of semistandard $\blam$-tableaux of
type $\bmu$ and set $\hSStd(\blam)=\bigcup_{\bmu}\hSStd(\blam,\bmu)$.
\end{Definition}

\begin{Example}\label{Ex:BigPositioning2}
Let $n=7$, $\ell=2$ and take the quiver from \autoref{E:ExampleQuiver}
with the indicated labeling. Let $\charge=(1,2)$ and $\rho=(1,2)$. Then
$\hell=23$, $\affine{\charge}=(1,2,14,28,42,\dots,294)$ and
$\affine{\brho}=(1,2,0^{7},1^{7},2^{7})$, where the exponents indicate
repeated entries. Take $(3,2|1^{2})$
and consider it as the $23$-partition
$\blam=(3,1|\emptyset|1^{2}|\emptyset|\dots|\emptyset)$.
Three $\blam$-tableaux of type $\lambda$ are
\begin{gather*}
\bT=
\begin{tikzpicture}[anchorbase,scale=0.75]
\draw[very thick] (-0.5,-1.5) to (-1,-2) to (-0.5,-2.5) to (0,-2);
\draw[very thick] (0,-1) to (-0.5,-1.5) to (0,-2) to (0.5,-1.5) to (0,-1);
\draw[very thick] (0.5,-1.5) to (1,-2) to (0.5,-2.5) to (0,-2);
\draw[very thick] (1,-2) to (1.5,-2.5) to (1,-3) to (0.5,-2.5);
\draw[very thick] (-0.5,-2.5) to (0,-3) to (0.5,-2.5);
\draw[very thick] (2.5,-1.5) to (2,-2) to (2.5,-2.5) to (3,-2);
\draw[very thick] (3,-1) to (2.5,-1.5) to (3,-2) to (3.5,-1.5) to (3,-1);
\node at (0,-1.5){\scalebox{0.7}{$0.94$}};
\node at (0,-2.5){\scalebox{0.7}{$0.92$}};
\node at (-0.5,-2){\scalebox{0.7}{$\text{-}0.07$}};
\node at (0.5,-2){\scalebox{0.7}{$1.93$}};
\node at (1,-2.5){\scalebox{0.7}{$2.92$}};
\node at (2.5,-2){\scalebox{0.65}{$12.84$}};
\node at (3,-1.5){\scalebox{0.65}{$13.85$}};
\end{tikzpicture}
,\quad
\bT^{\prime}=
\begin{tikzpicture}[anchorbase,scale=0.75]
\draw[very thick] (-0.5,-1.5) to (-1,-2) to (-0.5,-2.5) to (0,-2);
\draw[very thick] (0,-1) to (-0.5,-1.5) to (0,-2) to (0.5,-1.5) to (0,-1);
\draw[very thick] (0.5,-1.5) to (1,-2) to (0.5,-2.5) to (0,-2);
\draw[very thick] (1,-2) to (1.5,-2.5) to (1,-3) to (0.5,-2.5);
\draw[very thick] (-0.5,-2.5) to (0,-3) to (0.5,-2.5);
\draw[very thick] (2.5,-1.5) to (2,-2) to (2.5,-2.5) to (3,-2);
\draw[very thick] (3,-1) to (2.5,-1.5) to (3,-2) to (3.5,-1.5) to (3,-1);
\node at (0,-1.5){\scalebox{0.7}{$0.92$}};
\node at (0,-2.5){\scalebox{0.7}{$0.94$}};
\node at (-0.5,-2){\scalebox{0.7}{$\text{-}0.07$}};
\node at (0.5,-2){\scalebox{0.7}{$1.93$}};
\node at (1,-2.5){\scalebox{0.7}{$2.92$}};
\node at (2.5,-2){\scalebox{0.65}{$12.84$}};
\node at (3,-1.5){\scalebox{0.65}{$13.85$}};
\end{tikzpicture}
,\quad
\bT^{\prime\prime}=
\begin{tikzpicture}[anchorbase,scale=0.75]
\draw[very thick] (-0.5,-1.5) to (-1,-2) to (-0.5,-2.5) to (0,-2);
\draw[very thick] (0,-1) to (-0.5,-1.5) to (0,-2) to (0.5,-1.5) to (0,-1);
\draw[very thick] (0.5,-1.5) to (1,-2) to (0.5,-2.5) to (0,-2);
\draw[very thick] (1,-2) to (1.5,-2.5) to (1,-3) to (0.5,-2.5);
\draw[very thick] (-0.5,-2.5) to (0,-3) to (0.5,-2.5);
\draw[very thick] (2.5,-1.5) to (2,-2) to (2.5,-2.5) to (3,-2);
\draw[very thick] (3,-1) to (2.5,-1.5) to (3,-2) to (3.5,-1.5) to (3,-1);
\node at (0,-1.5){\scalebox{0.7}{$0.94$}};
\node at (0,-2.5){\scalebox{0.7}{$0.92$}};
\node at (-0.5,-2){\scalebox{0.7}{$2.92$}};
\node at (0.5,-2){\scalebox{0.7}{$1.93$}};
\node at (1,-2.5){\scalebox{0.7}{$\text{-}0.07$}};
\node at (2.5,-2){\scalebox{0.65}{$12.84$}};
\node at (3,-1.5){\scalebox{0.65}{$13.85$}};
\end{tikzpicture}
.
\end{gather*}
Note that, for example, $\bT(1,2,1)+1\approx 0.93<\bT(1,1,1)\approx 0.94$
and $\bT(1,1,2)\approx 1.93<\bT(1,1,1)+1\approx 1.94$. Similarly, one checks
that $\bT$ is a semistandard $\blam$-tableau of shape $\lambda$.
In contrast, $\bT^{\prime}$ is not semistandard, since {\eg}
$\bT^{\prime}(1,2,1)+1\approx 0.93>\bT^{\prime}(1,1,1)\approx 0.92$. The rightmost
tableau $\bT^{\prime\prime}$ is also not semistandard as
$\bT^{\prime\prime}(1,2,1)+1\approx 3.92>\bT^{\prime\prime}(1,1,1)\approx 0.94$.

Take now the $23$-partition of $7$ given by $\bmu=(3,2,1^{2}|\emptyset|\dots|\emptyset)$, which has
approximate coordinates $\hcoord(\bmu)\approx
\set{-2.09,-1.08,-0.07,0.92,\framebox{$0.94$},1.93,2.92}$. Consider
the $\blam$-tableaux of type $\bmu$ given by
\begin{gather*}
\bS=
\begin{tikzpicture}[anchorbase,scale=0.75]
\draw[very thick] (-0.5,-1.5) to (-1,-2) to (-0.5,-2.5) to (0,-2);
\draw[very thick] (0,-1) to (-0.5,-1.5) to (0,-2) to (0.5,-1.5) to (0,-1);
\draw[very thick] (0.5,-1.5) to (1,-2) to (0.5,-2.5) to (0,-2);
\draw[very thick] (1,-2) to (1.5,-2.5) to (1,-3) to (0.5,-2.5);
\draw[very thick] (-0.5,-2.5) to (0,-3) to (0.5,-2.5);
\draw[very thick] (2.5,-1.5) to (2,-2) to (2.5,-2.5) to (3,-2);
\draw[very thick] (3,-1) to (2.5,-1.5) to (3,-2) to (3.5,-1.5) to (3,-1);
\node at (0,-1.5){\scalebox{0.7}{$0.94$}};
\node at (0,-2.5){\scalebox{0.7}{$0.92$}};
\node at (-0.5,-2){\scalebox{0.7}{$\text{-}0.07$}};
\node at (0.5,-2){\scalebox{0.7}{$1.93$}};
\node at (1,-2.5){\scalebox{0.7}{$2.92$}};
\node at (2.5,-2){\scalebox{0.7}{$\text{-}2.09$}};
\node at (3,-1.5){\scalebox{0.7}{$\text{-}1.08$}};
\end{tikzpicture}
,\quad
\bS^{\prime}=
\begin{tikzpicture}[anchorbase,scale=0.75]
\draw[very thick] (-0.5,-1.5) to (-1,-2) to (-0.5,-2.5) to (0,-2);
\draw[very thick] (0,-1) to (-0.5,-1.5) to (0,-2) to (0.5,-1.5) to (0,-1);
\draw[very thick] (0.5,-1.5) to (1,-2) to (0.5,-2.5) to (0,-2);
\draw[very thick] (1,-2) to (1.5,-2.5) to (1,-3) to (0.5,-2.5);
\draw[very thick] (-0.5,-2.5) to (0,-3) to (0.5,-2.5);
\draw[very thick] (2.5,-1.5) to (2,-2) to (2.5,-2.5) to (3,-2);
\draw[very thick] (3,-1) to (2.5,-1.5) to (3,-2) to (3.5,-1.5) to (3,-1);
\node at (0,-1.5){\scalebox{0.7}{$\text{-}1.08$}};
\node at (0,-2.5){\scalebox{0.7}{$0.92$}};
\node at (-0.5,-2){\scalebox{0.7}{$\text{-}2.09$}};
\node at (0.5,-2){\scalebox{0.7}{$1.93$}};
\node at (1,-2.5){\scalebox{0.7}{$2.92$}};
\node at (2.5,-2){\scalebox{0.7}{$\text{-}0.07$}};
\node at (3,-1.5){\scalebox{0.7}{$0.94$}};
\end{tikzpicture}
,\quad
\bS^{\prime\prime}=
\begin{tikzpicture}[anchorbase,scale=0.75]
\draw[very thick] (-0.5,-1.5) to (-1,-2) to (-0.5,-2.5) to (0,-2);
\draw[very thick] (0,-1) to (-0.5,-1.5) to (0,-2) to (0.5,-1.5) to (0,-1);
\draw[very thick] (0.5,-1.5) to (1,-2) to (0.5,-2.5) to (0,-2);
\draw[very thick] (1,-2) to (1.5,-2.5) to (1,-3) to (0.5,-2.5);
\draw[very thick] (-0.5,-2.5) to (0,-3) to (0.5,-2.5);
\draw[very thick] (2.5,-1.5) to (2,-2) to (2.5,-2.5) to (3,-2);
\draw[very thick] (3,-1) to (2.5,-1.5) to (3,-2) to (3.5,-1.5) to (3,-1);
\node at (0,-1.5){\scalebox{0.7}{$\text{-}1.08$}};
\node at (0,-2.5){\scalebox{0.7}{$0.92$}};
\node at (-0.5,-2){\scalebox{0.7}{$\text{-}2.09$}};
\node at (0.5,-2){\scalebox{0.7}{$0.94$}};
\node at (1,-2.5){\scalebox{0.7}{$\text{-}0.07$}};
\node at (2.5,-2){\scalebox{0.7}{$1.93$}};
\node at (3,-1.5){\scalebox{0.7}{$2.92$}};
\end{tikzpicture}
.
\end{gather*}
Of these only $\bS$ is semistandard.
\end{Example}

Define an $R$-linear map $({}_{-}){}^{\star}\map{\WA(X)}{\WA(X)}$
by reflecting diagrams in the line $y=\frac{1}{2}$:
\begin{gather}\label{E:StarMap}
\left(
\begin{tikzpicture}[anchorbase,smallnodes,rounded corners]
\draw[ghost](1.5,1)node[above,yshift=-1pt]{$j$}--++(1,-1);
\draw[ghost](2.5,1)node[above,yshift=-1pt]{$i$}--++(-1,-1);
\draw[ghost](4,1)node[above,yshift=-1pt]{$k$}--++(0.5,-1);
\draw[solid](1,1)node[above,yshift=-1pt]{$\phantom{i}$}--++(1,-1)node[below]{$j$};
\draw[solid](2,1)--++(-1,-1)node[below]{$i$};
\draw[solid](3,1)--++(0.5,-1)node[below]{$k$};
\draw[redstring](1.25,1)--++(0,-1)node[below]{$\rho$};
\end{tikzpicture}
\right)^{\star}
=
\begin{tikzpicture}[anchorbase,smallnodes,rounded corners]
\draw[ghost](1.5,1)node[above,yshift=-1pt]{$i$}--++(1,-1);
\draw[ghost](2.5,1)node[above,yshift=-1pt]{$j$}--++(-1,-1);
\draw[ghost](4,0)--++(0.5,1)node[above,yshift=-1pt]{$k$};
\draw[solid](1,1)node[above,yshift=-1pt]{$\phantom{i}$}--++(1,-1)node[below]{$i$};
\draw[solid](2,1)--++(-1,-1)node[below]{$j$};
\draw[solid](3,0)node[below]{$k$}--++(0.5,1);
\draw[redstring](1.25,1)--++(0,-1)node[below]{$\rho$};
\end{tikzpicture}
.
\end{gather}
We write $\hcoord(\blam)=\set{x^{\blam}_{1}<\dots<x^{\blam}_{n}}$ in \autoref{D:Acoordinates}.
Then each tableau $\bT$ defines an associated permutation $w_{\bT}$
by tracing the map $\bT\map\blam\hcoord(\bmu)$ minimally. More explicitly,
let $y^{\ba}=y_{1}^{a_{1}}\dots y_{n}^{a_{n}}$ for $\ba=(a_{1},\dots,a_{n})\in\N^{n}$ and define:

\begin{Definition}\label{D:DST}
For $\bT\in\hSStd(\blam,\bmu)$ define the permutation $w_{\bT}\in\Sym$
by requiring that
\begin{gather*}
x^{\bmu}_{w_{\bT}(k)}=\bT(m,r,c)\text{ whenever }
x^{\blam}_{k}=\hcoord(m,r,c),
\quad\text{for $1\leq k\leq n$ and $(m,r,c)\in\blam$}.
\end{gather*}
Define
$D_{\bT}=D^{\hcoord(\blam),\bi^{\blam}}_{\hcoord(\bmu),\bi^{\bmu}}(w_{\bT})$,
a diagram in
$\Web^{\hcoord(\blam),\bi^{\blam}}_{\hcoord(\bmu),\bi^{\bmu}}$.
Given $\bS\in\hSStd(\blam,\bnu)$ and $\bT\in\hSStd(\blam,\bmu)$ set
\begin{gather*}
D_{\bS\bT}^{\ba}=D_{\bS}^{\star}y^{\ba}\1_{\blam}D_{\bT},\quad\text{ for }
\ba=(a_{1},\dots,a_{n})\in\N^{n}.
\end{gather*}
Set $D_{\bS\bT}=D_{\bS\bT}^{(0,\dots,0)}$.
\end{Definition}

The diagrams $D_{\bS\bT}^{\ba}$ are
well-defined only up to the choices of reduced expressions for $w_{\bS}$ and
$w_{\bT}$.

\begin{Example}\label{Ex:T}
Let $\bS$ be the leftmost $\blam$-tableau of type $\bmu$ in \autoref{Ex:BigPositioning2}.
The corresponding diagrams (ignoring the precise scaling, residues and ghost strings) are:
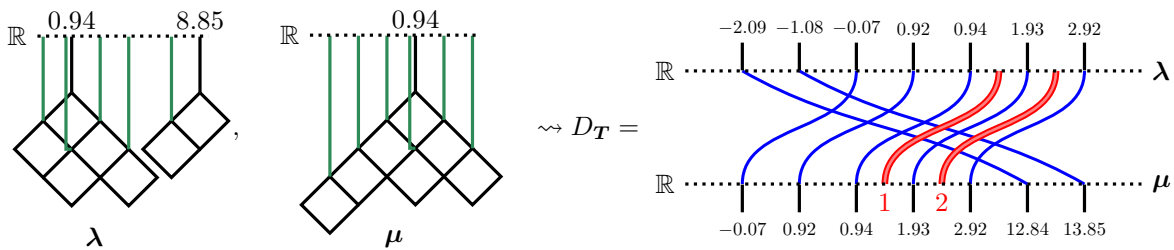
\begin{gather*}
\begin{tikzpicture}[anchorbase,scale=0.75]
\draw[very thick] (-0.5,-1.5) to (-1,-2) to (-0.5,-2.5) to (0,-2);
\draw[very thick] (0,-1) to (-0.5,-1.5) to (0,-2) to (0.5,-1.5) to (0,-1);
\draw[very thick] (0.5,-1.5) to (1,-2) to (0.5,-2.5) to (0,-2);
\draw[very thick] (1,-2) to (1.5,-2.5) to (1,-3) to (0.5,-2.5);
\draw[very thick] (-0.5,-2.5) to (0,-3) to (0.5,-2.5);
\draw[very thick] (1.75,-1.5) to (1.25,-2) to (1.75,-2.5) to (2.25,-2);
\draw[very thick] (2.25,-1) to (1.75,-1.5) to (2.25,-2) to (2.75,-1.5) to (2.25,-1);
\draw[very thick,spinach] (-0.5,-1.5) to (-0.5,-0.01);
\draw[very thick,black] (0,-1) to (0,-0.01)node[above,yshift=-1pt]{$0.94$};
\draw[very thick,spinach] (0,-2) to (-0.1,-2) to (-0.1,-0.01);
\draw[very thick,spinach] (0.5,-1.5) to (0.5,-0.01);
\draw[very thick,spinach] (1,-2) to (1,-0.01);
\draw[very thick,spinach] (1.75,-1.5) to (1.75,-0.01);
\draw[very thick,black] (2.25,-1) to (2.25,-0.01)node[above,yshift=-1pt]{$8.85$};
\draw[very thick,dotted] (-0.6,0)node[left]{$\R$} to (2.4,0);
\node at (0.4,-3.55) {$\blam$};
\end{tikzpicture}
,\quad
\begin{tikzpicture}[anchorbase,scale=0.75]
\draw[very thick] (-0.5,-1.5) to (-1,-2) to (-0.5,-2.5) to (0,-2);
\draw[very thick] (0,-1) to (-0.5,-1.5) to (0,-2) to (0.5,-1.5) to (0,-1);
\draw[very thick] (0.5,-1.5) to (1,-2) to (0.5,-2.5) to (0,-2);
\draw[very thick] (1,-2) to (1.5,-2.5) to (1,-3) to (0.5,-2.5);
\draw[very thick] (-0.5,-2.5) to (0,-3) to (0.5,-2.5);
\draw[very thick] (-1,-2) to (-1.5,-2.5) to (-1,-3) to (-0.5,-2.5);
\draw[very thick] (-1.5,-2.5) to (-2,-3) to (-1.5,-3.5) to (-1,-3);
\draw[very thick,spinach] (-0.5,-1.5) to (-0.5,-0.01);
\draw[very thick,black] (0,-1) to (0,-0.01)node[above,yshift=-1pt]{$0.94$};
\draw[very thick,spinach] (0,-2) to (-0.1,-2) to (-0.1,-0.01);
\draw[very thick,spinach] (0.5,-1.5) to (0.5,-0.01);
\draw[very thick,spinach] (1,-2) to (1,-0.01);
\draw[very thick,spinach] (-1,-2) to (-1,-0.01);
\draw[very thick,spinach] (-1.5,-2.5) to (-1.5,-0.01);
\draw[very thick,dotted] (-1.85,0)node[left]{$\R$} to (1.1,0);
\node at (-0.375,-3.6) {$\bmu$};
\end{tikzpicture}
\quad\rightsquigarrow
D_{\bT}
=
\begin{tikzpicture}[anchorbase,scale=0.75,yscale=-1]
\draw[solid](-1,0)to[out=90,in=270](-3,2);
\draw[solid](0,0)to[out=90,in=270](-2,2);
\draw[solid](1,0)to[out=90,in=270](-1,2);
\draw[solid](2,0)to[out=90,in=270](0,2);
\draw[solid](3,0)to[out=90,in=270](1,2);
\draw[solid](-3,0)to[out=30,in=210](2,2);
\draw[solid](-2,0)to[out=30,in=210](3,2);
\draw[redstring](1.5,0)node[above,yshift=-1pt]{$\phantom{i}$}to[out=90,in=270](-0.5,2)node[below]{$1$};
\draw[redstring](2.5,0)node[above,yshift=-1pt]{$\phantom{i}$}to[out=90,in=270](0.5,2)node[below]{$2$};
\draw[very thick,dotted] (-4,2)node[left]{$\R$} to (4,2)node[right]{$\bmu$};
\draw[very thick,dotted] (-4,0)node[left]{$\R$} to (4,0)node[right]{$\blam$};
\draw[very thick,black] (-3,-0.5)node[above,yshift=-1pt]{\scalebox{0.7}{$-2.09$}} to (-3,-0.01);
\draw[very thick,black] (-2,-0.5)node[above,yshift=-1pt]{\scalebox{0.7}{$-1.08$}} to (-2,-0.01);
\draw[very thick,black] (-1,-0.5)node[above,yshift=-1pt]{\scalebox{0.7}{$-0.07$}} to (-1,-0.01);
\draw[very thick,black] (0,-0.5)node[above,yshift=-1pt]{\scalebox{0.7}{$0.92$}} to (0,-0.01);
\draw[very thick,black] (1,-0.5)node[above,yshift=-1pt]{\scalebox{0.7}{$0.94$}} to (1,-0.01);
\draw[very thick,black] (2,-0.5)node[above,yshift=-1pt]{\scalebox{0.7}{$1.93$}} to (2,-0.01);
\draw[very thick,black] (3,-0.5)node[above,yshift=-1pt]{\scalebox{0.7}{$2.92$}} to (3,-0.01);
\draw[very thick,black] (-3,2.5)node[below]{\scalebox{0.7}{$-0.07$}} to (-3,2.01);
\draw[very thick,black] (-2,2.5)node[below]{\scalebox{0.7}{$0.92$}} to (-2,2.01);
\draw[very thick,black] (-1,2.5)node[below]{\scalebox{0.7}{$0.94$}} to (-1,2.01);
\draw[very thick,black] (0,2.5)node[below]{\scalebox{0.7}{$1.93$}} to (0,2.01);
\draw[very thick,black] (1,2.5)node[below]{\scalebox{0.7}{$2.92$}} to (1,2.01);
\draw[very thick,black] (2,2.5)node[below]{\scalebox{0.7}{$12.84$}} to (2,2.01);
\draw[very thick,black] (3,2.5)node[below]{\scalebox{0.7}{$13.85$}} to (3,2.01);
\end{tikzpicture}
\end{gather*}
That is, the diagram $D_{\bS}$ traces out the bijection $\bS\map\blam\hcoord(\bmu)$ with a minimal number of crossings. (We have
colored the tableaux only as a visual aid.) Consequently,
in the setup from \autoref{Ex:BigPositioning2}
we obtain
\begin{gather*}
D_{\bS\bS}^{(0,0,0,0,0,2,1)}
=
\left\{
\begin{aligned}
D_{\bS}^{\star}=
&\hspace*{0.17cm}\begin{tikzpicture}[anchorbase,smallnodes,rounded corners]
\draw[ghost](0.2,0)--++(0,1);
\draw[ghost](1.1,0)--++(0,1);
\draw[ghost](1.3,0)--++(0,1);
\draw[ghost](2.2,0)--++(0,1);
\draw[ghost](3.1,0)--++(0,1);
\draw[ghost](4.9,0)--++(0,0.35)--++(-6.9,0)--++(0,0.65);
\draw[ghost](6,0)--++(0,0.85)--++(-6.9,0)--++(0,0.15);
\draw[solid](-0.8,0)--++(0,1);
\draw[solid](0.1,0)--++(0,1);
\draw[solid](0.3,0)--++(0,1);
\draw[solid](1.2,0)--++(0,1);
\draw[solid](2.1,0)--++(0,1);
\draw[solid](3.9,0)--++(0,0.25)--++(-6.9,0)--++(0,0.75);
\draw[solid](5,0)--++(0,0.75)--++(-6.9,0)--++(0,0.25);
\draw[redstring](0.45,0)--++(0,1);
\draw[redstring](1.7,0)--++(0,1);
\draw[affine](5.2,0)--++(0,1);
\end{tikzpicture}
\\[-5pt]
y_{6}^{2}y_{7}\1_{\blam}=
&\begin{tikzpicture}[anchorbase,smallnodes,rounded corners]
\draw[white](3.72,0)--++(0,0.25)--++(-6.9,0)--++(0,0.75);
\draw[ghost](0.2,0)node[below]{$\phantom{i}$}--++(0,1)node[above,yshift=-0.05cm]{$0$};
\draw[ghost](1.1,0)node[below]{$\phantom{i}$}--++(0,1)node[above,yshift=-0.05cm]{$1$};
\draw[ghost](1.3,0)node[below]{$\phantom{i}$}--++(0,1)node[above,yshift=-0.05cm]{$1$};
\draw[ghost](2.2,0)node[below]{$\phantom{i}$}--++(0,1)node[above,yshift=-0.05cm]{$2$};
\draw[ghost](3.1,0)node[below]{$\phantom{i}$}--++(0,1)node[above,yshift=-0.05cm]{$0$};
\draw[ghost,dot=0.35,dot=0.65](4.9,0)node[below]{$\phantom{i}$}--++(0,1)node[above,yshift=-0.05cm]{$2$};
\draw[ghost,dot](6,0)node[below]{$\phantom{i}$}--++(0,1)node[above,yshift=-0.05cm]{$0$};
\draw[solid](-0.8,0)node[below]{$0$}--++(0,1)node[above,yshift=-1pt]{$\phantom{i}$};
\draw[solid](0.1,0)node[below]{$1$}--++(0,1)node[above,yshift=-1pt]{$\phantom{i}$};
\draw[solid](0.3,0)node[below]{$1$}--++(0,1)node[above,yshift=-1pt]{$\phantom{i}$};
\draw[solid](1.2,0)node[below]{$2$}--++(0,1)node[above,yshift=-1pt]{$\phantom{i}$};
\draw[solid](2.1,0)node[below]{$0$}--++(0,1)node[above,yshift=-1pt]{$\phantom{i}$};
\draw[solid,dot=0.35,dot=0.65](3.9,0)node[below]{$2$}--++(0,1)node[above,yshift=-1pt]{$\phantom{i}$};
\draw[solid,dot](5,0)node[below]{$0$}--++(0,1)node[above,yshift=-1pt]{$\phantom{i}$};
\draw[redstring](0.45,0)node[below]{$1$}--++(0,1)node[above,yshift=-1pt]{$\phantom{i}$};
\draw[redstring](1.7,0)node[below]{$2$}--++(0,1)node[above,yshift=-1pt]{$\phantom{i}$};
\draw[affine](5.2,0)node[below]{$0$}--++(0,1)node[above,yshift=-1pt]{$\phantom{i}$};
\end{tikzpicture}
\\[-5pt]
D_{\bS}=
&\hspace*{0.17cm}\begin{tikzpicture}[anchorbase,smallnodes,rounded corners,yscale=-1]
\draw[ghost](0.2,0)--++(0,1);
\draw[ghost](1.1,0)--++(0,1);
\draw[ghost](1.3,0)--++(0,1);
\draw[ghost](2.2,0)--++(0,1);
\draw[ghost](3.1,0)--++(0,1);
\draw[ghost](4.9,0)--++(0,0.35)--++(-6.9,0)--++(0,0.65);
\draw[ghost](6,0)--++(0,0.85)--++(-6.9,0)--++(0,0.15);
\draw[solid](-0.8,0)--++(0,1);
\draw[solid](0.1,0)--++(0,1);
\draw[solid](0.3,0)--++(0,1);
\draw[solid](1.2,0)--++(0,1);
\draw[solid](2.1,0)--++(0,1);
\draw[solid](3.9,0)--++(0,0.25)--++(-6.9,0)--++(0,0.75);
\draw[solid](5,0)--++(0,0.75)--++(-6.9,0)--++(0,0.25);
\draw[redstring](0.45,0)--++(0,1);
\draw[redstring](1.7,0)--++(0,1);
\draw[affine](5.2,0)--++(0,1);
\end{tikzpicture}
\end{aligned}
.
\right.
\end{gather*}
(We stack these together to obtain $D_{\bS\bS}^{(0,0,0,0,0,2,1)}$.)
Note that the strings in $\1_{\blam}$ are as far to the right as possible, with the affine red string keeping the rightmost solid $0$-string in check.
\end{Example}

We construct an affine cellular basis of the weighted
KLRW algebra $\WA[n](X)$. To this end, we choose the $Q$-polynomials as in \autoref{E:QPoly}. Recall that $\hcoord(\ell,1,n)$
is the maximal position that an $\ell$-partition of $n$ is allowed to have. Define the set
\begin{gather*}
\Affch=\set{\ba=(a_{1},\dots,a_{n})\in\N^{n}|a_{k}=0\text{ whenever }
\hcoord(\blam)_{k}\leq\hcoord(\ell,1,n)},
\end{gather*}
which will index the possible exponents in $D^{\ba}_{\bS\bT}$.
We consider the set
\begin{gather}\label{E:TBasis}
\BX=
\set[\big]{D^{\ba}_{\bS\bT}|\blam\in\hParts,\bS,\bT\in\hSStd(\blam),\ba\in\Affch}.
\end{gather}

\begin{Example}
Note that, if $\ba\in\Affch$, then $a_{k}>0$ only
if $\hcoord(\blam)_{k}$ is to the right of all possible coordinates for
$\ell$-partitions.
In particular, if $\bS$ is as in \autoref{Ex:T}, then $D_{\bS\bS}^{(0,0,0,0,0,2,1)}\in\BX$
whereas $D_{\bS\bS}^{(0,0,0,0,1,2,1)}\notin\BX$.
\end{Example}

\begin{Definition}\label{D:Dominance}
Let $\blam,\bmu\in\hParts$. Then $\blam$ \emph{is dominated by} $\bmu$,
written $\blam\ledom_{A}\bmu$, if there exists a bijection $d\colon\blam\to\bmu$ such that
$\hcoord(\alpha)\leq\hcoord(d(\alpha))$, for all $\alpha\in\blam$. Write
$\blam\ldom_{A}\bmu$ if $\blam\ledom_{A}\bmu$ and $\blam\ne\bmu$.
\end{Definition}

Consider $\hParts$ as a poset ordered by $\ledom_{A}$. We stress again that
$\hParts$ is a proper subset of the set of all $\hell$-partitions of $n$ whenever $n>0$. For $D_{\bS\bT}^{\ba}=D_{\bS}^{\star}y^{\ba}\1_{\blam}D_{\bT}$
we define $\deg\bS=\deg D_{\bS}$, $\deg\ba=\deg y^{\ba}\1_{\blam}$
and $\deg\bT=\deg D_{\bT}$.

\begin{Theorem}[{Bowman~\cite[Theorem 7.1]{Bo-many-cellular-structures}}, {Mathas-Tubbenhauer~\cite[Theorem 5B.7]{MaTu-klrw-algebras}}]\label{T:Basis}
The set $\BX$ is a homogeneous affine cellular basis of $\WA[n](X)$ with respect to the poset $(\hParts,\ledom_{A})$.
\end{Theorem}
\begin{proof}
    See \cite[Theorem 5B.7]{MaTu-klrw-algebras}.
\end{proof}

Quotienting out by the unsteady diagrams, \autoref{T:Basis} immediately implies:

\begin{Corollary}[{Mathas-Tubbenhauer~\cite[Corollary 5B.8]{MaTu-klrw-algebras}}]\label{C:TypeAKLR}
The set $\BX[{\WAc[n](X)}]=\set[\big]{D_{\bS\bT}|\blam\in\Parts,\\ \bS,\bT\in\SStd(\blam)}$
is a homogeneous cellular basis of
the cyclotomic weighted KLRW algebra $\WAc[n](X)$.
\end{Corollary}

For $\blam\in\hParts$ and $\bT\in\SStd(\blam)$
the graded cellular structure defines a
graded (left) cell module $\Delta(\blam,\bT)$ via
\begin{gather*}
\Delta_{\bT}(\blam)
=\bigl\langle D^{\ba}_{\bS\bT}|\bS\in\hSStd(\blam)\text{ and }\ba\in\Affch\bigr\rangle_{R}
\end{gather*}
with the $\WA(X)$-action defined modulo $\ldom$-higher order terms.
By \autoref{D:CellularAlgebra}, $\Delta_{\bT}(\blam)\cong\Delta_{\bT^{\prime}}(\blam)$ as $\WA(X)$-modules,
so we drop the second superscript. The $\Delta(\blam)$
are the \emph{cell modules} of $\WA(X)$.

We also need more general cell modules. To define them
let $B(\blam)=R[y^{\ba}|\ba\in\Affch]$, a polynomial subring of $R[y_{1},\dots,y_{n}]$.
By convention, we set $B(\blam)=R$ in the cyclotomic case.

Let $K$ be a simple $B(\blam)$-module.
Then the corresponding \emph{affine cell module} is
\begin{gather*}
\Delta(\blam,K)=\Delta(\blam)\otimes_{B(\lambda)}K.
\end{gather*}
Note that $\Delta(\blam,K)\cong\Delta(\blam)$ is a graded $\WA(X)$-module if and only
if $K\cong R$ is the trivial $B(\blam)$-module, since $R$ is the only graded $B(\lambda)$-module.

Each cell module $\Delta(\blam)$ has an associated
\emph{cellular pairing} determined by
\begin{gather*}
\langle
D^{\ba}_{\bS\bT},
D^{\bb}_{\bU\bV}
\rangle
=
r_{\bT\bU}
\leftrightsquigarrow
\langle
D^{\ba}_{\bS\bT},
D^{\bb}_{\bU\bV}
\rangle
=
\text{ coefficient of $1$ of }
\begin{tikzpicture}[anchorbase,scale=1]
\draw[line width=0.75,color=black,fill=cream] (0,1) to (0.25,0.5) to (0.75,0.5) to (1,1) to (0,1);
\node at (0.5,0.75){$\bU$};
\draw[line width=0.75,color=black,fill=cream] (0,1) to (0.25,1.5) to (0.75,1.5) to (1,1) to (0,1);
\node at (0.5,1.25){$\bT$};
\draw[line width=0.75,color=black,fill=cream] (0.25,0) to (0.25,0.5) to (0.75,0.5) to (0.75,0) to (0.25,0);
\node at (0.5,0.25){$\bb$};
\draw[line width=0.75,color=black,fill=cream] (0.25,1.5) to (0.25,2) to (0.75,2) to (0.75,1.5) to (0.25,1.5);
\node at (0.5,1.75){$\ba$};
\end{tikzpicture}
\text{ in }B(\blam).
\end{gather*}
See \cite[Section 2.2]{KoXi-affine-cellular} for the precise definition.
This can be extended to $\Delta(\blam,K)$ by using the identity
on $K$.
Let $\mathrm{rad}\Delta(\blam,K)$ be the associated
radical, and define
$L(\blam,K)=\Delta(\blam,K)/
\mathrm{rad}\Delta(\blam,K)$.

In the following theorem we denote grading shifts by the usual
notation using a variable $q$. Let $S_{B(\blam)}$
be a choice of simple $B(\blam)$-modules, up to isomorphism.

\begin{Theorem}[{Mathas-Tubbenhauer~\cite{MaTu-klrw-algebras}}]\label{T:SimplesA}
Suppose that $R$ is a field.
\begin{enumerate}
    \item The set $\set{L(\blam,K)|\blam\in\hParts,K\in S_{B(\blam)}}$ is a complete and non-redundant set of simple $\WA[n](X)$-modules.

    \item The set $\set{q^{s}L(\blam)|\blam\in\hParts,s\in\Z}$ is a complete and non-redundant set of graded simple $\WA[n](X)$-modules.

    \item The set $\set{L(\blam)|\blam\in\Parts}$ is a complete and non-redundant set of simple $\WAc[n](X)$-modules.

    \item The set $\set{q^{s}L(\blam)|\blam\in\Parts,s\in\Z}$ is a complete and non-redundant set of graded simple $\WAc[n](X)$-modules.

\end{enumerate}

\end{Theorem}
\begin{proof}
    See \cite[5C.1]{MaTu-klrw-algebras}.
\end{proof}

\subsection{Cellularity of KLRW algebra after subdivision}\label{cellularitymorita}
In section \autoref{Subdivision}, we will subdivide an ordinary KLRW algebra and get another KLRW algebra, whose ghost shift is always $1$ except one ghost shift very small. WLOG, we may assume this exceptional ghost shift is corresponding to the edge $1\mapsto 2$. 

We construct the idempotent diagram $1_\blam'$ exactly the same way with the normal KLRW algebra with unit ghost shift, with the exceptional change that we put the ghost $1$-string $\epsilon'$ unit to the right of the solid $1$-string. 

Use $1_\blam$ to denote the idempotent diagram in the KLRW algebra of unit ghost shift, then it is easy to see the relative positions are all strings remain unchanged from $1_\blam$ to $1_\blam'$. Indeed, one only needs to shrink the distance of solid $1$-strings and ghost $1$-strings from $1$ to $\epsilon'$, and adjust other strings not to change the relative positions and unit ghost shift.

Hence, all the cellularity results from \autoref{generalresult} also hold in this case, and we use the same notations in later sections with necessary explanation.

One can also try to construct an explicit isomorphism through straight lines isomorphisms between two families of KLRW algebras.
\section{Subdivision}\label{Subdivision}
The subdivision map is introduced by Mathas and Tubbenhauer in \cite[Section 4]{MaTu-klrw-algebras}. The map connects two KLRW algebras of type $A^{(1)}_{e}$ and $A^{(1)}_{e+1}$ and induces an isomorphism of the algebra of the former type with a quotient of the algebra of the latter type. It is likely we can make use of the map/isomorphism to conduct induction on the KLRW algebras of affine type A.

The aim of this section is to reintroduce the subdivision map. However, our definition differs from the original definition, see \cite[4C.1]{MaTu-klrw-algebras}. This new definition has the advantage that the subdivision isomorphism can be induced to the cyclotomic quotient, while the original definition in \cite{MaTu-klrw-algebras} does not have this property. We start with subdivision on the quivers.

\begin{Definition}\label{D:Subdvisionedge}
We call \Ddots{2} a \emph{subdivision} of the simply laced edge \Ddots{1}.
\end{Definition}

\begin{Definition}\label{D:Subdvision}
Let $\Gamma$ and $\Sub{\Gamma}$ be two quivers as in
\autoref{SS:Quiver}. Then $\Sub{\Gamma}$ is a
\emph{(simply laced) subdivision} of $\Gamma$ if $\Sub{\Gamma}$ is obtained from
$\Gamma$ by subdividing a finite number of simply laced edges and, potentially, relabeling the vertices.
\end{Definition}

By definition, any subdivision of $\Gamma$ is obtained by successively replacing edges \Ddots{1} in a quiver with \Ddots{2}. In particular, the orientations of the subdivided edges are compatible with the original orientation of the subdivided edge.

\begin{Examples}\label{Ex:Subdvision}
With respect to the quivers in \autoref{E:Quivers} we have:
\begin{enumerate}

\item If $e\geq d$, then the quiver $A^{(1)}_{e}$ is a subdivision of $A^{(1)}_{d}$.

\item Conversely, any subdivision of the quivers in (a) is a quiver of the same kind.

\end{enumerate}
\end{Examples}

\begin{Notation}
As a general rule, we place a bar above all of the associated Cartan data for the quiver
$\Sub{\Gamma}$. For example,
$\Sub{Q}^{+}=\bigoplus_{i\in\Sub{I}}\N\Sub{\alpha}_{i}$, and so on.
\end{Notation}

For $\Gamma$ we fix a subdivision $\Sub{\Gamma}$.
Subdivision determines two injective maps, denoted by the same symbol, $\Submap\map{I}{\Sub{I}}$
and $\Submap\map{E}{\Sub{E}}$, such
that if $r\colon i\to j$ is in $E$, then
\begin{gather*}
\begin{tikzcd}[ampersand replacement=\&,column sep=1.5cm]
\Submap(i)=i_{0}\ar[r,"\Submap(r)=r_{0}"]
\&
i_{1}\ar[r,"r_{1}"]
\&
\dots\ar[r,"r_{k-2}"]
\&
i_{k-1}\ar[r,"r_{k-1}"]
\&
i_{k}
=
\Submap(j)
\end{tikzcd}
\end{gather*}
is the subdivided edge in $\Sub{E}$. In particular, the edge $r\colon i\to j$ in $\Gamma$ is replaced with $k$ edges in $\Sub{\Gamma}$.

\begin{Definition}
Let $\beta=\sum_{i\in I_{\Gamma}}b_{i}\alpha_{i}\in Q^{+}_{\Gamma}$.
Define
\begin{gather*}
\Sub{\beta}=\sum_{i\in I_{\Gamma}}b_{i}\Sub{\alpha}_{S(i)}
\in Q^{+}_{\Sub{\Gamma}}.
\end{gather*}
\end{Definition}

\begin{Example}\label{Ex:SubdivQuiver}
Consider the quiver $\Gamma$ from \autoref{E:ExampleQuiver}, and the two subdivisions given by
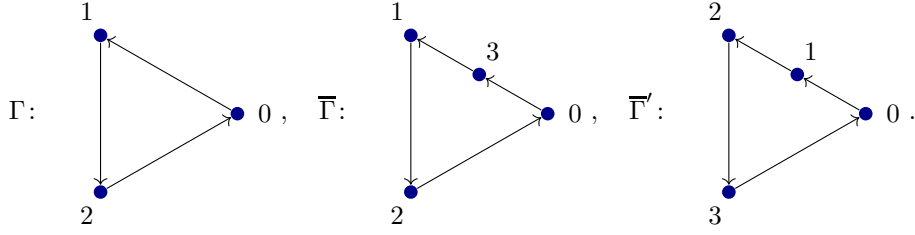
\begin{gather*}
\Gamma\colon\quad
\begin{tikzpicture}[scale=1.2,anchorbase]
\node[circle,inner sep=1.8pt,fill=DarkBlue] (0) at (360/3*0:1){};
\node at (360/3*0:1.3){$0$};
\node[circle,inner sep=1.8pt,fill=DarkBlue] (1) at (360/3*1:1){};
\node at (360/3*1:1.3){$1$};
\node[circle,inner sep=1.8pt,fill=DarkBlue] (2) at (360/3*2:1){};
\node at (360/3*2:1.3){$2$};
\draw[->](0)--(1);
\draw[->](1)--(2);
\draw[->](2)--(0);
\end{tikzpicture}
,\quad
\Sub{\Gamma}\colon\quad
\begin{tikzpicture}[scale=1.2,anchorbase]
\node[circle,inner sep=1.8pt,fill=DarkBlue] (0) at (360/3*0:1){};
\node at (360/3*0:1.3){$0$};
\node[circle,inner sep=1.8pt,fill=DarkBlue] (1) at (360/3*1:1){};
\node at (360/3*1:1.3){$1$};
\node[circle,inner sep=1.8pt,fill=DarkBlue] (2) at (360/3*2:1){};
\node at (360/3*2:1.3){$2$};
\node[circle,inner sep=1.8pt,fill=DarkBlue] (4) at (360/3*0.5:0.5){};
\node at (360/3*0.5:0.8){$3$};
\draw[->](0)--(4);
\draw[->](4)--(1);
\draw[->](1)--(2);
\draw[->](2)--(0);
\end{tikzpicture}
,\quad
\Sub{\Gamma}^{\prime}\colon\quad
\begin{tikzpicture}[scale=1.2,anchorbase]
\node[circle,inner sep=1.8pt,fill=DarkBlue] (0) at (360/3*0:1){};
\node at (360/3*0:1.3){$0$};
\node[circle,inner sep=1.8pt,fill=DarkBlue] (1) at (360/3*1:1){};
\node at (360/3*1:1.3){$2$};
\node[circle,inner sep=1.8pt,fill=DarkBlue] (2) at (360/3*2:1){};
\node at (360/3*2:1.3){$3$};
\node[circle,inner sep=1.8pt,fill=DarkBlue] (4) at (360/3*0.5:0.5){};
\node at (360/3*0.5:0.8){$1$};
\draw[->](0)--(4);
\draw[->](4)--(1);
\draw[->](1)--(2);
\draw[->](2)--(0);
\end{tikzpicture}
.
\end{gather*}
The subdivision map $\Submap$ sends $i$ to $i$, so all of the original edges keep their name, and $\Sub{\beta}=\beta$ using the Kac--Moody data for $\Sub{\Gamma}$ rather than $\Gamma$.
For $\Sub{\Gamma}^{\prime}$, the map $\Submap'$ sends $0\mapsto 0$, $1\mapsto 2$, and $2\mapsto 3$, so that the edges of $\Gamma$, and the subscripts of $\beta$, change accordingly.
\end{Example}
\begin{Remark}
    After subdivision, the labels in $\Sub{\Gamma}$ is not standard. For convenience, we may relabel the quiver by $0\mapsto 0$, $3\mapsto 1$, $1\mapsto 2$, $2\mapsto 3$, hence get the following quiver:\\
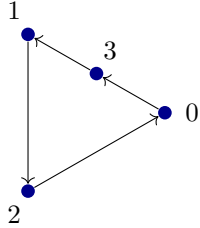
    \begin{center}
    \begin{tikzpicture}[scale=1.2,anchorbase]
        \node[circle,inner sep=1.8pt,fill=DarkBlue] (0) at (360/3*0:1){};
        \node at (360/3*0:1.3){$0$};
        \node[circle,inner sep=1.8pt,fill=DarkBlue] (1) at (360/3*1:1){};
        \node at (360/3*1:1.3){$1$};
        \node[circle,inner sep=1.8pt,fill=DarkBlue] (2) at (360/3*2:1){};
        \node at (360/3*2:1.3){$2$};
        \node[circle,inner sep=1.8pt,fill=DarkBlue] (4) at (360/3*0.5:0.5){};
        \node at (360/3*0.5:0.8){$3$};
        \draw[->](0)--(4);
        \draw[->](4)--(1);
        \draw[->](1)--(2);
        \draw[->](2)--(0);
    \end{tikzpicture}       
    \end{center}
    By abuse of notation we still denote it by $\Sub{\Gamma}$. In general, if we subdivide the edge $(i\mapsto j)$, fix the label $i$ and then add the label in the same way as in the original quiver.
\end{Remark}

Next, we introduce the subdivision map on the KLRW diagrams defined in \cite[Section 4]{MaTu-klrw-algebras}. Informally, if we subdivide the edge $i\mapsto i+1$, the subdivision map on diagrams is given by ``fattening'' the
ghost $i$-strings using the rule
\begin{gather}\label{E:Colors}
\begin{tikzpicture}[scale=1.2,anchorbase,smallnodes,rounded corners]
    \draw[ghost](1,0)--++(0,1)node[above,yshift=-1pt]{$i$};
    \draw[solid](0,0)node[below]{$i$}--++(0,1)node[above,yshift=-1pt]{$\phantom{i}$};
    \draw[redstring](0.5,0)node[below]{$\rho$}--++(0,1);
\end{tikzpicture}
\mapsto
\begin{tikzpicture}[scale=1.2,anchorbase,smallnodes,rounded corners]
    \draw[ghost](4.0,0)--++(0,1)node[above,yshift=-1pt]{$s(i)$};
    \draw[ghost,spinach](5.1,0)--++(0,1)node[above,yshift=-1pt,spinach]{$\Sub{i}$};
    \draw[solid](3.0,0)node[below]{$s(i)$}--++(0,1);
    \draw[solid,spinach](4.1,0)node[below,spinach]{$\Sub{i}$}--++(0,1);
    \draw[redstring](3.5,0)node[below]{$\rho$}--++(0,1);
    \draw[->] (4.1,-0.5)node[below]{added solid} to (4.1,-0.3);
    \draw[->] (5.1,1.5)node[above]{added ghost} to (5.1,1.3);
\end{tikzpicture}
\end{gather}
and keeping all other strings as they are, except that the $j$-strings are relabeled as $\Submap(j)$-strings. 

Notice that, the added solid $\phantom{i}$-string is always put to the right of the ghost $i$-string. Similarly, we can define the subdivision map by adding to the left consistently. But both two definitions will not induce a well-defined map to the cyclotomic
quotient. We begin with two examples:
\begin{Example}
 (1) Take $e=2=n,\ell=1$, $\blam=1^2$ and $\rho=1$, let $S$ be the subdivision map (to add on the left) at $0\mapsto 1$, so the subdivided quiver is $0\mapsto 3\mapsto 1$. It adds the solid $1$-string on the left of the ghost $0$-string. Then we have:\\
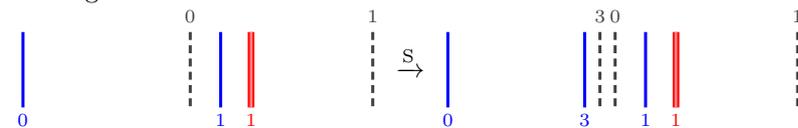
\begin{tikzpicture}[anchorbase,smallnodes,rounded corners]
    \draw[solid](1,1)node[above,yshift=-1pt]{$\phantom{i}$}--++(0,-1)node[below]{$1$};
    \draw[redstring](1.4,1)--++(0,-1)node[below]{$1$};
    \draw[solid](-1.6,1)node[above,yshift=-1pt]{$\phantom{i}$}--++(0,-1)node[below]{$0$};
    \draw[ghost](0.6,1)node[above,yshift=-1pt]{$0$}--++(0,-1)node[below]{$\phantom{i}$};
    \draw[ghost](3,1)node[above,yshift=-1pt]{$1$}--++(0,-1)node[below]{$\phantom{i}$};
\end{tikzpicture}
$\xrightarrow[]{\text{S}}$
\begin{tikzpicture}[anchorbase,smallnodes,rounded corners]
    \draw[solid](1,1)node[above,yshift=-1pt]{$\phantom{i}$}--++(0,-1)node[below]{$1$};
    \draw[redstring](1.4,1)--++(0,-1)node[below]{$1$};
    \draw[solid](-1.6,1)node[above,yshift=-1pt]{$\phantom{i}$}--++(0,-1)node[below]{$0$};
    \draw[ghost](0.6,1)node[above,yshift=-1pt]{$0$}--++(0,-1)node[below]{$\phantom{i}$};
    \draw[solid](0.2,1)node[above,yshift=-1pt]{$\phantom{i}$}--++(0,-1)node[below]{$3$};
    \draw[ghost](0.4,1)node[above,yshift=-1pt]{$3$}--++(0,-1)node[below]{$\phantom{i}$};
    \draw[ghost](3,1)node[above,yshift=-1pt]{$1$}--++(0,-1)node[below]{$\phantom{i}$};
\end{tikzpicture}\\
For simplicity, we relabel the quiver, so we replace $1$ by $2$, and replace $3$ by $1$, then the quiver is standard $0\mapsto 1\mapsto 2$ and the diagram becomes:
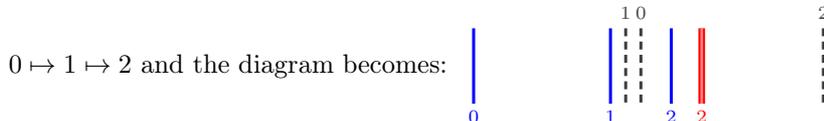
\begin{tikzpicture}[anchorbase,smallnodes,rounded corners]
    \draw[solid](1,1)node[above,yshift=-1pt]{$\phantom{i}$}--++(0,-1)node[below]{$2$};
    \draw[redstring](1.4,1)--++(0,-1)node[below]{$2$};
    \draw[solid](-1.6,1)node[above,yshift=-1pt]{$\phantom{i}$}--++(0,-1)node[below]{$0$};
    \draw[ghost](0.6,1)node[above,yshift=-1pt]{$0$}--++(0,-1)node[below]{$\phantom{i}$};
    \draw[solid](0.2,1)node[above,yshift=-1pt]{$\phantom{i}$}--++(0,-1)node[below]{$1$};
    \draw[ghost](0.4,1)node[above,yshift=-1pt]{$1$}--++(0,-1)node[below]{$\phantom{i}$};
    \draw[ghost](3,1)node[above,yshift=-1pt]{$2$}--++(0,-1)node[below]{$\phantom{i}$};
\end{tikzpicture}\\
However, the $0$- string can be pulled to the right arbitrarily far, hence unsteady.\\

(2) Let $S'$ be the subdivision map (to add on the right), add the solid $1$-string on the right of the ghost $0$-string. Take $\bmu=2$
$e=2=n,\ell=1$, and $\rho=0$. Then we have:\\
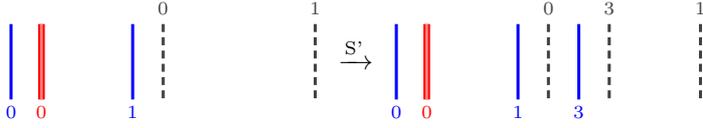
\begin{tikzpicture}[anchorbase,smallnodes,rounded corners]
    \draw[solid](1,1)node[above,yshift=-1pt]{$\phantom{i}$}--++(0,-1)node[below]{$0$};
    \draw[redstring](1.4,1)--++(0,-1)node[below]{$0$};
    \draw[solid](2.6,1)node[above,yshift=-1pt]{$\phantom{i}$}--++(0,-1)node[below]{$1$};
    \draw[ghost](3,1)node[above,yshift=-1pt]{$0$}--++(0,-1)node[below]{$\phantom{i}$};
    \draw[ghost](5,1)node[above,yshift=-1pt]{$1$}--++(0,-1)node[below]{$\phantom{i}$};
\end{tikzpicture}
$\xrightarrow[]{\text{S'}}$
\begin{tikzpicture}[anchorbase,smallnodes,rounded corners]
    \draw[solid](1,1)node[above,yshift=-1pt]{$\phantom{i}$}--++(0,-1)node[below]{$0$};
    \draw[redstring](1.4,1)--++(0,-1)node[below]{$0$};
    \draw[solid](2.6,1)node[above,yshift=-1pt]{$\phantom{i}$}--++(0,-1)node[below]{$1$};
    \draw[ghost](3,1)node[above,yshift=-1pt]{$0$}--++(0,-1)node[below]{$\phantom{i}$};
    \draw[ghost](5,1)node[above,yshift=-1pt]{$1$}--++(0,-1)node[below]{$\phantom{i}$};
    \draw[solid](3.4,1)node[above,yshift=-1pt]{$\phantom{i}$}--++(0,-1)node[below]{$3$};
    \draw[ghost](3.8,1)node[above,yshift=-1pt]{$3$}--++(0,-1)node[below]{$\phantom{i}$};\\
\end{tikzpicture}\\
However, by pulling the $3$-strings to the right, then pull the $1$-strings to the right, we see it is unsteady as well.
\end{Example}
This example demonstrates that it is not possible to add the new strings to only one particular side of the ghost $0$-string, as doing so cannot induce the subdivision isomorphism to the cyclotomic quotient. A more suitable definition should allow the side of the added string to depend on the residue pattern. {\color{Green}From now on, we will subdivide at the $0\mapsto 1$ place, though there is no difference. See \autoref{generalcase} for the general case.}
\begin{Definition}\label{closetuple}
    For an idempotent string diagram $Y$, if locally there are some strings of the following form:
    \begin{enumerate}
        \item (ghost $0$)(solid $1$)(ghost $0$)$\cdots$(solid $1$)(ghost $0$)(solid $1$)(ghost $0$), i.e.\\
        \begin{center}   
        \begin{tikzpicture}[anchorbase,smallnodes,rounded corners]
            \draw[solid](1,1)node[above,yshift=-1pt]{$\phantom{i}$}--++(0,-1)node[below]{$1$};
            \draw[ghost](0.4,1)node[above,yshift=-1pt]{$0$}--++(0,-1)node[below]{$\phantom{i}$};
            \draw[ghost](1.6,1)node[above,yshift=-1pt]{$0$}--++(0,-1)node[below]{$\phantom{i}$};
        \end{tikzpicture}
        $\cdots$
        \begin{tikzpicture}[anchorbase,smallnodes,rounded corners]
            \draw[solid](-0.2,1)node[above,yshift=-1pt]{$\phantom{i}$}--++(0,-1)node[below]{$1$};
            \draw[ghost](0.4,1)node[above,yshift=-1pt]{$0$}--++(0,-1)node[below]{$\phantom{i}$};
            \draw[solid](1,1)node[above,yshift=-1pt]{$\phantom{i}$}--++(0,-1)node[below]{$1$};
            \draw[ghost](1.6,1)node[above,yshift=-1pt]{$0$}--++(0,-1)node[below]{$\phantom{i}$};
        \end{tikzpicture}
        \end{center}
        \item (solid $1$)(ghost $0$)(solid $1$)$\cdots$(solid $1$)(ghost $0$), i.e.\\
        \begin{center}
        \begin{tikzpicture}[anchorbase,smallnodes,rounded corners]
            \draw[solid](1,1)node[above,yshift=-1pt]{$\phantom{i}$}--++(0,-1)node[below]{$1$};
            \draw[ghost](1.6,1)node[above,yshift=-1pt]{$0$}--++(0,-1)node[below]{$\phantom{i}$};
            \draw[solid](2.2,1)node[above,yshift=-1pt]{$\phantom{i}$}--++(0,-1)node[below]{$1$};
        \end{tikzpicture}
        $\cdots$
        \begin{tikzpicture}[anchorbase,smallnodes,rounded corners]
            \draw[solid](1,1)node[above,yshift=-1pt]{$\phantom{i}$}--++(0,-1)node[below]{$1$};
            \draw[ghost](1.6,1)node[above,yshift=-1pt]{$0$}--++(0,-1)node[below]{$\phantom{i}$};
        \end{tikzpicture}
        \end{center}
        \item (solid $1$)(ghost $0$)(solid $1$)$\cdots$(solid $1$)(ghost $0$)(solid $1$), i.e.\\
        \begin{center}
        \begin{tikzpicture}[anchorbase,smallnodes,rounded corners]
            \draw[solid](1,1)node[above,yshift=-1pt]{$\phantom{i}$}--++(0,-1)node[below]{$1$};
            \draw[ghost](1.6,1)node[above,yshift=-1pt]{$0$}--++(0,-1)node[below]{$\phantom{i}$};
            \draw[solid](2.2,1)node[above,yshift=-1pt]{$\phantom{i}$}--++(0,-1)node[below]{$1$};
        \end{tikzpicture}
        $\cdots$
        \begin{tikzpicture}[anchorbase,smallnodes,rounded corners]
            \draw[solid](1,1)node[above,yshift=-1pt]{$\phantom{i}$}--++(0,-1)node[below]{$1$};
            \draw[ghost](1.6,1)node[above,yshift=-1pt]{$0$}--++(0,-1)node[below]{$\phantom{i}$};
            \draw[solid](2.2,1)node[above,yshift=-1pt]{$\phantom{i}$}--++(0,-1)node[below]{$1$};
        \end{tikzpicture}
        \end{center}
        \item (ghost $0$)(solid $1$)$\cdots$(solid $1$)(ghost $0$)(solid $1$), i.e.\\
        \begin{center}
        \begin{tikzpicture}[anchorbase,smallnodes,rounded corners]
            \draw[ghost](1.6,1)node[above,yshift=-1pt]{$0$}--++(0,-1)node[below]{$\phantom{i}$};
            \draw[solid](2.2,1)node[above,yshift=-1pt]{$\phantom{i}$}--++(0,-1)node[below]{$1$};
        \end{tikzpicture}
        $\cdots$
        \begin{tikzpicture}[anchorbase,smallnodes,rounded corners]
            \draw[solid](1,1)node[above,yshift=-1pt]{$\phantom{i}$}--++(0,-1)node[below]{$1$};
            \draw[ghost](1.6,1)node[above,yshift=-1pt]{$0$}--++(0,-1)node[below]{$\phantom{i}$};
            \draw[solid](2.2,1)node[above,yshift=-1pt]{$\phantom{i}$}--++(0,-1)node[below]{$1$};
        \end{tikzpicture}
        \end{center}
    \end{enumerate}
    We call it a \textbf{close} $\boldsymbol{(0,1)}$-\textbf{tuple}. It is maximal if it can not be expanded, or there are no other ghost $0$-strings or solid $1$-string close to any string of tuple. The starting residue is the \textbf{rightmost} residue, while the ending residue is the \textbf{leftmost} residue. 
\end{Definition}
\begin{Remark}
    The definition of starting and ending residues may seem peculiar, but it’s chosen for the sake of consistency with later results. Alternatively, recall our 'pulling to the right' strategy for constructing the idempotent from a tableau: we place solid strings to the leftmost position and then pull them as far right as possible. Therefore, the right strings should correspond to earlier nodes.
\end{Remark}
\begin{Definition}\label{subdivisionidempotent}
     The $(t,\epsilon)$-subdivision map $S_{t,\epsilon}$ of an idempotent string diagram is defined as follows: for any (steady) idempotent diagram $Y$, first replace all $i$ by $i+1$, where $i\in\{2,\cdots,e\}$, then consider all maximal close $(0,1)$-tuple:
     \begin{enumerate}
         \item If it starts with $0$, ends with $0$, then relabel $1$ by $2$, add a solid $1$-string on the left of each ghost $0$-string by $t$-units, a ghost $1$-string on the right of the added solid $1$-string by $\epsilon$ units, i.e.\\
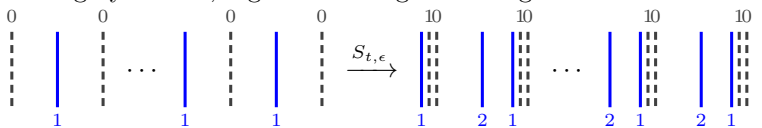
        \begin{tikzpicture}[anchorbase,smallnodes,rounded corners]
            \draw[solid](1,1)node[above,yshift=-1pt]{$\phantom{i}$}--++(0,-1)node[below]{$1$};
            \draw[ghost](0.4,1)node[above,yshift=-1pt]{$0$}--++(0,-1)node[below]{$\phantom{i}$};
            \draw[ghost](1.6,1)node[above,yshift=-1pt]{$0$}--++(0,-1)node[below]{$\phantom{i}$};
        \end{tikzpicture}
        $\cdots$
        \begin{tikzpicture}[anchorbase,smallnodes,rounded corners]
            \draw[solid](-0.2,1)node[above,yshift=-1pt]{$\phantom{i}$}--++(0,-1)node[below]{$1$};
            \draw[ghost](0.4,1)node[above,yshift=-1pt]{$0$}--++(0,-1)node[below]{$\phantom{i}$};
            \draw[solid](1,1)node[above,yshift=-1pt]{$\phantom{i}$}--++(0,-1)node[below]{$1$};
            \draw[ghost](1.6,1)node[above,yshift=-1pt]{$0$}--++(0,-1)node[below]{$\phantom{i}$};
        \end{tikzpicture}
        $\xrightarrow[]{S_{t,\epsilon}}$
        \begin{tikzpicture}[anchorbase,smallnodes,rounded corners]
            \draw[solid](0.2,1)node[above,yshift=-1pt]{$\phantom{i}$}--++(0,-1)node[below]{$1$};
            \draw[ghost](0.3,1)node[above,yshift=-1pt]{$1$}--++(0,-1)node[below]{$\phantom{i}$};
            \draw[solid](1.4,1)node[above,yshift=-1pt]{$\phantom{i}$}--++(0,-1)node[below]{$1$};
            \draw[ghost](1.5,1)node[above,yshift=-1pt]{$1$}--++(0,-1)node[below]{$\phantom{i}$};
            \draw[solid](1,1)node[above,yshift=-1pt]{$\phantom{i}$}--++(0,-1)node[below]{$2$};
            \draw[ghost](0.4,1)node[above,yshift=-1pt]{$0$}--++(0,-1)node[below]{$\phantom{i}$};
            \draw[ghost](1.6,1)node[above,yshift=-1pt]{$0$}--++(0,-1)node[below]{$\phantom{i}$};
        \end{tikzpicture}
        $\cdots$
        \begin{tikzpicture}[anchorbase,smallnodes,rounded corners]              
            \draw[solid](0.2,1)node[above,yshift=-1pt]{$\phantom{i}$}--++(0,-1)node[below]{$1$};
            \draw[ghost](0.3,1)node[above,yshift=-1pt]{$1$}--++(0,-1)node[below]{$\phantom{i}$};
            \draw[solid](-0.2,1)node[above,yshift=-1pt]{$\phantom{i}$}--++(0,-1)node[below]{$2$};
            \draw[ghost](0.4,1)node[above,yshift=-1pt]{$0$}--++(0,-1)node[below]{$\phantom{i}$};
            \draw[solid](1.4,1)node[above,yshift=-1pt]{$\phantom{i}$}--++(0,-1)node[below]{$1$};
            \draw[ghost](1.5,1)node[above,yshift=-1pt]{$1$}--++(0,-1)node[below]{$\phantom{i}$};
            \draw[solid](1,1)node[above,yshift=-1pt]{$\phantom{i}$}--++(0,-1)node[below]{$2$};
            \draw[ghost](1.6,1)node[above,yshift=-1pt]{$0$}--++(0,-1)node[below]{$\phantom{i}$};
        \end{tikzpicture}.
        \item If it starts with $0$, ends with $1$, then relabel $1$ by $2$, add a solid $1$-string on the left of each ghost $0$-string by $t$-units, a ghost $1$-string on the right of the added solid $1$-string by $\epsilon$ units, i.e.\\
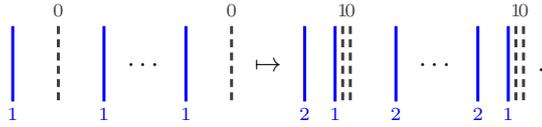
            \begin{tikzpicture}[anchorbase,smallnodes,rounded corners]
                   \draw[solid](1,1)node[above,yshift=-1pt]{$\phantom{i}$}--++(0,-1)node[below]{$1$};
                   \draw[ghost](1.6,1)node[above,yshift=-1pt]{$0$}--++(0,-1)node[below]{$\phantom{i}$};
                   \draw[solid](2.2,1)node[above,yshift=-1pt]{$\phantom{i}$}--++(0,-1)node[below]{$1$};
            \end{tikzpicture}
            $\cdots$
            \begin{tikzpicture}[anchorbase,smallnodes,rounded corners]
                   \draw[solid](1,1)node[above,yshift=-1pt]{$\phantom{i}$}--++(0,-1)node[below]{$1$};
                   \draw[ghost](1.6,1)node[above,yshift=-1pt]{$0$}--++(0,-1)node[below]{$\phantom{i}$};
            \end{tikzpicture}
            $\mapsto$
            \begin{tikzpicture}[anchorbase,smallnodes,rounded corners]
                   \draw[solid](1.4,1)node[above,yshift=-1pt]{$\phantom{i}$}--++(0,-1)node[below]{$1$};
                   \draw[ghost](1.5,1)node[above,yshift=-1pt]{$1$}--++(0,-1)node[below]{$\phantom{i}$};
                   \draw[solid](1,1)node[above,yshift=-1pt]{$\phantom{i}$}--++(0,-1)node[below]{$2$};
                   \draw[ghost](1.6,1)node[above,yshift=-1pt]{$0$}--++(0,-1)node[below]{$\phantom{i}$};
                   \draw[solid](2.2,1)node[above,yshift=-1pt]{$\phantom{i}$}--++(0,-1)node[below]{$2$};
            \end{tikzpicture}
            $\cdots$
            \begin{tikzpicture}[anchorbase,smallnodes,rounded corners]              
                   \draw[solid](1.4,1)node[above,yshift=-1pt]{$\phantom{i}$}--++(0,-1)node[below]{$1$};
                   \draw[ghost](1.5,1)node[above,yshift=-1pt]{$1$}--++(0,-1)node[below]{$\phantom{i}$};
                   \draw[solid](1,1)node[above,yshift=-1pt]{$\phantom{i}$}--++(0,-1)node[below]{$2$};
                   \draw[ghost](1.6,1)node[above,yshift=-1pt]{$0$}--++(0,-1)node[below]{$\phantom{i}$};
            \end{tikzpicture}.
        \item If it is a single solid $1$-string, relabel it by $2$, else if it starts with $1$, ends with $1$, then relabel $1$ by $2$, add a solid $1$-string on the right of each ghost $0$-string by $t$-units, a ghost $1$-string on the right of the added solid $1$-string by $\epsilon$ units, i.e.\\
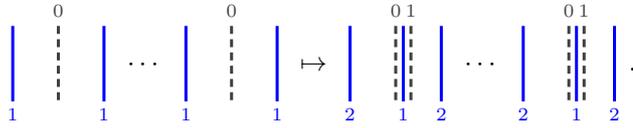
            \begin{tikzpicture}[anchorbase,smallnodes,rounded corners]
                   \draw[solid](1,1)node[above,yshift=-1pt]{$\phantom{i}$}--++(0,-1)node[below]{$1$};
                   \draw[ghost](1.6,1)node[above,yshift=-1pt]{$0$}--++(0,-1)node[below]{$\phantom{i}$};
                   \draw[solid](2.2,1)node[above,yshift=-1pt]{$\phantom{i}$}--++(0,-1)node[below]{$1$};
            \end{tikzpicture}
            $\cdots$
            \begin{tikzpicture}[anchorbase,smallnodes,rounded corners]
                   \draw[solid](1,1)node[above,yshift=-1pt]{$\phantom{i}$}--++(0,-1)node[below]{$1$};
                   \draw[ghost](1.6,1)node[above,yshift=-1pt]{$0$}--++(0,-1)node[below]{$\phantom{i}$};
                    \draw[solid](2.2,1)node[above,yshift=-1pt]{$\phantom{i}$}--++(0,-1)node[below]{$1$};
            \end{tikzpicture}
            $\mapsto$
            \begin{tikzpicture}[anchorbase,smallnodes,rounded corners]
                   \draw[solid](1.7,1)node[above,yshift=-1pt]{$\phantom{i}$}--++(0,-1)node[below]{$1$};
                   \draw[ghost](1.8,1)node[above,yshift=-1pt]{$1$}--++(0,-1)node[below]{$\phantom{i}$};
                   \draw[solid](1,1)node[above,yshift=-1pt]{$\phantom{i}$}--++(0,-1)node[below]{$2$};
                   \draw[ghost](1.6,1)node[above,yshift=-1pt]{$0$}--++(0,-1)node[below]{$\phantom{i}$};
                   \draw[solid](2.2,1)node[above,yshift=-1pt]{$\phantom{i}$}--++(0,-1)node[below]{$2$};
            \end{tikzpicture}
            $\cdots$
            \begin{tikzpicture}[anchorbase,smallnodes,rounded corners]              
                   \draw[solid](1.7,1)node[above,yshift=-1pt]{$\phantom{i}$}--++(0,-1)node[below]{$1$};
                   \draw[ghost](1.8,1)node[above,yshift=-1pt]{$1$}--++(0,-1)node[below]{$\phantom{i}$};
                   \draw[solid](1,1)node[above,yshift=-1pt]{$\phantom{i}$}--++(0,-1)node[below]{$2$};
                   \draw[ghost](1.6,1)node[above,yshift=-1pt]{$0$}--++(0,-1)node[below]{$\phantom{i}$};
                    \draw[solid](2.2,1)node[above,yshift=-1pt]{$\phantom{i}$}--++(0,-1)node[below]{$2$};
            \end{tikzpicture}.
        \item If it starts with $1$, ends with $0$, then relabel $1$ by $2$, add a solid $1$-string on the right of each ghost $0$-string by $t$-units, a ghost $1$-string on the right of the added solid $1$-string by $\epsilon$ units, i.e.\\
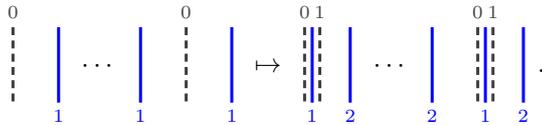
            \begin{tikzpicture}[anchorbase,smallnodes,rounded corners]
                   \draw[ghost](1.6,1)node[above,yshift=-1pt]{$0$}--++(0,-1)node[below]{$\phantom{i}$};
                   \draw[solid](2.2,1)node[above,yshift=-1pt]{$\phantom{i}$}--++(0,-1)node[below]{$1$};
            \end{tikzpicture}
            $\cdots$
            \begin{tikzpicture}[anchorbase,smallnodes,rounded corners]
                   \draw[solid](1,1)node[above,yshift=-1pt]{$\phantom{i}$}--++(0,-1)node[below]{$1$};
                   \draw[ghost](1.6,1)node[above,yshift=-1pt]{$0$}--++(0,-1)node[below]{$\phantom{i}$};
                    \draw[solid](2.2,1)node[above,yshift=-1pt]{$\phantom{i}$}--++(0,-1)node[below]{$1$};
            \end{tikzpicture}
            $\mapsto$
            \begin{tikzpicture}[anchorbase,smallnodes,rounded corners]
                   \draw[solid](1.7,1)node[above,yshift=-1pt]{$\phantom{i}$}--++(0,-1)node[below]{$1$};
                   \draw[ghost](1.8,1)node[above,yshift=-1pt]{$1$}--++(0,-1)node[below]{$\phantom{i}$};
                   \draw[ghost](1.6,1)node[above,yshift=-1pt]{$0$}--++(0,-1)node[below]{$\phantom{i}$};
                   \draw[solid](2.2,1)node[above,yshift=-1pt]{$\phantom{i}$}--++(0,-1)node[below]{$2$};
            \end{tikzpicture}
            $\cdots$
            \begin{tikzpicture}[anchorbase,smallnodes,rounded corners]              
                   \draw[solid](1.7,1)node[above,yshift=-1pt]{$\phantom{i}$}--++(0,-1)node[below]{$1$};
                   \draw[ghost](1.8,1)node[above,yshift=-1pt]{$1$}--++(0,-1)node[below]{$\phantom{i}$};
                   \draw[solid](1,1)node[above,yshift=-1pt]{$\phantom{i}$}--++(0,-1)node[below]{$2$};
                   \draw[ghost](1.6,1)node[above,yshift=-1pt]{$0$}--++(0,-1)node[below]{$\phantom{i}$};
                   \draw[solid](2.2,1)node[above,yshift=-1pt]{$\phantom{i}$}--++(0,-1)node[below]{$2$};
            \end{tikzpicture}.
     \end{enumerate}
\end{Definition}
\begin{Remark}
    \begin{itemize}
        \item The strategies of (a) and (b) are essentially the same, while the strategies of (c) and (d) are essentially the same. Thus, \textbf{the starting residue becomes the distinguished object.} We divide it into four cases, coinciding with \autoref{lambdaplus}.
        \item It is easy to see that both \autoref{closetuple} and \autoref{subdivisionidempotent} can be generalized to any $(i,i+1)$ case similarly. We will need this to define the general subdivision in \autoref{generalcase}.
    \end{itemize}
\end{Remark}
Let's extend the subdivision map to the entire algebra. To achieve this, we only need to define it for a basis of the KLRW algebra.
According to \autoref{standardbasis}, an element in a standard basis $\WABasis$ is of the form $\1_{\by,\bj}D(w)y_{1}^{a_{1}}\dots y_{n}^{a_{n}}\1_{\bx,\bi}$. We have already defined the subdivision map on the bottom part $\1_{\bx,\bi}$, while we leave the dots unchanged and do not add any dots to the newly added strings. Now, we need to define the map on the permutation part. We extend the added strings along the corresponding ghost $0$-strings.
\begin{Example} Here we revise the example \cite[Example 5A.22]{MaTu-klrw-algebras} to explain the effect of subdivision on $D(w)$. Let $D(w)$ be the following string diagram:
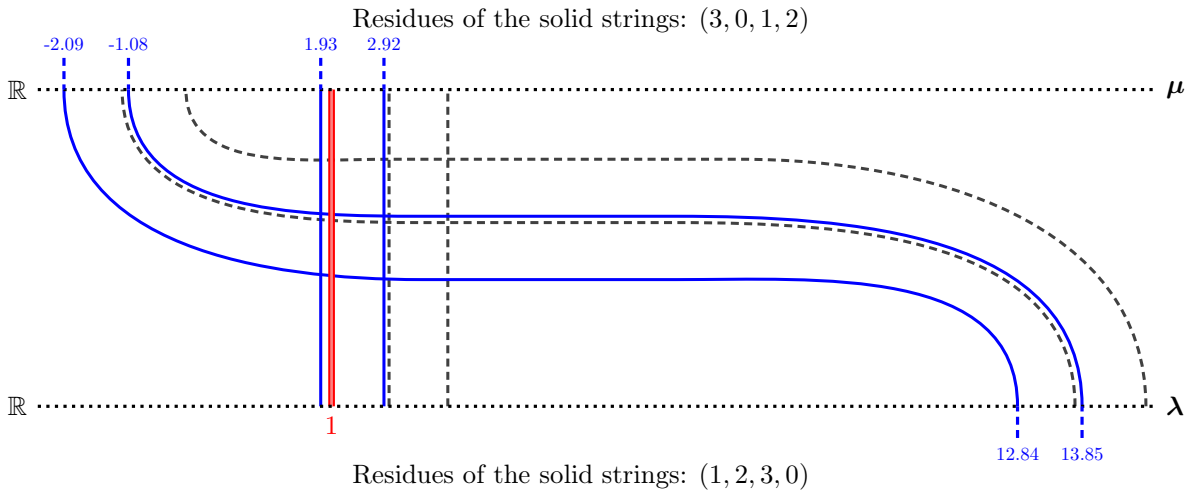
\begin{center}
    \begin{tikzpicture}[scale=0.84,anchorbase]
        \draw[ghost](3,0)to[out=90,in=270](3,5);
        \draw[ghost](3.92,0)to[out=90,in=270](3.92,5);
        \draw[ghost](13.74,0)to[out=90,in=0](7.5,2.9)to(3.5,2.9)to[out=180,in=270](-1.18,5);
        \draw[ghost](14.85,0)to[out=90,in=0](8.5,3.9)to(3.0,3.9)to[out=180,in=270](-0.18,5);
        \draw[solid](1.93,0)to[out=90,in=270](1.93,5);
        \draw[solid](2.92,0)to[out=90,in=270](2.92,5);
        \draw[solid](12.84,0)to[out=90,in=0](7.5,2)to(3.5,2)to[out=180,in=270](-2.09,5);
        \draw[solid](13.85,0)to[out=90,in=0](7.5,3)to(3.5,3)to[out=180,in=270](-1.08,5);
        \draw[redstring](2.1,5)node[above,yshift=-1pt]{$\phantom{i}$}to(2.1,0)node[below]{$1$};
        \draw[very thick,dotted] (-2.5,0)node[left]{$\R$} to (15,0)node[right]{$\blam$};
        \draw[very thick,dotted] (-2.5,5)node[left]{$\R$} to (15,5)node[right]{$\bmu$};
        \draw[very thick,blue,densely dashed] (-2.09,5.5)node[above,yshift=-1pt]{\scalebox{0.7}{$\text{-}2.09$}} to (-2.09,5);
        \draw[very thick,blue,densely dashed] (-1.08,5.5)node[above,yshift=-1pt]{\scalebox{0.7}{$\text{-}1.08$}} to (-1.08,5);
        \draw[very thick,blue,densely dashed] (1.93,5.5)node[above,yshift=-1pt]{\scalebox{0.7}{$1.93$}} to (1.93,5);
        \draw[very thick,blue,densely dashed] (2.92,5.5)node[above,yshift=-1pt]{\scalebox{0.7}{$2.92$}} to (2.92,5);
        \draw[very thick,blue,densely dashed] (12.84,-0.5)node[below]{\scalebox{0.7}{$12.84$}} to (12.84,0);
        \draw[very thick,blue,densely dashed] (13.85,-0.5)node[below]{\scalebox{0.7}{$13.85$}} to (13.85,0);
        \node at (6,-1.1){Residues of the solid strings: $(1,2,3,0)$};
        \node at (6,6.1){Residues of the solid strings: $(3,0,1,2)$};
    \end{tikzpicture}
\end{center}
Then after subdivision (and relabelling as we like to do), the image is as follows:
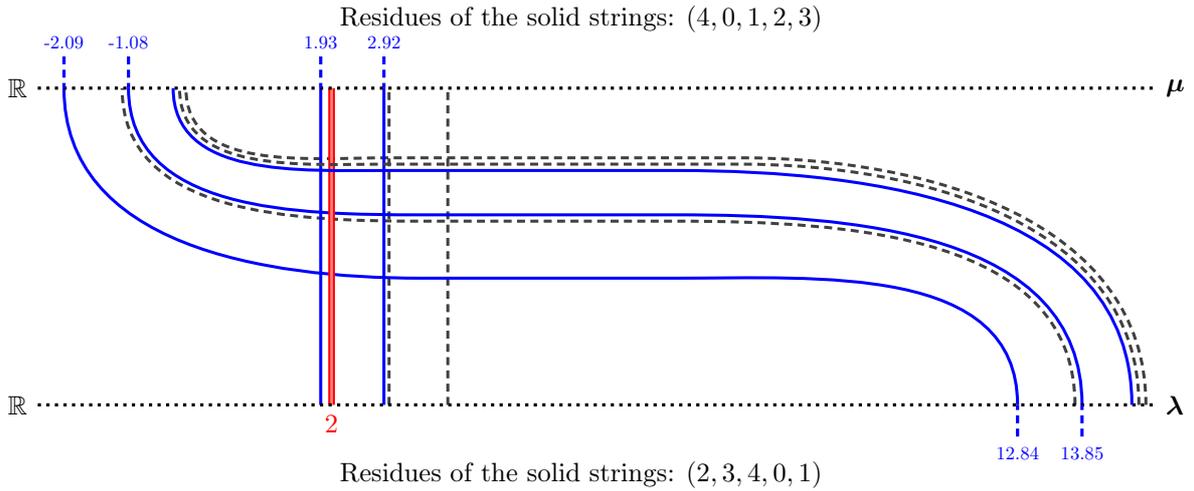
\begin{center}
    \begin{tikzpicture}[scale=0.84,anchorbase]
        \draw[ghost](3,0)to[out=90,in=270](3,5);
        \draw[ghost](3.92,0)to[out=90,in=270](3.92,5);
        \draw[ghost](13.74,0)to[out=90,in=0](7.5,2.9)to(3.5,2.9)to[out=180,in=270](-1.18,5);
        \draw[ghost](14.85,0)to[out=90,in=0](8.5,3.9)to(3.0,3.9)to[out=180,in=270](-0.18,5);
        \draw[solid](1.93,0)to[out=90,in=270](1.93,5);
        \draw[solid](2.92,0)to[out=90,in=270](2.92,5);
        \draw[solid](12.84,0)to[out=90,in=0](7.5,2)to(3.5,2)to[out=180,in=270](-2.09,5);
        \draw[solid](13.85,0)to[out=90,in=0](7.5,3)to(3.5,3)to[out=180,in=270](-1.08,5);
        \draw[redstring](2.1,5)node[above,yshift=-1pt]{$\phantom{i}$}to(2.1,0)node[below]{$2$};
        \draw[very thick,dotted] (-2.5,0)node[left]{$\R$} to (15,0)node[right]{$\blam$};
        \draw[very thick,dotted] (-2.5,5)node[left]{$\R$} to (15,5)node[right]{$\bmu$};
        \draw[very thick,blue,densely dashed] (-2.09,5.5)node[above,yshift=-1pt]{\scalebox{0.7}{$\text{-}2.09$}} to (-2.09,5);
        \draw[very thick,blue,densely dashed] (-1.08,5.5)node[above,yshift=-1pt]{\scalebox{0.7}{$\text{-}1.08$}} to (-1.08,5);
        \draw[very thick,blue,densely dashed] (1.93,5.5)node[above,yshift=-1pt]{\scalebox{0.7}{$1.93$}} to (1.93,5);
        \draw[very thick,blue,densely dashed] (2.92,5.5)node[above,yshift=-1pt]{\scalebox{0.7}{$2.92$}} to (2.92,5);
        \draw[very thick,blue,densely dashed] (12.84,-0.5)node[below]{\scalebox{0.7}{$12.84$}} to (12.84,0);
        \draw[very thick,blue,densely dashed] (13.85,-0.5)node[below]{\scalebox{0.7}{$13.85$}} to (13.85,0);
        \draw[solid](14.63,0)to[out=90,in=0](7.5,3.7)to(3.0,3.7)to[out=180,in=270](-0.38,5);
        \draw[ghost](14.74,0)to[out=90,in=0](8.5,3.8)to(3.0,3.8)to[out=180,in=270](-0.28,5);
        \node at (6,-1.1){Residues of the solid strings: $(2,3,4,0,1)$};
        \node at (6,6.1){Residues of the solid strings: $(4,0,1,2,3)$};
    \end{tikzpicture}
\end{center}
\end{Example}
The remaining part of this section follows immediately from \cite[Section 4]{MaTu-klrw-algebras}.
\begin{Definition}
    A subdivision is called \textbf{homogeneous} if $0<t\ll 1$ and $0<\epsilon\ll 1$.
\end{Definition}
\begin{Proposition}[{Mathas-Tubbenhauer~\cite{MaTu-klrw-algebras}}]
    Suppose the subdivision is homogeneous, then it preservers the degree of a string diagram, i.e. $\deg (D)=\deg S_{t,\epsilon}(D)$.
\end{Proposition}
\begin{proof}
    See \cite[Lemma 4C.9]{MaTu-klrw-algebras}.
\end{proof}
\begin{Definition}\label{D:SubWebster}
If $\overline{\Gamma}$ is a subdivision of $\Gamma$, then define
$\Oneg=\sum_{\bx\in X}\sum_{\bi\in I^\beta}\1_{\Sub{\bx},\overline{\bi}}$ and
\begin{gather}\label{E:IdemTrunc}
  \SA(X)=\Oneg\WAs(\Sub{X})\Oneg/
\Oneg\WAs(\overline{X})\1_{\text{bad}}\WAs(\overline{X})\Oneg.
\end{gather}
\end{Definition}

We identify a diagram in $\WA(\overline{X})$ with its image in $\SA(X)$.

\begin{Lemma}\label{L:KillDiagrams}
  Let $B\in\SA(X)$ be a bad diagram. Then $B=0$.
\end{Lemma}

\begin{proof}
See \cite[4C.14]{MaTu-klrw-algebras}.
\end{proof}

\begin{Theorem}[{Mathas-Tubbenhauer~\cite{MaTu-klrw-algebras}}]\label{T:SubDiv}
    Suppose the subdivision is homogeneous, then there is an isomorphism of graded algebras $S_{t,\epsilon }:\map{\WA(X)}{\SA(X)}$.
\end{Theorem}
\begin{proof}
    See \cite[Theorem 4D.2]{MaTu-klrw-algebras}. We remind the reader there are two parts of the proof. The first is to check it is a homomorphism of KLRW algebras, that is, the image also satisfies the KLRW relations. This can be verified locally and the essential condition is the homogeneity. The second part is to verify the bijection, the only difference is that this time we make use of both relations of \autoref{R:GhostSolid} while in \cite[Theorem 4D.2]{MaTu-klrw-algebras} only the left one is used.
\end{proof}
\begin{Corollary}
    The isomorphism $S_{t,\epsilon }:\map{\WA(X)}{\SA(X)}$ of \autoref{T:SubDiv} descends to an isomorphism of the cyclotomic quotients.
\end{Corollary}
\begin{proof}
    This follows from the fact being unsteady or not remain unchanged after the subdivision.
\end{proof}

\section{Equality of Decomposition Numbers}\label{equalityofdecompositionnumber}
Now that we have defined the subdivision map, an important question arises: what is the image of $1_\blam$ under $S_{t,\epsilon}$? It should be an idempotent diagram again, but does it appear in the cellular basis of the new KLRW algebra? This question is equivalent to asking whether there exists an $\ell$-partition $\blam^+$, such that $1_{\blam^+}=S_{t,\epsilon}(1_\blam)$. We provide an affirmative answer in this section. Let's first define the partition $\blam^+$. Indeed, we will provide two equivalent definitions.
\subsection{Box Definition}
In this section, we provide several definitions for the level $1$ case, i.e., for $1$-partitions. it is straightforward to generalize them to higher levels.
\begin{Definition}\label{strip}
    For $\lambda$ a partition of $n$, fix a charge $\rho\in I$, consider its box configuration filled with residues (with regard to $\rho$), denoted by $[\lambda]_\rho$. A $\boldsymbol{(0,1)}$-\textbf{strip} is a sequence of nodes $(A_1,A_2,\cdots,A_k)$where $A_i\in [\lambda]$ such that
    \begin{enumerate}
        \item $\res(A_i)=0$;
        \item $A_i$ and $A_{i+1}$ are (increasingly) adjacent, that is, assume $A_i=(m,r,c)$, then
        either $A_{i+1}=(m,r,c+1)$, or $A_{i+1}=(m,r+1,c)$; 
    \end{enumerate}
    A $(0,1)$-strip is \textbf{maximal} if the strip can't be completed into a larger one. Similarly, we define the (maximal) s$(1,0)$-strips.
\end{Definition}
\begin{Remark}
    \begin{itemize}
        \item Clearly, there are two types of $(0,1)$-strips and two types of $(1,0)$-strips, depending on the residue of the end node.
        \item We can similarly define (maximal) $(i,i+1)$-strips, but we will not use this until \autoref{generalcase}.
    \end{itemize}
\end{Remark}
\begin{Definition}\label{lambdaplus}
    For $\lambda$ a partition of $n$ and $\rho\in I$, take its box configuration filled with residues.  Replace all residue $i$ by $i+1$ where $i\in\{2,\cdots,e\}$. Then consider all maximal $(0,1)$-strips and maximal $(1,0)$-strips, step by step:
    \begin{enumerate}
        \item If the starting residue is $0$ and the ending residue is $0$, that is:\\
        \Tableau[scale=0.5]{{0,1,\\,\\,\\},{\\,0,1,\\,\\},{\\,\\,\cdots,\\,\\}, {\\,\\,\\,0,1},{\\,\\,\\,\\,0}}, replace it by \Tableau[scale=0.5]{{0,1,2,\\,\\,\\},{\\,0,1,2,\\,\\},{\\,\\,\cdots,\\,\\,\\}, {\\,\\,\\,0,1,2},{\\,\\,\\,\\,0,1}}.\\
        If it is single $0$-node, replace it by \Tableau[scale=0.5]{{0,1}}.
        \item If the starting residue is $0$ and the ending residue is $1$, that is:\\
        \Tableau[scale=0.5]{{0,1,\\,\\,\\},{\\,0,1,\\,\\},{\\,\\,\cdots,\\,\\}, {\\,\\,\\,0,1}}, replace it by \Tableau[scale=0.5]{{0,1,2,\\,\\,\\},{\\,0,1,2,\\,\\},{\\,\\,\cdots,\\,\\,\\}, {\\,\\,\\,0,1,2}}
        \item If the starting residue is $1$ and the ending residue is $1$, that is:\\
        \Tableau[scale=0.5]{{1,\\,\\,\\,\\},{0,1,\\,\\,\\},{\\,\cdots,\cdots,\\,\\}, {\\,\\,0,1,\\},{\\,\\,\\,0,1}}, replace it by \Tableau[scale=0.5]{{2,\\,\\,\\,\\},{1,2,\\,\\,\\},{0,1,\cdots,\\,\\}, {\\,\cdots,\cdots,\cdots,\\},{\\,\\,0,1,2},{\\,\\,\\,0}}. \\
        If there is just a single $1$-node, replace it by \Tableau[scale=0.5]{{2}}.
        \item If the starting residue is $1$ and the ending residue is $0$, that is:\\
        \Tableau[scale=0.5]{{1,\\,\\,\\},{0,1,\\,\\},{\\,\cdots,\cdots,\\},{\\,\\,0,1},{\\,\\,\\,0}}, replace it by \Tableau[scale=0.5]{{2,\\,\\,\\},{1,2,\\,\\},{0,\cdots,\cdots,\\},{\\,\cdots,1,2},{\\,\\,0,1},{\\,\\,\\,0}}     
    \end{enumerate}
    Hence, we will obtain a new box configuration (filled with residues), and let $\blam^+$ be the composition corresponding to this box configuration. The cases where $[\lambda]_\rho$ consists of one $0$-node or one $1$-node are called \textbf{trivial} cases.
\end{Definition}
\begin{Proposition}\label{welldefineoflambdaplus}
    $\blam^+$ is a partition.
\end{Proposition}
\begin{proof}
   Only need to show that for (a), (b), (c), and (d), every step of replacement still yields a partition. This is straightforward. Refer to \autoref{anotherproof} for another (seemingly more rigorous) proof.
\end{proof}
\begin{Remark}
    The definition of $\blam^+$ clearly depends on the residue pattern, hence it relies on the choice of $\brho$. This aligns with the notion of the subdivision map.
\end{Remark}
\begin{Example}
Assume $\ell=1$ and $e=3$. Fill the young diagram by residues:
    \begin{enumerate}
        \item If $\blam=(1)$ and $\rho=0$, then $\blam^+=(2)$.\\
                $[\blam]:$ \Tableau[scale=0.5]{{0}}, $[\blam^+]$: \Tableau[scale=0.5]{{0,1}}
        \item If $\blam=(1^2)$ and $\rho=1$, then $\blam^+=(1^3)$.\\
                $[\blam]:$ \Tableau[scale=0.5]{{1},{0}}, $[\blam^+]$: \Tableau[scale=0.5]{{2},{1},{0}}
        \item If $\blam=(2^2)$ and $\rho=0$, then $\blam^+=(3^2)$.\\
                $[\blam]:$ \Tableau[scale=0.5]{{0,1},{3,0}}, $[\blam^+]$: \Tableau[scale=0.5]{{0,1,2},{3,0,1}}
        \item If $\blam=(2,1)$ and $\rho=1$, then $\blam^+=(2,1^2)$.\\
                $[\blam]:$ \Tableau[scale=0.5]{{1,2},{0}}, $[\blam^+]$: \Tableau[scale=0.5]{{2,3},{1},{0}}
        \item If $\blam=(2,2)$ and $\rho=1$, then $\blam^+=(2^2,1)$.\\
                $[\blam]:$ \Tableau[scale=0.5]{{1,2},{0,1}}, $[\blam^+]$: \Tableau[scale=0.5]{{2,3},{1,2},{0}}
        \item If $\blam=(2,2,2)$ and $\rho=1$, then $\blam^+=(2^4)$.\\
                $[\blam]:$ \Tableau[scale=0.5]{{1,2},{0,1},{3,0}}, $[\blam^+]$: \Tableau[scale=0.5]{{2,3},{1,2},{0,1},{4,0}}
        \item\label{emptyrunnerexample} If $\blam=(4,3,2)$ and $\rho=3$, then $\blam^+=(5,4,2)$.\\
                $[\blam]:$ \Tableau[scale=0.5]{{3,0,1,2},{2,3,0},{1,2}}, $[\blam^+]$: \Tableau[scale=0.5]{{4,0,1,2,3},{3,4,0,1},{2,3}}
        \item\label{fullrunnerexample}  If $\blam=(4,4,4,4)$ and $\rho=1$, then $\blam^+=(5,4,4,4,3)$.\\
                $[\blam]:$ \Tableau[scale=0.5]{{1,2,3,0},{0,1,2,3},{3,0,1,2},{2,3,0,1}}, $[\blam^+]$: \Tableau[scale=0.5]{{2,3,4,0,1},{1,2,3,4},{0,1,2,3},{4,0,1,2},{3,4,0}}
        \item  If $\blam=(5,3,3,2,1)$ and $\rho=2$, then $\blam^+=(6,3,3,2,2,1)$.\\
                $[\blam]:$ \Tableau[scale=0.5]{{2,3,0,1,2},{1,2,3},{0,1,2},{3,0},{2}}, $[\blam^+]$: \Tableau[scale=0.5]{{3,4,0,1,2,3},{2,3,4},{1,2,3},{0,1},{4,0},{3}}
        \item  If $\blam=(5,3,3,2,1)$ and $\rho=3$, then $\blam^+=(6,4,3,2,1,1)$.\\
                $[\blam]:$ \Tableau[scale=0.5]{{3,0,1,2,3},{2,3,0},{1,2,3},{0,1},{3}}, $[\blam^+]$: \Tableau[scale=0.5]{{4,0,1,2,3,4},{3,4,0,1},{2,3,4},{1,2},{0},{4}}
        \item If $\blam=(8,7,5,5,4,3,2,2,1)$ and $\rho=0$, then $\blam^+=(10,9,6,6,4,3,2,2,2,1,1)$.\\
                $[\blam]:$ \Tableau[scale=0.5]{{0,1,2,3,0,1,2,3},{3,0,1,2,3,0,1},{2,3,0,1,2},{1,2,3,0,1},{0,1,2,3},{3,0,1},{2,3},{1,2},{0}}, $[\blam^+]$: \Tableau[scale=0.5]{{0,1,2,3,4,0,1,2,3,4},{4,0,1,2,3,4,0,1,2},{3,4,0,1,2,3},{2,3,4,0,1,2},{1,2,3,4},{0,1,2},{4,0},{3,4},{2,3},{1},{0}}
        \item If $\blam=(1^9)$ and $\rho=0$, then $\blam^+=(2,1^{10})$.\\
                $[\blam]:$ \Tableau[scale=0.5]{{0},{3},{2},{1},{0},{3},{2},{1},{0}}, $[\blam^+]$: \Tableau[scale=0.5]{{0,1},{4},{3},{2},{1},{0},{4},{3},{2},{1},{0}}
        
    \end{enumerate}
\end{Example}
We have some easy facts:
\begin{Lemma}\label{boxlemma}
    For $\lambda$ a partition of $n$ and $\rho\in I$ fixed, form the box configuration of $\lambda$ filled with residues with regard to $\rho$, it defines a $[\lambda^+]$:
    \begin{enumerate}
        \item Let $k(\lambda)$ be the number of maximal $(1,0)$-strips in $[\lambda]$, then $\ell(\lambda^+)=\ell(\lambda)+k(\lambda)$.
        \item There is a natural bijection between maximal $(0,1)$-strips in $[\lambda]$ and maximal $(0,1)$-strips in $[\lambda^+]$. 
        \item All maximal $(1,0)$-strips, the starting node is in the first column.
        \item All maximal $(0,1)$-strips, the starting node is in the first row.
    \end{enumerate}
\end{Lemma}
\begin{proof}
    Obvious from construction.
\end{proof}
\subsection{Abacus Definition}
There is another definition of $\blam^+$, using tthe $e'$-abacus configurations, where $e'=e+1$. 
\begin{Remark}
    It is common to use $A^{(1)}_{e}$ to denote the quiver with $e$ vertices, but in \cite{MaTu-klrw-algebras}, they use $A^{(1)}_{e-1}$ to denote the same quiver. We adopt the notations from \cite{MaTu-klrw-algebras}, so the number of vertices in the quiver of type $A^{(1)}_{e}$ is equal to $e'=e+1$. Hope this will not cause any confusion.
\end{Remark}

Fix a partition $\lambda$ of $n$ and $\rho\in I$, we identify $I$ to a representative set $\{0,1,\cdots,e\}$ of $\mathbb{Z}/(e+1)\mathbb{Z}$, so $\rho$ can be regarded as an element in $\Z$.

The beta numbers of $\lambda$ corresponds to $\rho$ are defined to be an infinite sequence: $\beta_i:=\lambda_i+\rho-i$ where $i=1,2,\cdots$. Denote this sequence by $\beta_\rho(\lambda)$.

Take $N$ a negative integer, such that $\rho-Ne'\ge\ell(\lambda)$. Now we can form $e'$ runners, labeled by $\{0,1,\cdots,e\}$ from left to right. For each runner, it is labelled by integers increasing downwards, called \textbf{level}. In each such position, we can put a bead. Assume a position is on $i$-runner at level $j$, we say it is in $je'+i$ position.

For each $i$, if $\beta_i=a_ie'+b_i$, put a bead to $b_i$-runner at level $a_i$ . Stop just before $a_i<N$. Now we get an abacus configuration for $\lambda$. Denote it by $\Ab(\lambda,\rho)$.
\begin{Remark}
    \begin{enumerate}[label=(\alph*)]
        \item Given a sequence of beta numbers, we can recover the partition and charge. Therefore, providing an abacus configuration is equivalent to providing a partition $\lambda$ and a charge $\rho$.
        \item Strictly speaking, this abacus configuration is called a truncated abacus configuration. The traditional abacus configuration allows infinitely many beads to extend to the $-\infty$ level. Everything in this paper can be translated using an infinite abacus configuration without difficulty.
        \item It is easy to generalize the above to multipartitions. Denote is by $\Ab(\blam,\brho)$.
        \item We say two abacus configurations are equivalent if they represent the same partition. In particular, if we maintain the same $\rho$ and $\lambda$ but vary $N$, we obtain two equivalent abacus configurations. This is useful in \autoref{addingfromtheleft} and \autoref{generalcase}, where we need to move all beads of an abacus configuration to the left or right. resulting in another equivalent abacus configuration. 
    \end{enumerate}  
\end{Remark}
\begin{Example}
    Assume $\ell=1$ and $e=3$. If $\rho=(4,4,4,4)$ and $\rho=1$, then take $N=-2$ and $\Ab(\lambda,\rho)$ is:
    \begin{equation}
        \abacus(lmmr,bbbb,bnnn,nbbb,bnnn)
    \end{equation}
    If $\blam=(1^3,(3,1))$ and $\brho=(2,1)$, then take $\mathbf{N}=(-1,-1)$ and $\Ab(\blam,\brho)$ is:
    \begin{equation}
        \abacus(lmmr,bbbn,bbbn),\hspace{0.5cm} \abacus(lmmr,bbbn,bnnb)
    \end{equation}
\end{Example}
\begin{Lemma}\label{N0}
    \begin{enumerate}
        \item There exists maximal $N_0$ such that $N_0\geq N$ for all possible $N$.
        \item  For some $M>-N$, every position above level $M$ is occupied by a bead.
        \item There are $\rho-Ne'$ beads in $\Ab(\lambda,\rho)$.
    \end{enumerate}
\end{Lemma}
\begin{proof}
    \begin{enumerate}[label=(\alph*)]
        \item Assume the length of $\lambda$ is $t$, then $\beta_{t+1}=a_te'+b_t$ for some $a_t\in\Z$ and $0\leq b_t<e'$. Define
        \begin{equation}\label{N0def}
            N_0=\begin{cases}
                a_t+1 & \text{if } b_t=e'-1\\
                a_t & \text{else}
            \end{cases}
        \end{equation}
        Then $N_0$ is the maximal element.
        \item This follows immediately from the construction by taking $M=N_0$.
        \item If $N_0=a_t+1$, then the number of beads is equal to $t+(N_0-N)e'=t+(a_t+1)e'-Ne'=t+1+\beta_{t+1}-Ne'=\rho-Ne'$ where the last equality holds because $\beta_{t+1}=\rho-t-1$; If $N_0=a_t$, then the number of beads is equal to $t+(N_0-N)e'+b_t+1=t+a_te'+b_t+1-Ne'=t+1+\beta_{t+1}-Ne'=\rho-Ne'$.
    \end{enumerate}
\end{proof}

For an abacus configuration $\Ab(\lambda,\rho)$, if we can add one runner with some beads on it, the resulting new abacus configuration will correspond to a new partition. 
\begin{Definition}\label{abacusdeflambdaplus}
    For $\lambda$ a partition of $n$ and $\rho\in I$, let $k(\lambda):=\#\{i| i\equiv \rho\mod(e'), \lambda_{i+1}>0\}$. Insert a new runner to the rightmost, then add $k(\lambda)+N_0-N$ beads from the top of this new runner. That is, for the new abacus, the added beads occupy the positions: $N(e'+1)+e',(N+1)(e'+1)+e',\cdots,(N_0+k(\lambda)-1)(e'+1)+e'$.  
    This new abacus configuration corresponds to a new partition and a new charge, denoted by $\lambda_+$ and $\rho_+$, respectively.
\end{Definition}
\begin{Lemma}\label{newrho}
    $\rho_+=\rho+N_0+k(\lambda)$.
\end{Lemma}
\begin{proof}
    When insert from the rightmost as in \ref{abacusdeflambdaplus}, the new abacus has $\rho-Ne'+k(\lambda)+N_0-N$ beads, which should be equal to $\rho_+-N(e'+1)$ by \autoref{N0}. Hence, we have $\rho_+=\rho+N_0+k(\lambda)$. We warn here it is possible $\rho_+<0$, but this is allowed.  
\end{proof}
\begin{Remark}
    The motivation behind defining $k(\lambda)$ is that it equals the number of non-trivial $(0,1)$-strips in the first column of $[\lambda]$, see \autoref{quantresults}. Later, in \autoref{newdeflambdaplus}, we can eliminate the integer $k(\lambda)$. 
\end{Remark}
\begin{Example}\label{abacusexample1}
    \begin{enumerate}
        \item If $\lambda=(4,3,2)$ and $\rho=3$, then $N_0=0$ and we take $N=-1$. Notice that $k(\lambda)=0$. The abacus configuration $\Ab(\lambda,\rho)$ is:
        \begin{equation}
            \abacus(lmmr,bbbb,nnbn,bnbn)
        \end{equation}
        Since $k(\lambda)=0$, we add a runner to the rightmost and put $k(\lambda)+N_0-N=1$ beads from the top on this runner, hence we get:
        \begin{equation}
            \abacus(lmmmr,bbbbo,nnbnn,bnbnn)
        \end{equation}
        As $\rho_+=\rho+N_0+k(\lambda)=3$, we translate the abacus to the partition $\lambda_+=(5,4,2)$. Compare this example with \ref{emptyrunnerexample}.
        \item If $\lambda=(4,4,4,4)$ and $\rho=1$, then $N_0=-1$ and we take $N=-2$. Also notice $k(\lambda)=1$.
        The abacus configuration $\Ab(\lambda,\rho)$ is:
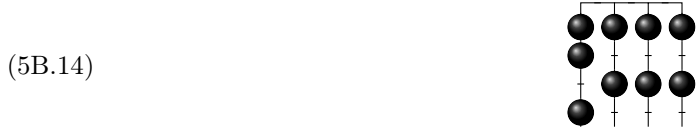
        \begin{equation}
            \abacus(lmmr,bbbb,bnnn,nbbb,bnnn)
        \end{equation}
        Since $k(\lambda)>0$, we add a runner to the rightmost and put $N_0-N+k(\lambda)=2$ beads from the top on this runner, hence we get:
        \begin{equation}
            \abacus(lmmmr,bbbbo,bnnno,nbbbn,bnnnn)
        \end{equation}
        As $\rho_+=\rho+N_0+k(\lambda)=1$, we translate the abacus to the partition $\lambda_+=(5,4,4,4,3)$. Compare this example with \ref{fullrunnerexample}.
    \end{enumerate}
\end{Example}
\subsection{Equivalence of Two Definitions}
Our goal is to show that $\lambda_+=\lambda^+$ if $\rho$ is fixed. Therefore, the two definition\textemdash one using box configuration and the other use abacus configuration---are equivalent. From now on, fix $\rho$ and $\lambda$. We begin with a few lemmas. It is recommended for the reader to first read \autoref{equivdefinition} to avoid confusion.
\begin{Lemma}\label{quantresults}
    \begin{enumerate}
        \item For $i$-th row of $[\lambda]$, the starting node has residue $\rho-i+1\mod e'$, and the end node has residue $\rho-i+\lambda_i=\beta_i\mod e'$. So if a bead in $\Ab(\lambda,\rho)$ is on runner~$j$, then the last node in the corresponding row of $\lambda$ has residue $j$, if this row is non-empty.
        \item The number of $(0,1)$-strips in the first column of $[\lambda]$ is equal to $k(\lambda)=\#\{i| i\equiv \rho\mod(e'), \lambda_{i+1}>0\}$.
    \end{enumerate}    
\end{Lemma}
\begin{proof}
    \begin{enumerate}[label=(\alph*)]
        \item The first node in row~$i$ of $\lambda$ has residue $\rho+1-i\mod e'$. There are $\lambda_i$ nodes in the $i$-th row, so the end node is located at $(i,\lambda_i)$ and has residue $\lambda_i-i+\rho$.
        \item A node in the first column has residue $1$ if and only if $\rho+1-i\equiv 1\mod e'$, i.e. $\rho\equiv i\mod e'$. It is part of a $(0,1)$-strip if and only if there exists a node below it, i.e. $\lambda_{i+1}>0$. 
    \end{enumerate}
\end{proof}
\begin{Lemma}\label{easyrowresult}
    Form the box configuration $[\lambda]$ and the abacus configuration $\Ab(\lambda,\rho)$, then
    \begin{enumerate}
        \item The number of rows of $\lambda^+$ is equal to $t+k(\lambda)$ where $t$ is the number of rows of $\lambda$.
        \item The number of rows of $\lambda_+$ is also equal to $t+k(\lambda)$.
    \end{enumerate}
\end{Lemma}
\begin{proof}
    \begin{enumerate}[label=(\alph*)]
        \item This is clear since any $(1,0)$-strip appears in the first column and each of them add one row to $\lambda^+$.
        \item Immediately, since we indeed add $k(\lambda)$ beads contributing to the nonzero parts of the partition.
    \end{enumerate}
\end{proof}
\begin{Lemma}\label{fullrunnerabacusresult}
    Form the abacus configuration $\Ab(\lambda,\rho)$ and write any $\beta_i=a_ie'+b_i$ with $0\leq b_i\leq e'-1$. For each $1\leq j\leq k(\lambda)$, let $s_{j}$ be the integer such that $\lambda_{s_j}$ corresponds to the rightmost bead on level $N_0+k(\lambda)-j$ in $\Ab(\lambda,\rho)$. Then the new partition $\lambda_+$ is of the form $(\lambda_1+a_1-N_0-k(\lambda),\cdots,\lambda_{s_1-1}+a_{s_1-1}-N_0-k(\lambda),\lambda_{s_1}',\lambda_{s_1},\cdots,\lambda_{s_2-1},\lambda_{s_2+1}',\lambda_{s_2},\cdots,\cdots,\lambda_{s_{k(\lambda)}+k(\lambda)-1}',\lambda_{s_{k(\lambda)}},\cdots,\lambda_t)$. To be explicit, this partition has $k(\lambda)+1$ parts:
    \begin{itemize}
        \item the first part consists of $(\lambda_1+a_1-N_0-k(\lambda),\cdots,\lambda_i+a_i-N_0-k(\lambda),\cdots,\lambda_{s_1-1}+a_{s_1-1}-N_0-k(\lambda))$, which correspond to the beads of level $a_i> N_0+k(\lambda)-1$.
        \item for $1\leq j\leq k(\lambda)$, the $j$-th part consists of $(\lambda_{s_j+j-1}',\lambda_{s_1},\lambda_{s_1+1},\cdots,\lambda_{s_{j+1}-1})$, which correspond to the beads on $a_i=N_0+k(\lambda)-j$ level. Here $\lambda_{s_j+j-1}'=(N_0+k(\lambda)-j+1)e'+s_j-1-\rho$ corresponding to the new added bead on the rightmost of $N_0+k(\lambda)-j$ level, and it is the $s_j+j-1$-th bead in $\Ab(\lambda_+,\rho_+)$.
    \end{itemize}
\end{Lemma}
\begin{proof}
    By definition, we add $N_0-N+k(\lambda)$ beads from the top of the added runner on the right-hand side. It is easily seen the top $N_0-N$ beads are not part of the partition, since together with original beads they occupy the whole level. This implies we can assume $N=N_0$.

    Let $\lambda_+=(\beta_1',\cdots,\beta_i',\cdots,\beta_{t+k(\lambda})$ where $t$ is the length of $\lambda$.
    When $a_i>N_0+k(\lambda)-1$, then the order of $\beta_i'$ remain unchanged, that is, the $i$-th bead in $\Ab(\lambda,\rho)$ is the $i$-th bead in $\Ab(\lambda_+,\rho_+)$. So we have $\beta_i'=a_i(e'+1)+b_i=\beta_i+a_i$ and $\beta_i'=\lambda_i'+\rho_+-i$, recall $\rho_+=\rho+N_0+k(\lambda)$, this gives use $\lambda_i'=\lambda_i+a_i-N_0-k(\lambda)$.

    Consider the bead on $a_i=N_0+k(\lambda)-j$ level, the rightmost bead in $\Ab(\lambda,\rho)$, which corresponds to $\lambda_{s_j}$, is not the rightmost in $\Ab(\lambda_+,\rho_+)$ anymore since there is an added bead on the rightmost. By recounting, it is easy to see, the added bead corresponds to $\lambda_{s_j-j+1}'$, while all other beads on this level will increase index by $j=N_0+k(\lambda)-a_i$, so we have $\beta_{i+j}'=a_i(e'+1)+b_i=\beta_i+a_i$ and $\beta_{i+j}'=\lambda_{i+j}'-i-j+\rho_+$. This gives us $\lambda_{i+j}'=\lambda_i$. Hence completes the proof.
\end{proof}
\begin{Remark}\label{0nodecaseremark}
    If $k(\lambda)=0$, then only the first part of \autoref{fullrunnerabacusresult} appear in $\lambda_+$.
\end{Remark}
\begin{Remark}
    It is easy to use \autoref{fullrunnerabacusresult} to write codes to compute $\lambda^+$.
\end{Remark}
\begin{Lemma}\label{count0nodeagain}
    Form the box configuration $[\lambda]$, then either of the following will happen:
    \begin{enumerate}   
        \item If $\rho=0$, then $N_0=-k(\lambda)-1$.
        \item If $1\leq \rho\leq e'-1$, then $N_0=-k(\lambda)$.
    \end{enumerate}
\end{Lemma}
\begin{proof}
    \begin{enumerate}[label=(\alph*)]
        \item $\rho=0$, assume $t$ is the length of $\lambda$. $\beta_{t+1}=0-t-1+\rho=-t-1=a_{t+1}e'+b_{t+1}$.  The starting residue of the $t$-th row is equal to $1-t\mod e'$. The residue pattern of the first column, starts with $0$, then reaches $e'-1,\cdots,2,1$. By definition of $k(\lambda)$, this should repeat $k(\lambda)$ times and the remaining nodes will start with $0$, and hit at $\rho'$ where $1\leq \rho'\leq e'-1$ is the last starting residue. This implies the number of rows $t$ is also equal to $k(\lambda)e'+e'-\rho'+1$.
        
        If $b_{t+1}=e'-1$, then by \autoref{N0def}, $N_0=a_{t+1}+1$. $t=-1-a_{t+1}e'-e'+1=-N_0e'$, so the last starting residue is $1$ and this forces $k(\lambda)e'+e'=-N_0e'\implies N_0=-k(\lambda)-1$
        
        If $1\leq b_{t+1}<e'-1$, $N_0=a_{t+1}$ by \autoref{N0def}. The same argument gives us $t=-1-b_{t+1}-N_0e'$ and the last starting residue is $b_{t+1}+2$. Now $k(\lambda)e'+(e'-b_{t+1}-2+1)=t\implies N_0=-k(\lambda)-1$.
        \item $\rho>0$, assume $t$ is the length of $\lambda$. $\beta_{t+1}=0-t-1+\rho=-t-1+\rho=a_{t+1}e'+b_{t+1}$. The starting residue of the last node is equal to $\rho+1-t\mod e'$. Consider the residue pattern of the first column, it starts with $\rho$, then $\rho-1,\cdots,1,0$. Then it starts with $e'-1,e'-2,\cdots,1,0$, this should repeat $k(\lambda)-1$ times by definition of $k(\lambda)$. At last, it starts with $e'-1$ and decrease until it hits $\rho+1-t\mod e'$. This also implies the number of rows is equal to $\rho+1+(k(\lambda)-1)e'+e'-\rho'$ where $1\leq \rho'\leq e'-1$ is the last starting residue.

        If $b_{t+1}=e'-1$, then by \autoref{N0def}, $N_0=a_{t+1}+1$. This gives us $t=\rho-(a_{t+1}+1)e'=\rho-N_0e'$. The last starting residue is equal to $\rho+1-t\equiv 1\mod e'$. So we count the nodes in the first column and get $t=\rho+1+(k(\lambda)-1)e'+e'-1\implies N_0=-k(\lambda)$.

        If $b_{t+1}<e'-1$, then by \autoref{N0def}, $N_0=a_{t+1}$. This gives us $t=\rho-1+b_{t+1}-N_0e'$ and the last starting residue is $b_{t+1}+2$. Hence the equality $t=\rho+1+(k(\lambda)-1)e'+e'-(b_{t+1+2})\implies N_0=-k(\lambda)$.
    \end{enumerate}
\end{proof}
This result is quite surprising, it means we can get rid of the integer $k(\lambda)$ in the \autoref{abacusdeflambdaplus}  of $\lambda_+$ and restate it as follows:
\begin{Definition}\label{newdeflambdaplus}
    For $\lambda$ a partition of $n$ and $\rho\in I$. Insert a new runner to the rightmost, then add 
    $N':=\begin{cases}
        -1-N & \text{if }\rho=0\\
        -N & \text{else}
    \end{cases}$
    beads from the top of this new runner, that is, for the new abacus, the added beads occupy the positions: $N(e'+1)+e',(N+1)(e'+1)+e',\cdots, (N+N'-1)(e'+1)+e'$.   
    This new abacus configuration corresponds to a new partition and new charge, denoted by $\lambda_+$ and $\rho_+$, respectively.
\end{Definition}
\begin{Corollary}
    $\rho_+=\begin{cases}
                \rho-1=-1 & \text{if }\rho=0\\
                \rho & \text{else }
            \end{cases}$.
\end{Corollary}
\begin{proof}
    This follows from \autoref{count0nodeagain} and \autoref{newrho}.
\end{proof}
By \autoref{easyrowresult}, we observe the number of the rows of $\lambda_+$ and $\lambda^+$ are equal. Therefore, we can establish a bijection between them based on their natural order. Now form the box configuration $[\lambda_+]$ and fill it with the residues, with the starting residue $\rho^+$ at $(1,1)$ coinciding with that of $[\lambda^+]$. It is evident to see there are only two cases:
\begin{equation}
    \rho^+=\begin{cases}
                0 & \text{if }\rho=0\\
                \rho+1 & \text{else, i.e. }1\leq \rho\leq e'-1
            \end{cases}
\end{equation}
\begin{Corollary}\label{endingresidue}
    The ending residue of the rows in corresponding to the new added beads in $[\lambda_+]$ is $0$.
\end{Corollary}
\begin{proof}
    Let $\lambda_+=(\lambda_1',\cdots,\lambda_{t+k(\lambda}')$, if $\rho=0$, the starting residue of the $j$-th row of $[\lambda_+]$ is $1-j\mod e'+1$ and the ending residue of $j$-th row is $1-j+\lambda_j'\mod e'+1$. Recall $\beta_j'=a_j(e'+1)+b_j=\lambda_j'-j+\rho_+$, in particular, take $b_j=e'$ gives the ending residue is $e'-k(\lambda)-N_0\mod e'+1$, now the conclusion follows from \autoref{count0nodeagain}. If $\rho>0$, the starting residue of the $j$-th row is $1-j+\rho^+=1-j+\rho+1=2-j+\rho\mod e'+1$ and the ending residue of this row is $2-j+\rho+\lambda_j-1\mod e'+1$. Play the same game gives us the ending residue is $e'-N_0-k(\lambda)+1\mod e'+1$, completes the proof by \autoref{count0nodeagain}.
\end{proof}

\begin{Lemma}\label{fullrunnerboxresult}
     Assume $k(\lambda)>0$, form the box configuration $[\lambda]$. Then $\lambda_+=\lambda^+$.
\end{Lemma}
\begin{proof}
    By \autoref{fullrunnerabacusresult}, we only need to show in this case $\lambda^+$ has the same expression. There are two cases again, when $\rho=0$ or when $\rho>0$. {\color{Green}Here we only prove the latter case, since just as \autoref{count0nodeagain}, the proofs of the two cases are very similar.}

    Assume the length of $\lambda$ is $t$, for any $1\leq j\leq t$, let $N(j)$ be the number of $0$-nodes in the $j$-th row. Let $\lambda^+=(\lambda_1'',\cdots,\lambda_t'')$. There are $s:=k(\lambda)$ non-trivial $(0,1)$ strips (starting with 1) in the first column, record the row number of the corresponding row with starting residue $1$ by $k_1,\cdots,k_s$. Clearly, $k_1<k_2<\cdots<k_s$. Also let $m_1,\cdots, m_s$ be the number of rows where the $i$-th maximal $(0,1)$-strips end, that is, each maximal $(0,1)$-strip starts from $k_i$-th row and ends at $m_i$-th row. We have $m_1\leq m_2\cdots\leq m_s$ and $m_i>k_i$ for each $1\leq i\leq s$.
    
    Define the following two functions:
    \begin{equation}
        g(j)=
        \begin{cases}
            0 & \text{if }j\leq m_i, \forall i\\
            \max\{i|m_i<j\} & \text{else}
        \end{cases}
    \end{equation}  
    \begin{equation}
        f(j)=
        \begin{cases}
            0 & \text{if }1\leq j\leq k_1\\
            1 & \text{if }k_1+1\leq j\leq k_2\\
            \cdots & \cdots\\
            s-1 & \text{if }k_{s-1}\leq j\leq k_s\\
            s & \text{if }k_s+1\leq j\leq t
        \end{cases}
    \end{equation}
    
    It is clear $g(j)$ is to record when a $(1,0)$-strip ends, while $f(j)$ is to record when it starts. The function $f(j)$ and $g(j)$ are defined on the set $\{1,2,\cdots,t\}$. We may extend it to the whole set of positive integers by letting each (empty) row $k>t$ have starting residue and ending residue $1-k+\rho\mod e'$.

    The starting residue $\rho_j$ of $j$-th row is $1-j+\rho\mod e'$, but we want to find a better expression such that $0\leq \rho_j\leq e'-1$. This can be done by considering the residue patterns of the first column. The starting residue of the first row is $\rho$, then goes down $\rho-1,\cdots,1,0$. 
    It then starts $e'-1,e'-2,\cdots,1,0$, which happen $k(\lambda)-1=s-1$ times. At last, it starts with $e'-1$ and stop at $\rho_{t+1}\equiv 1-t+\rho\mod e'$, which is the $t$-th node in the first column. So it is easy to see $k_i=\rho+(i-1)e'$ and 
    \begin{equation}\label{startingresidueequation}
        \rho_j=
        \begin{cases}
            1+\rho-j+f(j)e' & \text{if }j\neq k_i+1,\forall i\\
            0 & \text{else}
        \end{cases}
    \end{equation}.
    
    We claim 
    \begin{equation}\label{originalrows}
        \lambda_{j+g(j)}''=\lambda_j+N(j)+g(j)-f(j),1\leq j\leq t
    \end{equation}
    \begin{equation}\label{newaddedrows}
        \lambda_{m_i+i}''=
        \begin{cases}
            (1-f(j))e'-\rho+j &\text{if }j\neq k_i+1,\forall i\\
            e'+1 & \text{else}
        \end{cases},1\leq i\leq s=k(\lambda)
    \end{equation}
    This follows from these facts:
    \begin{itemize}
        \item The corresponding row in $[\lambda^+]$ to $j$-th row of $[\lambda]$ has a number of nodes equal to the sum of $\lambda_j$ and the number of $0$-nodes minus the number of $0$-nodes contained in any maximal $(1,0)$-strip. That is, $\lambda_j+N(j)+g(j)-f(j)$.
        \item A maximal $(0,1)$-strip increase the number of rows by $1$, the added row is $1$-row below the strip, that is $m_i+1$. This new row ends with $0$ and the nodes before it are equal to the number of nodes of the orginal $m_i$-th row which are before or equal to this $0$-node.
        \item Assume $j=m_i$, It is easy to see the next $g(j)$ rows in $[\lambda^+]$ must ends with $0$ by \autoref{lambdaplus}.
    \end{itemize}

    It remains to show \autoref{originalrows} and \autoref{newaddedrows} coincide with \autoref{fullrunnerabacusresult}. By \autoref{count0nodeagain}, when $\rho>0$, we have $s=N_0=-k(\lambda)$. By \autoref{endingresidue} and the expressions in \autoref{fullrunnerabacusresult}, we see the index of $\lambda^+=(\lambda_1'',\cdots,\lambda_{t+s}'')$ is the same as $\lambda_+=(\lambda_1',\cdots,\lambda_{t+s}')$. Only need to show $\lambda_i'=\lambda_i''$ for each $i$.

    If $\rho_j\neq 0$, then by \autoref{startingresidueequation}, $\rho_j=1+\rho-j+f(j)e'$, assume $\lambda_j-j+\rho=\beta_j=a_je'+b_j$, then $b_j$ is the ending residue of this row. We have $\lambda_j-(e'-\rho_j)=(a_j+f(j)-1)e'+b_j+1$, hence number of $0$-nodes in this row is equal to $N_j=a_j+f(j)$. If $\rho_j=0$, then $j=k_i+1=\rho+1+f(j)e'-e'$ for some $i$. This time we have $\lambda_j=(a_i+f(j)-1)e'+b_i+1$ and again this implies there are $a_i+f(j)$ $0$-nodes in this row. We have shown $N(j)=a_j+f(j)$. Before the $m_1+1$-th row, $g(j)=0$, this implies $\lambda_+=\lambda^+$ by looking at \autoref{originalrows} and first part of \autoref{fullrunnerabacusresult}. 

    The $m_1+1$-th row is very special. From this row (possibly earlier but no later), all the $0$-nodes are coming from the maximal $(1,0)$-strip. This is because those $(0,1)$-strip are located to the right of the rightmost $(1,0)$-strip. On the $m_{1}+1$-th row, at least the rightmost $(1,0)$-strip does not exist anymore. Therefore, the claim follows.
    
    Let $j\geq m_1+1$, we compare the number of $0$-nodes for adjacent rows. Clearly, the increased $0$-node comes from the $0$-node of a new $(1,0)$-strip and decreased $0$-nodes come from the vanishing of $(1,0)$-strip. So we have $N(j)-N(j+1)=f(j)-f(j+1)+g(j+1)-g(j)$. As $N(j)=a_j+f(j)$, this implies $a_j+g(j)$ is always a constant for any $j\geq m_1+1$. Take $j=t+1$, then $\beta_{t+1}=-t-1+\rho=a_{t+1}e'+b_{t+1}$. If $b_{t+1}=e'-1$, then $a_{t+1}+1=N_0=-k(\lambda)=-s$. In this case, the residue pattern of the first column from $k_s$-row is just $1,0,e'-1,\cdots,2,1$. Hence $g(t+1)=s+1$ and $a_{t+1}+g(t+1)=0$. {\color{Green} (Here is where we need to extend the function g(j).)} If $b_{t+1}<e'-1$, then $a_{t+1}=N_0=-s$ and $g_{t+1}=s$ imply $a_{t+1}+g(t+1)=0$. So in any case we have $a_j+g(j)=0$. Compare with \autoref{fullrunnerabacusresult} shows the proof is complete except for the rows corresponding to the newly added bead, that is, rows with ending residue $0$.

    Consider the $j=m_i$-th row, assume $g(j+1)=x$, then there are $x$ $(1,0)$-strip end in this row. We have $\lambda_{m_r+r}''=\lambda_{m_i+r}''=$ for $i\leq r\leq i+x-1$. When $\rho_j\neq 0$, then we know $\lambda_{m_i+r}''=e'-\rho_j+(f(m_i)-1)e'+1=(1-r)e'+m_i-\rho$ by \autoref{startingresidueequation}. When $\rho_j=0$, then $\lambda_{m_i+r}''=(f(m_i)-r)e'-1$. In this case, $m_i=k_{i'}+1=\rho+1+(i'-1)e'$ for some $i'\geq i$ and hence $f(m_i)=i'$. So $\lambda_{m_i+r}=(i'-r)e'+1=(1-r)e'+(i'-1)e'+1=(1-r)e'+m_i-\rho$. Compare with \autoref{fullrunnerabacusresult}, since $s_i-1=m_i$ by definition, we complete the whole proof.
\end{proof}
\begin{Corollary}\label{lambdaplusequality0case}
    If $k(\lambda)=0$, then $\lambda_+=\lambda^+$.
\end{Corollary}
\begin{proof}
    See \autoref{0nodecaseremark}, this follows from the proof of \autoref{fullrunnerboxresult}. Later we give another explicit proof in \autoref{secondproof}.
\end{proof}
\begin{Theorem}\label{equivdefinition}
    $\lambda_+=\lambda^+$.
\end{Theorem}
\begin{proof}
    It follows from \autoref{fullrunnerboxresult} and \autoref{lambdaplusequality0case}.
\end{proof}

\begin{Corollary}\label{equivdefinitionhighlevel}
    Let $\blam$ be an $\ell$-partition of $n$ and $\brho\in I^\ell$, then $\blam_+=\blam^+$.
\end{Corollary}
\begin{Remark}\label{anotherproof}
    We remind the reader that \autoref{equivdefinition} and \autoref{equivdefinitionhighlevel} give another proof of \autoref{welldefineoflambdaplus}: $\lambda^+$ can be viewed as a composition at first. As $\lambda_+$ is always a partition, the equality shows $\lambda^+$ is also a partition.
\end{Remark}
Hence, from now on, we can freely use these two equivalent definitions of the new partition $\blam^+$. Our goal is to show this new partition corresponds to the new idempotent diagram after subdividing $1_\blam$.
\subsection{Adding From the Left}\label{addingfromtheleft}
There is another abacus definition of $\lambda_+$, which is equivalent to \autoref{abacusdeflambdaplus}. Moreover, we give another proof of 
\autoref{lambdaplusequality0case}.
\begin{Definition}\label{emptyrunner}
    Fix $\lambda$ and $\rho$ and form the abacus configuration $\Ab(\lambda,\rho)$. If $k(\lambda)=0$, Insert a new runner to the leftmost, then add $N_0-N+1$ beads from the top of this new runner, that is, for the new abacus, the added beads occupy the positions: $N(e'+1),(N+1)(e'+1),\cdots,N_0(e'+1)$. This new abacus configuration, it corresponds to a new partition and new charge, denoted by $\lambda_+'$ and $\rho_+'$, respectively.
\end{Definition}
\begin{Lemma}
    $\rho_+'=\rho+N_0+1$.
\end{Lemma}
\begin{proof}
    When insert from the leftmost as in \ref{emptyrunner}, the new abacus has $\rho-Ne'+N_0-N+1$ beads, which should be equal to $\rho_+'-N(e'+1)$. Hence, we have $\rho_+'=\rho+N_0+1$.
\end{proof}
\begin{Example}
    If $\lambda=(4,3,2)$ and $\rho=3$, then $N_0=0$ and we take $N=-1$. Notice that $k(\lambda)=0$. The abacus configuration $\Ab(\lambda,\rho)$ is:
        \begin{equation}
            \abacus(lmmr,bbbb,nnbn,bnbn)
        \end{equation}
        Since $k(\lambda)=0$, we add a runner to the leftmost and put $N_0-N+1=2$ beads from the top on this runner, hence we get:
        \begin{equation}
            \abacus(lmmmr,obbbb,onnbn,nbnbn)
        \end{equation}
        As $\rho_+'=\rho+N_0+1=4$, we translate the abacus to the partition $\lambda_+'=(5,4,2)$. Compare this example with \ref{abacusexample1}.
\end{Example}
\begin{Proposition}\label{movebeadsleft0}
    When $k(\lambda)=0$, $\lambda_+=\lambda_+'$.
\end{Proposition}
\begin{proof}
    Denote the beta numbers correspond to $\lambda_+$ and $\lambda_+'$ by $\beta_i'$ and $\beta_i''$, respectively. Similarly for the partition parts $\lambda_i'$ and $\lambda_i''$. Both strategies don't change the order of beads, as well as the number of beads. Assume $\beta_i=a_ie'+b_i$ for each $i$, then we have:
    \begin{equation}
        \lambda_i'-i+\rho_+=\beta_i'=a_i(e'+1)+b_i
    \end{equation}
    \begin{equation}
        \lambda_i''-i+\rho_+'=\beta_i''=a_i(e'+1)+b_i+1
    \end{equation}
    \begin{equation}
        \rho_+'=\rho_++1
    \end{equation}
    Those three equations give the expected result.
\end{proof}
From now on, we may identify the $\lambda_+$ and $\lambda_+'$ and only use the former notation. 
\begin{Lemma}\label{count0node}
    Assume $k(\lambda)=0$, form the box configuration $[\lambda]$, then:
    \begin{enumerate}
        \item If $k(\lambda)=0$, assume $\lambda=(\lambda_1,\cdots,\lambda_t)$, then $\lambda_+=(\lambda_1+a_1-N_0,\cdots,\lambda_t+a_t-N_0)$ where $\beta_i=a_ie'+b_i$, i.e. $i$-th bead is put on the $a_i$-runner and $N_0$ is the maximal number such that all positions of level above $N_0$ are occupied by beads.
        \item If $\rho=0$, the length $t$ of $\lambda$ is at most $e'$. $N_0=-1$, and the number of $0$-nodes in the $i$-th row is equal to $a_i+1$.
        \item If $1\leq \rho\leq e'-1$, the length $t$ of $\lambda$ is at most $\rho$. $N_0=0$, and the number of $0$-nodes in the $i$-th row is equal to $a_i$.
    \end{enumerate}
\end{Lemma}
\begin{proof}
    In this proof, whenever we write some number in the form $a_ie'+b_i$ or $a_i'e'+b_i'$, we automatically assume $0\leq b_i,b_i' \leq e'-1$.
    \begin{enumerate}[label=(\alph*)]
        \item The order of beads is not changed, since the added beads are on the top left. If $\beta_i=\lambda_i+\rho-i=a_ie'+b_i$, then $\beta_i''=a_i(e'+1)+b_i+1=\beta_i+a_i+1=\lambda_i+\rho-i+a_i+1$. As $\beta_i''=\lambda_i''+\rho_+'-i$, we have $\lambda_i''=\lambda_i+a_i-N_0$. 
        \item When $\rho=0$, the starting residue of $e'$-th row is equal to $1$, and the starting residue of $(e'+1)$-row is equal to $0$, which yield $1$ to $k(\lambda)$. So there can be at most $e'$ rows and
        the residue of $i$-th row is equal to $1-i+e'$. Let $\beta_{t+1}=a_{t+1}e'+b_{t+1}=-t-1\in [-e'-1,-2]$, then either 
            $\begin{cases}
                a_{t+1}=-1 &\\
                b_{t+1}\in [0,e'-2]
            \end{cases}$
        or 
            $\begin{cases}
                a_{t+1}=-2 &\\
                b_{t+1}=e'-1
            \end{cases}$. 
        Both cases give $N_0=-1$ by \autoref{N0def}. 

        Let's count the number of $0$-nodes in $i$-th row. If $i=1$, then the starting residue is $0$, assume $\lambda_1=a_1'e'+b_1'$, then there are $a_1'$ $0$-nodes if $b_1'=0$, $(a_1'+1)$ $0$-nodes if $b_1'>0$. As $\lambda_i=\beta_i+i-\rho=a_ie'+b_i+i-\rho$, we have $\lambda_1=a_1e'+b_1+1$. If $b_1=e'-1$, then
            $\begin{cases}
                a_1'=a_1+1 &\\
                b_1'=0
            \end{cases}$.
        Else $\begin{cases}
                a_1'=a_1 &\\
                b_1'=b_1+1
            \end{cases}$. In both cases, there are $(a_1+1)$ $0$-nodes in total. 
        If $i>1$, then the starting residue is $1-i+e'$, the first $e'-(1-i+e')=i-1$ nodes don't have reside $0$, assume $\lambda_i-i+1=a_1'e'+b_1'$, then there are $a_1'$ $0$-nodes if $b_1'=0$, $(a_1'+1)$ $0$-nodes if $b_1'>0$, just as above. As $a_ie'+b_i+i=\beta_i+i-0=\lambda_i$, we have $a_ie'+b_i+1=a_i'e'+b_i'$. Just as above case, this implies there are in total $(a_1+1)$ $0$-nodes in $i$-th row.  
        \item When $1\leq \rho\leq e'-1$, the number of rows $t$ is at most $\rho$ since the starting residue of $\rho$-th row is $1$ and the starting residue of $(\rho+1)$-th row is $0$, if they exist. $\beta_{t+1}=-t-1+\rho\in [-1,\rho-2]$. Let $\beta_{t+1}=a_{t+1}e'+b_{t+1}$, then 
        either 
            $\begin{cases}
                a_{t+1}=-1 &\\
                b_{t+1}=e'-1
            \end{cases}$
        or 
            $\begin{cases}
                a_{t+1}=0 &\\
                b_{t+1}\in [0,\rho-2]
            \end{cases}$. 
        Both cases give $N_0=0$ by \autoref{N0def}. (The latter case can happen only when $\rho>1$.)

        Let's count the number of $0$-nodes in $i$-th row. The starting residue of $i$-th row is $\rho+1-i$, clearly the first $e'-(\rho+1-i)$ nodes don't have residue $0$, so let $\lambda_i-e'+\rho+1-i=a_i'e'+b_i'$. If $b_i'=0$, there are $a_i'$ $0$-nodes, while if $b_i'>0$, then there are $(a_i'+1)$ $0$-nodes. As $\lambda_i=\beta_i+i-\rho=a_ie'+b_i+i-\rho$, we have $a_i'e'+b_i'=(a_i-1)e'+b_i+1$. Essentially the same argument as above will give us there are in total $a_i$ $0$-nodes in the $i$-th row.
    \end{enumerate}
\end{proof}
\begin{Proposition}\label{secondproof}
    When $k(\lambda)=0$, $\lambda_+=\lambda^+$.
\end{Proposition}
\begin{proof}
    If $k(\lambda)=0$, there is no maximal $(1,0)$-strip, hence the number of rows in $[\lambda^+]$ is equal to that in $[\lambda]$. Moreover, any $0$-node in a row of $[\lambda]$ will add one column to the same row in $[\lambda^+]$, i.e. $\lambda^+_i=\lambda_i+\#\{\text{ number of }0\text{-nodes in }i\text{-th row}\}$. Now the result follows from \autoref{count0node}.
\end{proof}
\begin{Remark}
    Note \autoref{secondproof} gives a second proof of \autoref{lambdaplusequality0case}.
\end{Remark}
The result in this section is really not surprising, considering the following lemma, we can do the same thing when $k(\lambda)>0$:
\begin{Lemma}\label{movebeadsleft}
    Take an abacus configuration, it gives a partition $\lambda$ and a charge $\rho$. Move all beads to the left by $1$ gives a new abacus configurations corresponding to the same $\lambda$.
\end{Lemma}
\begin{proof}
    The proof goes exactly like \autoref{movebeadsleft0}. The order and distance of all beads are unchanged by moving beads to left by $1$, so it automatically give the same partition, but with different $\rho$. 
\end{proof}
\begin{Corollary}\label{movebeadsright}
    Move all the beads to the right by $1$ (and add an essential bead to the northwest) give the same partition. \hfill$\square$
\end{Corollary}
\begin{Theorem}\label{addfromleft}
    Given a partition $\lambda$ of $n$, then $\lambda^+$ is given by putting a runner to the leftmost and add $k(\lambda)+N_0-N+1$ beads from the top of this runner.
\end{Theorem}
\begin{proof}
    This follows immediately from \autoref{equivdefinition}, \autoref{movebeadsleft} and \autoref{movebeadsright}.
\end{proof}
\begin{Example}
    As in \autoref{abacusexample1}, take $\lambda=(4,4,4,4)$ and $\rho=1$, then $N_0=-1$ and we take $N=-2$. Also notice $k(\lambda)=1$.
        The abacus configuration $\Ab(\lambda,\rho)$ is:
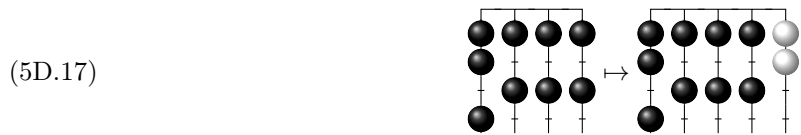
        \begin{equation}
            \abacus(lmmr,bbbb,bnnn,nbbb,bnnn)
            \mapsto
            \abacus(lmmmr,bbbbo,bnnno,nbbbn,bnnnn)
        \end{equation}
        So $\lambda_+=(5,4,4,4,3)$ and $\rho_+=1$. We may move the beads of the second abacus to the right by $1$. This gives us the following abacus configuration:
        \begin{equation}
            \abacus(lmmmr,obbbb,obnnn,onbbb,nbnnn)
        \end{equation}
        It corresponds to the partition $\lambda_+=(5,4,4,4,3)$ and $\rho_+'=2$.
\end{Example}
If the partition $\lambda$ and $\rho$ are fixed, it is not true to move all beads to the left by any length. This is because the last non-zero parts might disappear if the $N$ is not small enough. (recall $N<0$) To deal with this problem, we shall always take $|N|$ to be large enough.
\begin{Remark}\label{correctcomment}
    We have seen now adding from the leftmost is really the same as adding from the rightmost, however, do not forget we're subdividing at the $0\mapsto 1$ edge. For general subdivision at other edges, it is better to adopt 'adding from the left' for consistence. See \autoref{generalcase}.
\end{Remark}

\subsection{Idempotent Correspondence}
Take a quiver $\Gamma$ of type $A^{(1)}_e$, it defines a root datum. 
Fix a positioning set $X$ and $\beta\in Q^+$ such that $\beta=\sum_{i\in I}b_i\alpha_i$ and $\sum_{i\in I}b_i=n$. Also fix $\rho\in I^\ell$ and an asymptotic charge $\charge\in\mathbb{Z}^\ell$,
form the KLRW algebra $\WA(X)$ and $\WA[n](X):=\bigoplus\limits_{\beta\in Q_n^+}\WA(X)$. 

Take $\blam$ a $\ell$-partition for $n$, we can construct the idempotent $1_\lambda\in\WA[n](X)$ and corresponding cell module $\triangle(\blam):=1_\blam\WA(X)$, simple head $L(\blam):=\triangle(\blam)/rad(\triangle(\blam))$. Let $\res(\blam)$ be the residue set of $\blam$ with multiplicities, define $\res(1_\blam):=\sum\limits_{i\in\res(\blam)}\alpha_i$. Assume $\res(1_\blam)=\beta$.

Denote $S_\Gamma:\WA(X)\to W^\rho_{\overline{\beta}}(\overline{X})$ the subdivision map, then $1_\blam^+:=S_\Gamma(1_\blam)$ is an idempotent in $W^\rho_{\overline{\beta}}(\overline{X})$.
\begin{Lemma}\label{correspondencefourcases}
    Maximal $(0,1)$-strips of type $(a)(b)(c)(d)$ in \autoref{lambdaplus} corresponds to the maximal close $(0,1)$-tuple of type $(a)(b)(c)(d)$ in \autoref{closetuple}, respectively.
\end{Lemma}
\begin{proof}
    We only need to show that maximal $(0,1)$-strips correspond to maximal $(0,1)$-tuple in the string diagram. Then the fours cases correspond to each other is self-clear. The correspondence now follows from our strategy of drawing diagrams:
    \begin{itemize}
        \item in the same row, the string corresponding to later node is put further to the right.
        \item the first node in $r$-th row  is stopped by either the solid string or the ghost string of the first node in $r-1$-th row.
        \item solid $1$-strings and ghost $0$-strings block each other.
    \end{itemize}
\end{proof}
\begin{Proposition}\label{idempotentequal}
    $1_\blam^+=1_{\blam^+}$.
\end{Proposition}
\begin{proof}
    By \autoref{correspondencefourcases}, we only need to show in the four cases, the claim holds.
    For instance, let's show $(c)$. It suffices to show the local pattern coincide.
    
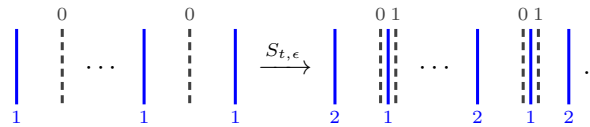
\begin{center}        
        \begin{tikzpicture}[anchorbase,smallnodes,rounded corners]
            \draw[solid](1,1)node[above,yshift=-1pt]{$\phantom{i}$}--++(0,-1)node[below]{$1$};
            \draw[ghost](1.6,1)node[above,yshift=-1pt]{$0$}--++(0,-1)node[below]{$\phantom{i}$};
        \end{tikzpicture}
        $\cdots$
        \begin{tikzpicture}[anchorbase,smallnodes,rounded corners]
            \draw[solid](1,1)node[above,yshift=-1pt]{$\phantom{i}$}--++(0,-1)node[below]{$1$};
            \draw[ghost](1.6,1)node[above,yshift=-1pt]{$0$}--++(0,-1)node[below]{$\phantom{i}$};
            \draw[solid](2.2,1)node[above,yshift=-1pt]{$\phantom{i}$}--++(0,-1)node[below]{$1$};
        \end{tikzpicture}
        $\xrightarrow[]{S_{t,\epsilon}}$
        \begin{tikzpicture}[anchorbase,smallnodes,rounded corners]
            \draw[solid](1.7,1)node[above,yshift=-1pt]{$\phantom{i}$}--++(0,-1)node[below]{$1$};
            \draw[ghost](1.8,1)node[above,yshift=-1pt]{$1$}--++(0,-1)node[below]{$\phantom{i}$};
            \draw[solid](1,1)node[above,yshift=-1pt]{$\phantom{i}$}--++(0,-1)node[below]{$2$};
            \draw[ghost](1.6,1)node[above,yshift=-1pt]{$0$}--++(0,-1)node[below]{$\phantom{i}$};
        \end{tikzpicture}
        $\cdots$
        \begin{tikzpicture}[anchorbase,smallnodes,rounded corners]    
            \draw[solid](1.7,1)node[above,yshift=-1pt]{$\phantom{i}$}--++(0,-1)node[below]{$1$};
            \draw[ghost](1.8,1)node[above,yshift=-1pt]{$1$}--++(0,-1)node[below]{$\phantom{i}$};
            \draw[solid](1,1)node[above,yshift=-1pt]{$\phantom{i}$}--++(0,-1)node[below]{$2$};
            \draw[ghost](1.6,1)node[above,yshift=-1pt]{$0$}--++(0,-1)node[below]{$\phantom{i}$};
            \draw[solid](2.2,1)node[above,yshift=-1pt]{$\phantom{i}$}--++(0,-1)node[below]{$2$};
        \end{tikzpicture}.
    \end{center}
    We remind the reader that solid $2$-strings should be blocked by ghost $1$-strings, ghost $0$-strings are blocked by solid $1$-strings. After the subdivision, the idempotent diagram $1_\blam^+$ can be further transformed since the solid $2$-strings can be pulled further to the right. Then it becomes(zoom a bit):
    
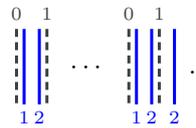
\begin{center}
        \begin{tikzpicture}[anchorbase,smallnodes,rounded corners]
            \draw[solid](1.7,1)node[above,yshift=-1pt]{$\phantom{i}$}--++(0,-1)node[below]{$1$};
            \draw[ghost](2.0,1)node[above,yshift=-1pt]{$1$}--++(0,-1)node[below]{$\phantom{i}$};
            \draw[solid](1.9,1)node[above,yshift=-1pt]{$\phantom{i}$}--++(0,-1)node[below]{$2$};
            \draw[ghost](1.6,1)node[above,yshift=-1pt]{$0$}--++(0,-1)node[below]{$\phantom{i}$};
        \end{tikzpicture}
        $\cdots$
        \begin{tikzpicture}[anchorbase,smallnodes,rounded corners]    
            \draw[solid](1.7,1)node[above,yshift=-1pt]{$\phantom{i}$}--++(0,-1)node[below]{$1$};
            \draw[ghost](2,1)node[above,yshift=-1pt]{$1$}--++(0,-1)node[below]{$\phantom{i}$};
            \draw[solid](1.9,1)node[above,yshift=-1pt]{$\phantom{i}$}--++(0,-1)node[below]{$2$};
            \draw[ghost](1.6,1)node[above,yshift=-1pt]{$0$}--++(0,-1)node[below]{$\phantom{i}$};
            \draw[solid](2.2,1)node[above,yshift=-1pt]{$\phantom{i}$}--++(0,-1)node[below]{$2$};
        \end{tikzpicture}.\\
    \end{center}
    While from $\lambda$ to $\lambda^+$, we use the box configuration to describe the transform:\\
    \hspace*{5cm}\Tableau[scale=0.5]{{1,\\,\\,\\,\\},{0,1,\\,\\,\\},{\\,\cdots,\cdots,\\,\\}, {\\,\\,\cdots,\cdots,\\},{\\,\\,\\,0,1}}
        $\xrightarrow[]{S_{t,\epsilon}}$
        \Tableau[scale=0.5]{{2,\\,\\,\\,\\},{1,2,\\,\\,\\},{0,\cdots,\cdots,\\,\\}, {\\,\cdots,\cdots,\cdots,\\},{\\,\\,\cdots,1,2},{\\,\\,\\,0}}.\\
    So the diagram of $1_{\lambda^+}$ (locally of the $(0,1)$-strip together with adjacent $2$-nodes) should be:
    
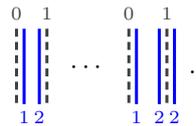
\begin{center}
        \begin{tikzpicture}[anchorbase,smallnodes,rounded corners]
            \draw[solid](1.7,1)node[above,yshift=-1pt]{$\phantom{i}$}--++(0,-1)node[below]{$1$};
            \draw[ghost](2.0,1)node[above,yshift=-1pt]{$1$}--++(0,-1)node[below]{$\phantom{i}$};
            \draw[solid](1.9,1)node[above,yshift=-1pt]{$\phantom{i}$}--++(0,-1)node[below]{$2$};
            \draw[ghost](1.6,1)node[above,yshift=-1pt]{$0$}--++(0,-1)node[below]{$\phantom{i}$};
        \end{tikzpicture}
        $\cdots$
        \begin{tikzpicture}[anchorbase,smallnodes,rounded corners]
            \draw[solid](1.7,1)node[above,yshift=-1pt]{$\phantom{i}$}--++(0,-1)node[below]{$1$};
            \draw[ghost](2.1,1)node[above,yshift=-1pt]{$1$}--++(0,-1)node[below]{$\phantom{i}$};
            \draw[solid](2,1)node[above,yshift=-1pt]{$\phantom{i}$}--++(0,-1)node[below]{$2$};
            \draw[ghost](1.6,1)node[above,yshift=-1pt]{$0$}--++(0,-1)node[below]{$\phantom{i}$};
            \draw[solid](2.2,1)node[above,yshift=-1pt]{$\phantom{i}$}--++(0,-1)node[below]{$2$};
        \end{tikzpicture}.\\
    \end{center}
    Compare the two diagrams we see $1_\lambda^+=1_{\lambda^+}$. The case (d) is essentially the same with (c), so we omit it. (The next \autoref{typedeg}, however, is in the case (d).)

    Next, we consider cases (a) and (b). Again, they're essentially the same, so we only show (a). If there is only a single $0$-node, the claim is obvious, so let's consider the non-trivial case:
    
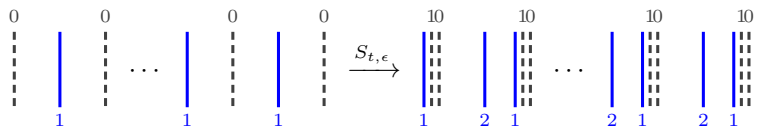
\begin{center}
        \begin{tikzpicture}[anchorbase,smallnodes,rounded corners]
            \draw[solid](1,1)node[above,yshift=-1pt]{$\phantom{i}$}--++(0,-1)node[below]{$1$};
            \draw[ghost](0.4,1)node[above,yshift=-1pt]{$0$}--++(0,-1)node[below]{$\phantom{i}$};
            \draw[ghost](1.6,1)node[above,yshift=-1pt]{$0$}--++(0,-1)node[below]{$\phantom{i}$};
        \end{tikzpicture}
        $\cdots$
        \begin{tikzpicture}[anchorbase,smallnodes,rounded corners]
            \draw[solid](-0.2,1)node[above,yshift=-1pt]{$\phantom{i}$}--++(0,-1)node[below]{$1$};
            \draw[ghost](0.4,1)node[above,yshift=-1pt]{$0$}--++(0,-1)node[below]{$\phantom{i}$};
            \draw[solid](1,1)node[above,yshift=-1pt]{$\phantom{i}$}--++(0,-1)node[below]{$1$};
            \draw[ghost](1.6,1)node[above,yshift=-1pt]{$0$}--++(0,-1)node[below]{$\phantom{i}$};
        \end{tikzpicture}
        $\xrightarrow[]{S_{t,\epsilon}}$
        \begin{tikzpicture}[anchorbase,smallnodes,rounded corners]
            \draw[solid](0.2,1)node[above,yshift=-1pt]{$\phantom{i}$}--++(0,-1)node[below]{$1$};
            \draw[ghost](0.3,1)node[above,yshift=-1pt]{$1$}--++(0,-1)node[below]{$\phantom{i}$};
            \draw[solid](1.4,1)node[above,yshift=-1pt]{$\phantom{i}$}--++(0,-1)node[below]{$1$};
            \draw[ghost](1.5,1)node[above,yshift=-1pt]{$1$}--++(0,-1)node[below]{$\phantom{i}$};
            \draw[solid](1,1)node[above,yshift=-1pt]{$\phantom{i}$}--++(0,-1)node[below]{$2$};
            \draw[ghost](0.4,1)node[above,yshift=-1pt]{$0$}--++(0,-1)node[below]{$\phantom{i}$};
            \draw[ghost](1.6,1)node[above,yshift=-1pt]{$0$}--++(0,-1)node[below]{$\phantom{i}$};
        \end{tikzpicture}
        $\cdots$
        \begin{tikzpicture}[anchorbase,smallnodes,rounded corners]              
            \draw[solid](0.2,1)node[above,yshift=-1pt]{$\phantom{i}$}--++(0,-1)node[below]{$1$};
            \draw[ghost](0.3,1)node[above,yshift=-1pt]{$1$}--++(0,-1)node[below]{$\phantom{i}$};
            \draw[solid](-0.2,1)node[above,yshift=-1pt]{$\phantom{i}$}--++(0,-1)node[below]{$2$};
            \draw[ghost](0.4,1)node[above,yshift=-1pt]{$0$}--++(0,-1)node[below]{$\phantom{i}$};
            \draw[solid](1.4,1)node[above,yshift=-1pt]{$\phantom{i}$}--++(0,-1)node[below]{$1$};
            \draw[ghost](1.5,1)node[above,yshift=-1pt]{$1$}--++(0,-1)node[below]{$\phantom{i}$};
            \draw[solid](1,1)node[above,yshift=-1pt]{$\phantom{i}$}--++(0,-1)node[below]{$2$};
            \draw[ghost](1.6,1)node[above,yshift=-1pt]{$0$}--++(0,-1)node[below]{$\phantom{i}$};
        \end{tikzpicture}.
    \end{center}
    We pull the solid $2$-strings and solid $1$-strings further to the right and get:
    
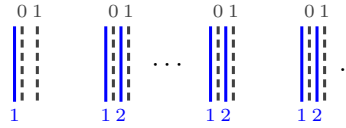
\begin{center}
        \begin{tikzpicture}[anchorbase,smallnodes,rounded corners]
            \draw[solid](0.3,1)node[above,yshift=-1pt]{$\phantom{i}$}--++(0,-1)node[below]{$1$};
            \draw[ghost](0.6,1)node[above,yshift=-1pt]{$1$}--++(0,-1)node[below]{$\phantom{i}$};
            \draw[ghost](0.4,1)node[above,yshift=-1pt]{$0$}--++(0,-1)node[below]{$\phantom{i}$};
            \draw[solid](1.5,1)node[above,yshift=-1pt]{$\phantom{i}$}--++(0,-1)node[below]{$1$};
            \draw[ghost](1.8,1)node[above,yshift=-1pt]{$1$}--++(0,-1)node[below]{$\phantom{i}$};
            \draw[solid](1.7,1)node[above,yshift=-1pt]{$\phantom{i}$}--++(0,-1)node[below]{$2$};
            \draw[ghost](1.6,1)node[above,yshift=-1pt]{$0$}--++(0,-1)node[below]{$\phantom{i}$};
        \end{tikzpicture}
        $\cdots$
        \begin{tikzpicture}[anchorbase,smallnodes,rounded corners]              
            \draw[solid](0.3,1)node[above,yshift=-1pt]{$\phantom{i}$}--++(0,-1)node[below]{$1$};
            \draw[ghost](0.6,1)node[above,yshift=-1pt]{$1$}--++(0,-1)node[below]{$\phantom{i}$};
            \draw[solid](0.5,1)node[above,yshift=-1pt]{$\phantom{i}$}--++(0,-1)node[below]{$2$};
            \draw[ghost](0.4,1)node[above,yshift=-1pt]{$0$}--++(0,-1)node[below]{$\phantom{i}$};
            \draw[solid](1.5,1)node[above,yshift=-1pt]{$\phantom{i}$}--++(0,-1)node[below]{$1$};
            \draw[ghost](1.8,1)node[above,yshift=-1pt]{$1$}--++(0,-1)node[below]{$\phantom{i}$};
            \draw[solid](1.7,1)node[above,yshift=-1pt]{$\phantom{i}$}--++(0,-1)node[below]{$2$};
            \draw[ghost](1.6,1)node[above,yshift=-1pt]{$0$}--++(0,-1)node[below]{$\phantom{i}$};
        \end{tikzpicture}.
    \end{center}

    While from $\lambda$ to $\lambda^+$, we use the box configuration to describe the transform:\\
    \hspace*{5cm}\Tableau[scale=0.5]{{0,1,\\,\\,\\},{\\,0,1,\\,\\},{\\,\\,\cdots,\\,\\}, {\\,\\,\\,0,1},{\\,\\,\\,\\,0}}
    $\xrightarrow[]{S_{t,\epsilon}}$
    \Tableau[scale=0.5]{{0,1,2,\\,\\,\\},{\\,0,1,2,\\,\\},{\\,\\,\cdots,\\,\\,\\}, {\\,\\,\\,0,1,2},{\\,\\,\\,\\,0,1}}.\\
    So the diagram of $1_{\lambda^+}$ (locally of the $(0,1)$-strip together with adjacent $2$-nodes) should be:
    
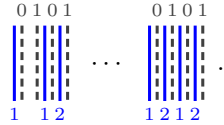
\begin{center}
        \begin{tikzpicture}[anchorbase,smallnodes,rounded corners]
            \draw[solid](1.1,1)node[above,yshift=-1pt]{$\phantom{i}$}--++(0,-1)node[below]{$1$};
            \draw[ghost](1.4,1)node[above,yshift=-1pt]{$1$}--++(0,-1)node[below]{$\phantom{i}$};
            \draw[ghost](1.2,1)node[above,yshift=-1pt]{$0$}--++(0,-1)node[below]{$\phantom{i}$};
            \draw[solid](1.5,1)node[above,yshift=-1pt]{$\phantom{i}$}--++(0,-1)node[below]{$1$};
            \draw[ghost](1.8,1)node[above,yshift=-1pt]{$1$}--++(0,-1)node[below]{$\phantom{i}$};
            \draw[solid](1.7,1)node[above,yshift=-1pt]{$\phantom{i}$}--++(0,-1)node[below]{$2$};
            \draw[ghost](1.6,1)node[above,yshift=-1pt]{$0$}--++(0,-1)node[below]{$\phantom{i}$};
        \end{tikzpicture}
        $\cdots$
        \begin{tikzpicture}[anchorbase,smallnodes,rounded corners]              
            \draw[solid](1.1,1)node[above,yshift=-1pt]{$\phantom{i}$}--++(0,-1)node[below]{$1$};
            \draw[ghost](1.4,1)node[above,yshift=-1pt]{$1$}--++(0,-1)node[below]{$\phantom{i}$};
            \draw[solid](1.3,1)node[above,yshift=-1pt]{$\phantom{i}$}--++(0,-1)node[below]{$2$};
            \draw[ghost](1.2,1)node[above,yshift=-1pt]{$0$}--++(0,-1)node[below]{$\phantom{i}$};
            \draw[solid](1.5,1)node[above,yshift=-1pt]{$\phantom{i}$}--++(0,-1)node[below]{$1$};
            \draw[ghost](1.8,1)node[above,yshift=-1pt]{$1$}--++(0,-1)node[below]{$\phantom{i}$};
            \draw[solid](1.7,1)node[above,yshift=-1pt]{$\phantom{i}$}--++(0,-1)node[below]{$2$};
            \draw[ghost](1.6,1)node[above,yshift=-1pt]{$0$}--++(0,-1)node[below]{$\phantom{i}$};
        \end{tikzpicture}.
    \end{center}
    Compare them give the equality $1_\lambda^+=1_\lambda^+$, which completes the proof.
\end{proof}

\begin{Example}\label{typedeg}
    Here we give an explicit example to show what happens in each step. Let $\lambda=(2,2,2)$ and $\rho=1$, then the box configuration $[\lambda]$ is:
    \begin{center}
        \Tableau[scale=0.5]{{1,2},{0,1},{e,0}}
    \end{center}
    the corresponding idempotent diagram is(omit the ghost $2$-string and solid $e$-string):
    
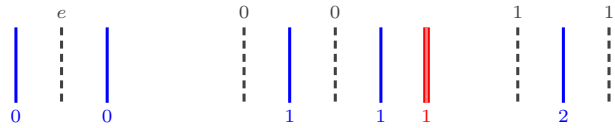
\begin{center}
        \begin{tikzpicture}[anchorbase,smallnodes,rounded corners]                
            \draw[redstring](0,1)--++(0,-1)node[below]{$1$};
            \draw[solid](-5.4,1)node[above,yshift=-1pt]{$\phantom{i}$}--++(0,-1)node[below]{$0$};        
            \draw[ghost](-4.8,1)node[above,yshift=-1pt]{$e$}--++(0,-1)node[below]{$\phantom{i}$}; 
            \draw[solid](-4.2,1)node[above,yshift=-1pt]{$\phantom{i}$}--++(0,-1)node[below]{$0$};
            \draw[ghost](-2.4,1)node[above,yshift=-1pt]{$0$}--++(0,-1)node[below]{$\phantom{i}$};
            \draw[solid](-1.8,1)node[above,yshift=-1pt]{$\phantom{i}$}--++(0,-1)node[below]{$1$};
            \draw[ghost](-1.2,1)node[above,yshift=-1pt]{$0$}--++(0,-1)node[below]{$\phantom{i}$};
            \draw[solid](-0.6,1)node[above,yshift=-1pt]{$\phantom{i}$}--++(0,-1)node[below]{$1$};        
            \draw[ghost](2.4,1)node[above,yshift=-1pt]{$1$}--++(0,-1)node[below]{$\phantom{i}$};
            \draw[solid](1.8,1)node[above,yshift=-1pt]{$\phantom{i}$}--++(0,-1)node[below]{$2$};
            \draw[ghost](1.2,1)node[above,yshift=-1pt]{$1$}--++(0,-1)node[below]{$\phantom{i}$};
        \end{tikzpicture}\\
    \end{center}
    Apply the subdivision, we have the diagram:
    
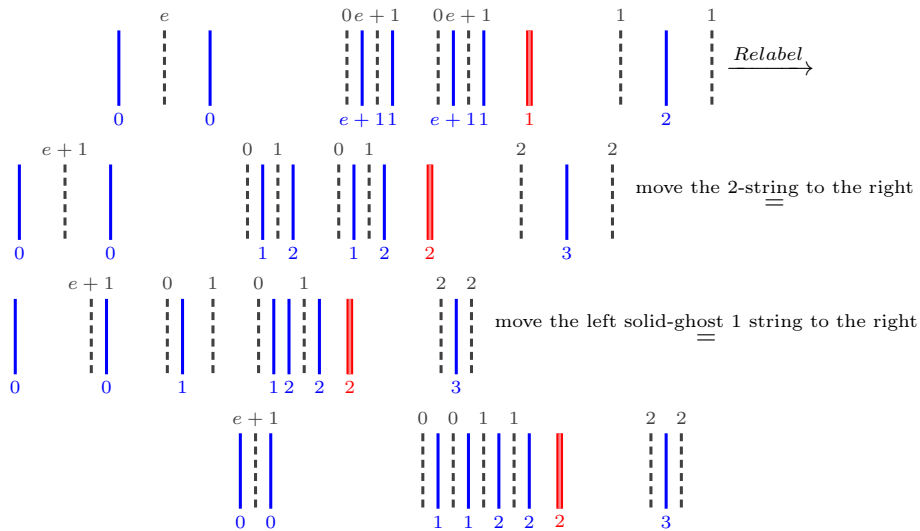
\begin{center}
        \begin{tikzpicture}[anchorbase,smallnodes,rounded corners]                   
            \draw[redstring](0,1)--++(0,-1)node[below]{$1$};
            \draw[solid](-5.4,1)node[above,yshift=-1pt]{$\phantom{i}$}--++(0,-1)node[below]{$0$};        
            \draw[ghost](-4.8,1)node[above,yshift=-1pt]{$e$}--++(0,-1)node[below]{$\phantom{i}$}; 
            \draw[solid](-4.2,1)node[above,yshift=-1pt]{$\phantom{i}$}--++(0,-1)node[below]{$0$};
            \draw[ghost](-2.4,1)node[above,yshift=-1pt]{$0$}--++(0,-1)node[below]{$\phantom{i}$};
            \draw[solid](-2.2,1)node[above,yshift=-1pt]{$\phantom{i}$}--++(0,-1)node[below]{$e+1$};        
            \draw[ghost](-2,1)node[above,yshift=-1pt]{$e+1$}--++(0,-1)node[below]{$\phantom{i}$};
            \draw[solid](-1.8,1)node[above,yshift=-1pt]{$\phantom{i}$}--++(0,-1)node[below]{$1$};
            \draw[ghost](-1.2,1)node[above,yshift=-1pt]{$0$}--++(0,-1)node[below]{$\phantom{i}$};
            \draw[solid](-1,1)node[above,yshift=-1pt]{$\phantom{i}$}--++(0,-1)node[below]{$e+1$};        
            \draw[ghost](-0.8,1)node[above,yshift=-1pt]{$e+1$}--++(0,-1)node[below]{$\phantom{i}$};
            \draw[solid](-0.6,1)node[above,yshift=-1pt]{$\phantom{i}$}--++(0,-1)node[below]{$1$};        
            \draw[ghost](2.4,1)node[above,yshift=-1pt]{$1$}--++(0,-1)node[below]{$\phantom{i}$};
            \draw[solid](1.8,1)node[above,yshift=-1pt]{$\phantom{i}$}--++(0,-1)node[below]{$2$};
            \draw[ghost](1.2,1)node[above,yshift=-1pt]{$1$}--++(0,-1)node[below]{$\phantom{i}$};
        \end{tikzpicture}$\xrightarrow[]{Relabel}$\\
        \begin{tikzpicture}[anchorbase,smallnodes,rounded corners]                   
            \draw[redstring](0,1)--++(0,-1)node[below]{$2$};
            \draw[solid](-5.4,1)node[above,yshift=-1pt]{$\phantom{i}$}--++(0,-1)node[below]{$0$};        
            \draw[ghost](-4.8,1)node[above,yshift=-1pt]{$e+1$}--++(0,-1)node[below]{$\phantom{i}$}; 
            \draw[solid](-4.2,1)node[above,yshift=-1pt]{$\phantom{i}$}--++(0,-1)node[below]{$0$};
            \draw[ghost](-2.4,1)node[above,yshift=-1pt]{$0$}--++(0,-1)node[below]{$\phantom{i}$};
            \draw[solid](-2.2,1)node[above,yshift=-1pt]{$\phantom{i}$}--++(0,-1)node[below]{$1$};        
            \draw[ghost](-2,1)node[above,yshift=-1pt]{$1$}--++(0,-1)node[below]{$\phantom{i}$};
            \draw[solid](-1.8,1)node[above,yshift=-1pt]{$\phantom{i}$}--++(0,-1)node[below]{$2$};
            \draw[ghost](-1.2,1)node[above,yshift=-1pt]{$0$}--++(0,-1)node[below]{$\phantom{i}$};
            \draw[solid](-1,1)node[above,yshift=-1pt]{$\phantom{i}$}--++(0,-1)node[below]{$1$};        
            \draw[ghost](-0.8,1)node[above,yshift=-1pt]{$1$}--++(0,-1)node[below]{$\phantom{i}$};
            \draw[solid](-0.6,1)node[above,yshift=-1pt]{$\phantom{i}$}--++(0,-1)node[below]{$2$};        
            \draw[ghost](2.4,1)node[above,yshift=-1pt]{$2$}--++(0,-1)node[below]{$\phantom{i}$};
            \draw[solid](1.8,1)node[above,yshift=-1pt]{$\phantom{i}$}--++(0,-1)node[below]{$3$};
            \draw[ghost](1.2,1)node[above,yshift=-1pt]{$2$}--++(0,-1)node[below]{$\phantom{i}$};
        \end{tikzpicture} $\stackrel{\text{move the 2-string to the right}}{=}$\\
        
        \begin{tikzpicture}[anchorbase,smallnodes,rounded corners]                         
            \draw[redstring](0,1)--++(0,-1)node[below]{$2$};
            \draw[solid](-4.4,1)node[above,yshift=-1pt]{$\phantom{i}$}--++(0,-1)node[below]{$0$};        
            \draw[ghost](-3.4,1)node[above,yshift=-1pt]{$e+1$}--++(0,-1)node[below]{$\phantom{i}$}; 
            \draw[solid](-3.2,1)node[above,yshift=-1pt]{$\phantom{i}$}--++(0,-1)node[below]{$0$};
            \draw[ghost](-2.4,1)node[above,yshift=-1pt]{$0$}--++(0,-1)node[below]{$\phantom{i}$};
            \draw[solid](-2.2,1)node[above,yshift=-1pt]{$\phantom{i}$}--++(0,-1)node[below]{$1$};        
            \draw[ghost](-1.8,1)node[above,yshift=-1pt]{$1$}--++(0,-1)node[below]{$\phantom{i}$};
            \draw[ghost](-1.2,1)node[above,yshift=-1pt]{$0$}--++(0,-1)node[below]{$\phantom{i}$};
            \draw[solid](-1,1)node[above,yshift=-1pt]{$\phantom{i}$}--++(0,-1)node[below]{$1$}; 
            \draw[solid](-0.8,1)node[above,yshift=-1pt]{$\phantom{i}$}--++(0,-1)node[below]{$2$};
            \draw[ghost](-0.6,1)node[above,yshift=-1pt]{$1$}--++(0,-1)node[below]{$\phantom{i}$};
            \draw[solid](-0.4,1)node[above,yshift=-1pt]{$\phantom{i}$}--++(0,-1)node[below]{$2$};        
            \draw[ghost](1.2,1)node[above,yshift=-1pt]{$2$}--++(0,-1)node[below]{$\phantom{i}$};
            \draw[solid](1.4,1)node[above,yshift=-1pt]{$\phantom{i}$}--++(0,-1)node[below]{$3$};
            \draw[ghost](1.6,1)node[above,yshift=-1pt]{$2$}--++(0,-1)node[below]{$\phantom{i}$};
        \end{tikzpicture} $\stackrel{\text{move the left solid-ghost 1 string to the right}}{=}$\\
        \begin{tikzpicture}[anchorbase,smallnodes,rounded corners]                         
            \draw[redstring](0,1)--++(0,-1)node[below]{$2$};
            \draw[solid](-4.2,1)node[above,yshift=-1pt]{$\phantom{i}$}--++(0,-1)node[below]{$0$};        
            \draw[ghost](-4,1)node[above,yshift=-1pt]{$e+1$}--++(0,-1)node[below]{$\phantom{i}$}; 
            \draw[solid](-3.8,1)node[above,yshift=-1pt]{$\phantom{i}$}--++(0,-1)node[below]{$0$};
            \draw[ghost](-1.8,1)node[above,yshift=-1pt]{$0$}--++(0,-1)node[below]{$\phantom{i}$};
            \draw[solid](-1.6,1)node[above,yshift=-1pt]{$\phantom{i}$}--++(0,-1)node[below]{$1$};        
            \draw[ghost](-1,1)node[above,yshift=-1pt]{$1$}--++(0,-1)node[below]{$\phantom{i}$};
            \draw[ghost](-1.4,1)node[above,yshift=-1pt]{$0$}--++(0,-1)node[below]{$\phantom{i}$};
            \draw[solid](-1.2,1)node[above,yshift=-1pt]{$\phantom{i}$}--++(0,-1)node[below]{$1$}; 
            \draw[solid](-0.8,1)node[above,yshift=-1pt]{$\phantom{i}$}--++(0,-1)node[below]{$2$};
            \draw[ghost](-0.6,1)node[above,yshift=-1pt]{$1$}--++(0,-1)node[below]{$\phantom{i}$};
            \draw[solid](-0.4,1)node[above,yshift=-1pt]{$\phantom{i}$}--++(0,-1)node[below]{$2$};        
            \draw[ghost](1.2,1)node[above,yshift=-1pt]{$2$}--++(0,-1)node[below]{$\phantom{i}$};
            \draw[solid](1.4,1)node[above,yshift=-1pt]{$\phantom{i}$}--++(0,-1)node[below]{$3$};
            \draw[ghost](1.6,1)node[above,yshift=-1pt]{$2$}--++(0,-1)node[below]{$\phantom{i}$};
        \end{tikzpicture}\\
    \end{center}
    On the other hand, the idempotent diagram corresponding to $[\lambda^+]=\Tableau[scale=0.7]{{2,3},{1,2},{0,1},{e+1,0}}$ is(omit the ghost $3$-string and solid $e$-string):\\
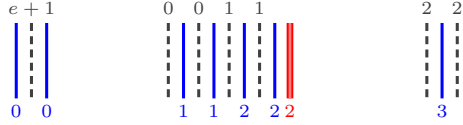
    \begin{center}
        \begin{tikzpicture}[anchorbase,smallnodes,rounded corners]                         
            \draw[redstring](0,1)--++(0,-1)node[below]{$2$};
            \draw[solid](-3.2,1)node[above,yshift=-1pt]{$\phantom{i}$}--++(0,-1)node[below]{$0$};        
            \draw[ghost](-3.4,1)node[above,yshift=-1pt]{$e+1$}--++(0,-1)node[below]{$\phantom{i}$}; 
            \draw[solid](-3.6,1)node[above,yshift=-1pt]{$\phantom{i}$}--++(0,-1)node[below]{$0$};
            \draw[ghost](-1.6,1)node[above,yshift=-1pt]{$0$}--++(0,-1)node[below]{$\phantom{i}$};
            \draw[solid](-1.4,1)node[above,yshift=-1pt]{$\phantom{i}$}--++(0,-1)node[below]{$1$};        
            \draw[ghost](-0.8,1)node[above,yshift=-1pt]{$1$}--++(0,-1)node[below]{$\phantom{i}$};
            \draw[ghost](-1.2,1)node[above,yshift=-1pt]{$0$}--++(0,-1)node[below]{$\phantom{i}$};
            \draw[solid](-1,1)node[above,yshift=-1pt]{$\phantom{i}$}--++(0,-1)node[below]{$1$}; 
            \draw[solid](-0.6,1)node[above,yshift=-1pt]{$\phantom{i}$}--++(0,-1)node[below]{$2$};
            \draw[ghost](-0.4,1)node[above,yshift=-1pt]{$1$}--++(0,-1)node[below]{$\phantom{i}$};
            \draw[solid](-0.2,1)node[above,yshift=-1pt]{$\phantom{i}$}--++(0,-1)node[below]{$2$};        
            \draw[ghost](1.8,1)node[above,yshift=-1pt]{$2$}--++(0,-1)node[below]{$\phantom{i}$};
            \draw[solid](2,1)node[above,yshift=-1pt]{$\phantom{i}$}--++(0,-1)node[below]{$3$};
            \draw[ghost](2.2,1)node[above,yshift=-1pt]{$2$}--++(0,-1)node[below]{$\phantom{i}$};
        \end{tikzpicture}\\
    \end{center}
    
Compare the two idempotent diagrams, we see they're equal. 
\end{Example}
As a consequence, in the KLRW algebra $W^\rho_{\overline{\beta}}(\overline{X})$, $1_\blam^+=1_{\blam^+}$ is an idempotent appearing in the cellular basis. So it defines a cell module $\triangle^+(\blam^+):=1_\blam^+W^\rho_{\overline{\beta}}(\overline{X})$ and corresponding simple head $L^+(\blam^+)$, which is the quotient of $\triangle^+(\blam^+)$ by its radical. We have the graded decomposition number $[\triangle^+(\blam^+):L^+(\bmu^+)]_q$. 

\begin{Proposition}\label{quotientideal}
    Let $J:=\Oneg\WAs(\overline{X})\1_{\text{bad}}\WAs(\overline{X})\Oneg$,  $\blam$ a $\ell$-partition of $n$ and $\brho$ a charge. Fill $[\blam]$ and $[\bmu]$ with residues with regard to $\brho$, if there is at most one $0$-node, then $[\triangle^+(\blam^+):L^+(\bmu^+)]_q=[\triangle'(\blam):L'(\bmu)]_q$ where we use $'$ to denote cell modules and simple heads of $\SA(X)$.
\end{Proposition}
\begin{proof}
    This is clear, as when there is at most one $0$-node in $[\blam]$, there can be at most one $1$-node in $[\blam^+]$, this means there can't be any bad diagram, hence $\triangle^+(\blam^+)\cap J=0=L^+(\bmu^+)\cap J$ and the conclusion follows.
\end{proof}
    
\begin{Theorem}\label{decompositionequal}
    With the same assumption as \autoref{quotientideal}, $[\triangle(\blam):L(\bmu)]_q=[\triangle^+(\blam^+):L^+(\bmu^+)]_q$.
\end{Theorem}
\begin{proof}
    Because $S_{t,\epsilon}$ induces an isomorphism between $\WA(X)$ and $\SA(X)$, the cell modules and simple heads are one-to-one corresponding. As a consequence, we have $[\triangle(\blam):L(\bmu)]_q=[\triangle'(\blam):L'(\bmu)]_q$. By \autoref{quotientideal}, 
    the conclusion then follows.
\end{proof}
\begin{Remark}\label{conjecture}
    {\color{Green}While it is too optimistic, our conjecture is that the assumption of \autoref{decompositionequal} can be removed. However, the challenge is to show when the intersection of the cell module and the quotient ideal $J$ is nonzero, the \autoref{quotientideal} still holds, which is a much more cumbersome task.}
\end{Remark}
\section{General Subdivision}\label{generalcase}
Recall in \autoref{Subdivision} and \autoref{equalityofdecompositionnumber}, we assume we are subdividing the edge $0\mapsto 1$ of the quiver of type $A^{(1)}_{e}$. In this section, we briefly discuss subdividing at another edge, say $i\mapsto i+1$.

As in \autoref{equalityofdecompositionnumber}, there are two ways to define the partition $\lambda^+$, one is to use the box configuration and the other is to use the abacus configuration. 

Let $\lambda$ be a partition of $n$ and $\rho\in I$ be a charge, form the box configuration filled with residues with regard to $\rho$ and denote it by $[\lambda]_\rho$. To make use of the results of previous sections, we relabel the quiver by replacing $i,i+1,\cdots,e,0,\cdots,i-1$ by $0,1,\cdots, e-i,e-i+1,\cdots,i-1$, respectively, with the orientation unchanged. 

Since we have relabeled the quiver, the charge $\rho$ is also changed to: 
\begin{equation}
    \rho'=
\begin{cases}
    \rho-i & \text{if }\rho\geq i\\
    \rho+e+1-i & \text{else}
\end{cases}
\end{equation}
Corresponding to the new $\rho'$, we form the box configuration filled with residues: $[\lambda]_{\rho'}$. It is easy to see, the effect of relabeling the quiver changes  $[\lambda]_\rho$ to $[\lambda]_{\rho'}$. See \autoref{generalexample}. 

Now all $i$-nodes/strings become $0$-nodes/strings and all $(i,i+1)$-strips become $(0,1)$-strips, we can apply the previous construction of subdivision (at $0\mapsto 1$ edge) and get $\blam^+, 1_{\blam^+},1^+_{\blam}$. Back to $[\lambda]_\rho$, we can define $\lambda^+$ as the following:
\begin{Definition}
    Given $\lambda$ and $\rho$. Replace any $0$ and $1$ by $i$ and $i+1$, respectively in \autoref{lambdaplus}, and make essential changes to other residues. The new box configuration gives $\lambda^+$.
\end{Definition}

Let's consider the abacus definition of $\lambda_+$. To unify the strategies in general cases, we choose 'adding from the left' as in \autoref{addingfromtheleft}. After relabeling, the new abacus configuration is $\Ab(\lambda,\rho')$. Add a new runner to the leftmost and put $k(\lambda)+N_0-N+1$ beads on this runner from the top. This abacus configuration gives the partition $\lambda_+$.

Consider the difference of $\Ab(\lambda,\rho)$ and $\Ab(\lambda,\rho')$. If $\rho\geq i$, $\rho'=\rho-i$. Then $\Ab(\lambda,\rho)$ can be obtained by moving all beads of $\Ab(\lambda,\rho')$ to the right by $i$ positions. If $\rho<i$, $\rho'=\rho+e+1-i$. Then $\Ab(\lambda,\rho)$ can be obtained by moving all beads of $\Ab(\lambda,\rho')$ to the left by $e+1-i$ positions, given the condition $|N|$ is large enough. See the comment above \autoref{correctcomment}.  In any case, adding a runner to the leftmost of $\Ab(\lambda,\rho')$ is equivalent to adding a runner between $(i-1)$-runner and $i$-runner of $\Ab(\lambda,\rho)$.

However, the number of beads we put to the new runner is different. This is because, $N_0$ might be different for $(\lambda,\rho)$ and $(\lambda,\rho')$. Nevertheless, by the same argument in previous sections, we have the following definition:
\begin{Definition}
    With assumptions above, add a new runner between $(i-1)$-runner and $i$-runner of $\Ab(\lambda,\rho)$, and add $k(\lambda)+N_0-N$ beads to this runner from the top. The new abacus gives a partition, denoted by $\lambda_+$. 
\end{Definition}
\begin{Example}\label{generalexample}
    Fix $e=3,\rho=2$, $\lambda=(4,3,2)$, so the box configuration $[\lambda]$ filled with residues is as following:
            \begin{center}
                \Tableau[scale=0.7]{{2,3,0,1},{1,2,3},{0,1}}
            \end{center}
    \begin{enumerate}
        \item  Assume we are subdividing the edge $2\mapsto 3$, the box configuration after relabeling is as the following:
            \begin{center}
                \Tableau[scale=0.7]{{0,1,2,3},{3,0,1},{2,3}}
            \end{center}
        As $\rho'=\rho-2=0$, the abacus configuration $\Ab(\lambda,\rho')$ is as following:
            \begin{center}
                0\hspace{1mm}  1\hspace{1mm}  2\hspace{1mm}  3\hspace{1mm}\\  
              \abacus(lmmr,bbbb,bnnb,nbnb)
            \end{center}
        Notice that $N_0=-1$ and we take $N=-2$. Add a runner to the leftmost and put $k(\lambda)+N_0-N+1=2$ beads from the top, we have the following new abacus $\Ab(\lambda_+,\rho_+')$:
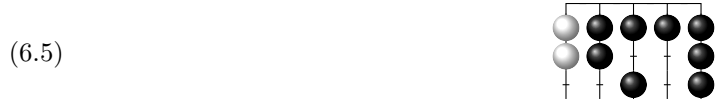
            \begin{equation}\label{ab:1}
                \abacus(lmmmr,obbbb,obnnb,nnbnb)
            \end{equation}
        It is easy to get $\rho_+=(5,4,2)$. If we directly write the abacus $\Ab(\lambda,\rho)$:
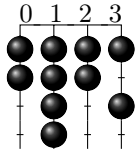
            \begin{center}
                0\hspace{1mm}  1\hspace{1mm}  2\hspace{1mm}  3\hspace{1mm}\\
                \abacus(lmmr,bbbb,bbbn,nbnb,nbnn)
            \end{center}
        Note this time we still take $N_0=-1,N=-2$. Put a new runner between $1$-runner and $2$-runner with $k(\lambda)+N_0-N+1=2$ bead from the top, we have the following $\Ab(\lambda_+,\rho_+)$:
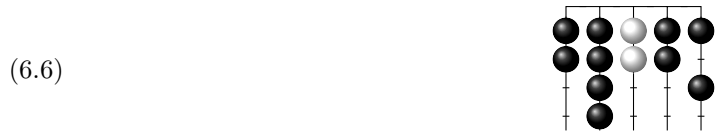
            \begin{equation}\label{ab:2}
                \abacus(lmmmr,bbobb,bbobn,nbnnb,nbnnn)
            \end{equation}
        It is clear \autoref{ab:1} is got from \autoref{ab:2} by moving beads to the right by $2$ (and add two essential beads from the northwest), hence denote the same partition $\lambda_+$.
        \item  Assume we are subdividing the edge $3\mapsto 0$, the box configuration after relabeling is as the following:
            \begin{center}
                \Tableau[scale=0.7]{{3,0,1,2},{2,3,0},{1,2}}
            \end{center}
        As $\rho'=\rho+e+1-3=3$, the abacus configuration $\Ab(\lambda,\rho')$ is as the following:
            \begin{center}
                0\hspace{1mm}  1\hspace{1mm}  2\hspace{1mm}  3\hspace{1mm} \\
              \abacus(lmmr,bbbb,bbbb,nnbn,bnbn)
            \end{center}
        Notice that $N_0=0$ and we take $N=-2$. Add a runner to the leftmost and put $k(\lambda)+N_0-N+1=3$ beads from the top, we have the following new abacus $\Ab(\lambda_+,\rho_+')$:
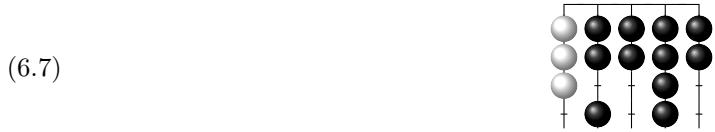
            \begin{equation}\label{ab:3}
                \abacus(lmmmr,obbbb,obbbb,onnbn,nbnbn)
            \end{equation}
        It is easy to get $\rho_+=(5,4,2)$. If we directly write the abacus $\Ab(\lambda,\rho)$:
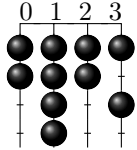
            \begin{center}
                0\hspace{1mm}  1\hspace{1mm}  2\hspace{1mm}  3\hspace{1mm} \\
                \abacus(lmmr,bbbb,bbbn,nbnb,nbnn)
            \end{center}
        Note this time we take $N_0=-1,N=-2$. Put a new runner between $2$-runner and $3$-runner with $k(\lambda)+N_0-N+1=2$ bead from the top, we have the following $\Ab(\lambda_+,\rho_+)$:
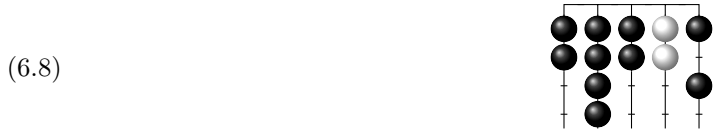
            \begin{equation}\label{ab:4}
                \abacus(lmmmr,bbbob,bbbon,nbnnb,nbnnn)
            \end{equation}
        It is clear \autoref{ab:4} is got from \autoref{ab:3} by moving beads to the left by $2$ (and add two essential beads from the northwest), hence denote the same partition $\lambda_+$.
    \end{enumerate}
\end{Example}
It is clear all above constructions generalize to $\ell$-partitions. We summarize the results in general case here:
\begin{Proposition}\label{equivdefgeneral}
    $\blam^+=\blam_+$.
\end{Proposition}
\begin{proof}
    This follows from the above arguments, \autoref{equivdefinition} and \autoref{addfromleft}.
\end{proof}
\begin{Proposition}\label{idempotentequalgeneral}
     With the same assumptions as in \autoref{equivdefgeneral},  $1^+_\blam=1_{\blam^+}$.
\end{Proposition}   
\begin{proof}
    This follows from the above arguments and \autoref{idempotentequal}.
\end{proof}
\begin{Theorem}\label{decompositionequalgeneral}
    With the same assumptions as in \autoref{equivdefgeneral}, and with regard to the charge $\brho$, the number of $i$-nodes in  $[\blam]$ is no more than $1$, then we have 
    \begin{equation}
        [\triangle(\blam):L(\bmu)]_q=[\triangle^+(\blam^+):L^+(\bmu^+)]_q
    \end{equation}
\end{Theorem}
\begin{proof}
    It follows from the above arguments and \autoref{decompositionequal}.
\end{proof}
\section{Further Topics}\label{furthertopics}
It is natural to think about the followings questions:
\begin{enumerate}[label=(\Alph*)]
    \item\label{othertypes} How about the cases in other (affine) types, like $C^{(1)}_{e}$, $A^{(2)}_{2e}$, $D^{(2)}_{e+1}$?
    \item\label{algorithm} How to make use of \autoref{decompositionequal} to develop an algorithm of computing decomposition numbers of KLRW algebras (in affine type A)?
    \item\label{comparison} There are empty runner removal theorem of Hekce (and Schur) algebras developed by Gordon and Mathas in \cite{GM-runnerremoval}, full runner removal theorem of Ariki-Koike algebras developed by Fayers, Dell'Arciprete in \cite{Fayers-runnerremoval, Dell-runner-removal}. Use the idempotent truncation from KLRW algebras to KLR algebras and then use the Brundan-Kleshchev isomorphism, can we give a comparison of the results in this paper with theirs?
    \item\label{conjectureques}To prove \autoref{conjecture} or weaken the condition of \autoref{decompositionequal}.
\end{enumerate}

We will not delve into these topics extensively in this paper, as doing so would make it overly lengthy. However, we provide some insights into these questions here.

For \ref{othertypes}, the construction of $\lambda^+$ is rather similar in these three types. The main obstruction is to deal with (seemingly) different cellular structure. For \ref{algorithm}, our aim is to start with a small quiver and gradually move to a larger quiver. However, directly generalizing \autoref{decompositionequal} is not feasible because the subdivided KLRW algebra has a ghost shift not equal to $1$, making it non-standard. Nevertheless, the canonical isomorphism can be applied to esolve this issue, and we plan to provide numerous examples. For \ref{comparison}, the task is more difficult. The idempotent truncation is outlined in \cite[Section 7]{Bo-many-cellular-structures}, however, to explicitly understand what happens (to simple/cell modules, decomposition nubmers, etc) in this process is already a very interesting and hard question, both in type $A^{(1)}_{e}$ and in other types. Additionally, translating the runner removal theorem using the Brundan-Kleshchev isomorphism is non-trivial. For \ref{conjectureques}, we need to fully analyze the structure of $J$ to prove it, or we can utilize \ref{comparison} to construct some counterexamples and conjecture the boundary condition.


\end{document}